\documentclass[onefignum,onetabnum]{siamart171218}

\usepackage{lipsum}
\usepackage{amsfonts}
\usepackage{graphicx}
\usepackage{epstopdf}
\usepackage{algorithmic}
\ifpdf
  \DeclareGraphicsExtensions{.eps,.pdf,.png,.jpg}
\else
\DeclareGraphicsExtensions{.eps}
\fi

\newsiamremark{remark}{Remark}
\newsiamremark{hypothesis}{Hypothesis}
\crefname{hypothesis}{Hypothesis}{Hypotheses}
\newsiamthm{claim}{Claim}
\newsiamthm{prop}{Proposition}

\headers{Layer Potential Methods for Doubly-Periodic Harmonic Functions}{B. Kim and B. Osting}

\title{Layer Potential Methods for \\ Doubly-Periodic Harmonic Functions\thanks{date: \today
\funding{B. Kim and B. Osting acknowledge partial support from NSF DMS-2136198 and DMS-2513175.}}}

\author{Bohyun Kim\thanks{Department of Mathematics, University of Utah, Salt Lake City, UT
  (\email{bohyun.kim@utah.edu}, \email{osting@math.utah.edu}).}
\and Braxton Osting\footnotemark[2]}

\usepackage{amsopn}

\makeatletter
\newcommand*{\addFileDependency}[1]{
  \typeout{(#1)}
  \@addtofilelist{#1}
  \IfFileExists{#1}{}{\typeout{No file #1.}}
}
\makeatother

\ifpdf
\hypersetup{
  pdftitle={Layer Potential Methods for Doubly-Periodic Harmonic Functions},
  pdfauthor={B. Kim and B. Osting}
}
\fi

\newtheorem{assumption}{Assumption}[section]
\newtheorem{rmk}{Remark}[section]
\newcommand{\onehalf}{\tfrac{1}{2}}
\usepackage{color}

\renewcommand{\S}{\mathrm{S}}
\newcommand{\K}{\mathrm{K}}
\newcommand{\M}{\mathrm{M}}
\newcommand{\X}{\mathrm{X}}

\usepackage{mathtools}
\usepackage{pgfplots}
\usepackage{bbm}

\begin{document}

\maketitle

\begin{abstract}
We develop and analyze layer potential methods to represent harmonic functions on finitely-connected tori (i.e., doubly-periodic harmonic functions).  The layer potentials are expressed in terms of a doubly-periodic and non-harmonic Green's function that can be explicitly written in terms of the Jacobi theta function or a modified Weierstrass sigma function. Extending results for finitely-connected Euclidean domains, we prove that the single- and double-layer potential operators are compact linear operators and derive the relevant limiting properties at the boundary. We show that when the boundary has more than one connected component, the Fredholm operator of the second kind associated with the double-layer potential operator has a non-trivial null space, which can be explicitly constructed. Finally, we apply our developed theory to obtain solutions to the Dirichlet and Neumann boundary value problems, as well as the Steklov eigenvalue problem. We present numerical results using Nystr\"om discretizations and find approximate solutions to these problems in several numerical examples. Our method avoids a lattice sum of the free-space Green's function, is shown to be spectrally convergent, and exhibits a faster convergence rate than the method of particular solutions for problems on tori with irregularly shaped holes. 
\end{abstract}



\section{Introduction}
\label{sec: introduction}
Doubly-periodic boundary value problems (BVPs) involving the Laplace, Helmholtz, Maxwell, and Stokes operators arise in a variety of physical and engineering problems;  examples include:
\begin{enumerate}
\item[(i)] Electrostatics in doubly-periodic composite materials~\cite{Greengard_1994,Helsing_1995,Otani_2006} and electromagnetic wave propagation in photonic crystals~\cite{Ammari_2018,Joannopoulos_2011, Kuchment_2016}
\item[(ii)] Thermal conductivity in doubly-periodic composite materials~\cite{kapanadze2014exact},  and  
\item[(iii)] Fluid flow through a doubly-periodic array of obstacles~\cite{emersleben1925darcysche,Greengard_2004,Hasimoto1959, Hasimoto2008, OGATA2003411, Pozrikidis_1996,Van_De_Vorst_1996}. 
\end{enumerate} 
Despite various approaches explored over the past century, the accurate approximation of solutions to such problems has remained a challenge, both analytically and computationally. 
Layer potential methods~\cite{Barnett_2018,folland2020introduction,greenbaum1993laplace,Greengard_2004,Helsing_1995,kress1989linear,MIKHLINdirichlet} are widely used for studying spatially-homogeneous elliptic BVPs due to their high accuracy~\cite{hao2011high} and versatility, as they can be extended to 
the Helmholtz equation~\cite{kress1991boundary}, 
the heat equation~\cite{folland2020introduction,kress1989linear}, 
Maxwell's equations~\cite{Ammari_2018}, 
the Laplace-Beltrami operator on the sphere~\cite{Gemmrich2008},
and more general elliptic
and parabolic problems~\cite{baderko1997parabolic,barton2017layer,Cioranescu1999}, as well as accommodate quasi-periodic boundary conditions~\cite{Ammari_2018,Barnett_2010}.

Typically, layer potential methods for solving harmonic BVPs on doubly-periodic domains rely on estimating a periodic Green's function $G_{\textrm{per}}(z)$ and representing the solution as a convolution of $G_{\textrm{per}}$ with an unknown boundary density $\phi$. The solution to the BVP is then obtained by solving a linear integral equation for $\phi$. Two common approaches to efficiently approximating solutions to the integral equation are lattice sums~\cite{Berman1994,borwein2013lattice,Cazeaux_2015} and Ewald-based methods~\cite{Ewald1921,LINDBO2011}. These methods are intuitive, but implementations with near linear time complexity, require (i) careful treatment of error terms arising from lattice summation and (ii) the implementation of numerical techniques such as the Fast Multipole Method \cite{GREENGARD1987325}. Other notable approaches include the Generalized Multipole Technique~\cite{hafner1990generalized}, a unified integral equation method~\cite{Barnett_2018}, and the method of particular solutions for doubly-periodic settings~\cite{kaoharmonic2023}. 

In our work, \textit{we are particularly interested in solving the Laplace equation on doubly-periodic, finitely-connected 2D domains (i.e., finitely-connected flat tori) with the Dirichlet, Neumann, or Steklov boundary conditions} by leveraging an explicit  Green's function representation~\cite{Bergweiler2016, Hasimoto2008, kaoharmonic2023,lin2010elliptic} based on Weierstrass functions. Using this representation, we develop an alternative layer potential approach for harmonic BVPs. The lattice sum in the proposed method is effectively incorporated into the evaluation of the Green's function, which can be efficiently evaluated using an exponentially convergent series. We validate our approach with rigorous proofs and numerical experiments, demonstrating its efficiency in solving the Dirichlet and Neumann BVPs, as well as Steklov eigenvalue problems (EVPs). In a variety of numerical experiments, we show that the method is spectrally convergent and can accurately compute solutions to machine precision using a small number of discretization points.

\subsubsection*{Geometry and problem setup}
\begin{figure}[t!]
\begin{center}
\includegraphics[width=0.48\textwidth]{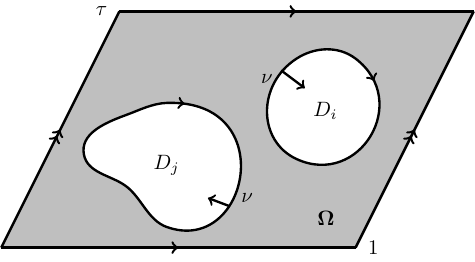}
\includegraphics[width=0.51\textwidth]{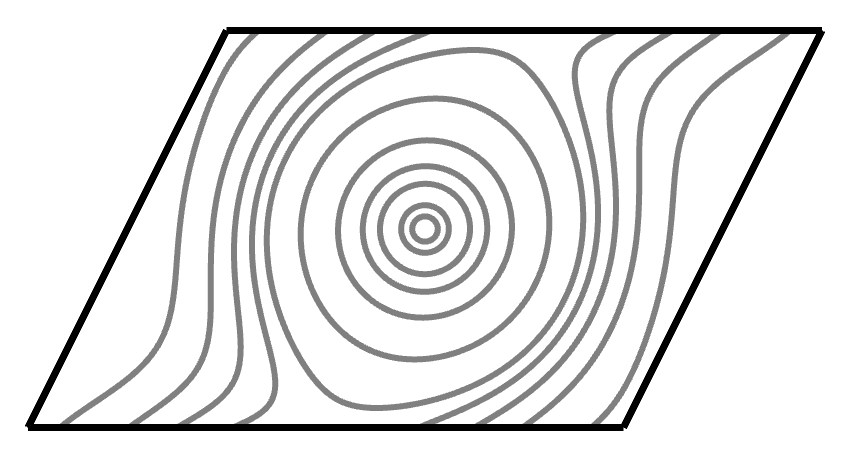}
\end{center}
\caption{\textbf{(left)} Illustration of a finitely-connected flat torus $\Omega$ as in~\eqref{e:Omega}.\textbf{~(right)} Level sets of the Green's function $z \mapsto G\left(z - \frac{1}{2} (1 + \tau) \right)$ from~\eqref{eqn: green} on the torus $\mathbb{T}_\tau$ with $\tau = \frac{1}{3} + \frac{2}{3}i$.}
\label{fig:fun domain}
\end{figure}

Consider a torus $\mathbb{T}_\tau = \mathbb{C}/\Lambda_\tau$, where $\Lambda_\tau = \mathbb{Z}+  \mathbb{Z}\tau$ is a lattice with periods $1 \text{ and } \tau = a + bi\, (b> 0)$, with area  $|\mathbb{T}_\tau| = b$. 
Typical examples include the \emph{square torus} ($\tau = i$) and the \emph{equilateral torus} ($\tau = \frac{1}{2} + \frac{\sqrt{3}}{2}i$). Let \(\Omega\) denote a \emph{finitely-connected torus} with $M$ holes:
\begin{equation}
\label{e:Omega}
    \Omega := \mathbb{T}_\tau \setminus \overline D, \qquad 
    D:= \bigcup_{j=1}^M D_j, 
\end{equation}
where $D_j$ are ``holes'' removed from the torus. Throughout, we assume the following. 
\begin{assumption} \label{assumption1}
The holes \(D_j \subset \mathbb C \), for \(j = 1, 2, \ldots, M\), are open, simply-connected sets with disjoint closures.  We assume that each hole \( D_j \) has a \( C^2 \) boundary, \( \partial D_j \). Let \( \nu(z) \) be the unit normal vector at \( z \in \partial \Omega \) pointing into \( D_j \), such that \( \partial \Omega = \cup_{j =1}^M \partial D_j = \partial D\) has a positive (counterclockwise) orientation with respect to $\Omega$. 
\end{assumption}
These definitions and assumptions are illustrated in Fig.~\ref{fig:fun domain}(left). Throughout, we identify $z = x+iy \in \mathbb C$ with $(x,y) \in \mathbb R^2$ where $x = \Re z$ and $y = \Im z$ denote the real and imaginary parts of $z$, respectively.

Our goal is to develop layer potential methods to solve the Dirichlet BVP, the Neumann BVP, and the Steklov EVP, respectively: 
\begin{align}
    \Delta u &= 0  \quad \text{in } \Omega,\qquad\quad
         u = g  \quad \text{on } \partial \Omega, 
    \label{eqn: dirichlet}\\
        \Delta u &= 0  \quad \text{in } \Omega,\qquad\quad
         \partial_\nu u  = g  \quad \text{on } \partial \Omega, 
    \label{eqn: neumann}\\
        \Delta u &= 0  \quad \text{in } \Omega,\qquad\quad
          \partial_\nu u  = \sigma u \quad \text{on } \partial \Omega.
    \label{eqn: steklov}
\end{align}
Here, 
$\Delta =  \partial^2_{x}+\partial^2_{y}$ 
is the Laplacian and $\partial_\nu$ denotes the normal derivative at the boundary. In~\eqref{eqn: dirichlet}, we assume $g \in C(\partial \Omega)$. In~\eqref{eqn: neumann}, we assume $g \in C_0(\partial \Omega):= \{ g \in C(\partial \Omega) \colon \int_{\partial \Omega} g(\xi)\, |d\xi|  = 0\}$. In~\eqref{eqn: steklov}, $\sigma\geq 0$ and $u \not\equiv 0$ are the Steklov eigenvalues and eigenfunctions, respectively. 
Problems \eqref{eqn: dirichlet}--\eqref{eqn: steklov} are  challenging  for $M \geq 2$ holes due to the combined effects of periodicity and multiple connectivity.

\begin{remark} \label{rem:MultConnDom}
Harmonic functions on multiply-connected domains often exhibit behavior absent in simply-connected ones. This can be described and interpreted from several different viewpoints:
(i) analytic functions and the logarithmic conjugation theorem~\cite{Axler_1986,kaoharmonic2023,MIKHLINdirichlet,Trefethen_2018},
(ii) flux in an ideal flow~\cite{kundu2012,Lamb_1895}, and 
(iii) the distinction between closed and exact forms~\cite{cohn2014,nahon2024}. 
Such behavior makes the development of numerical methods on multiply-connected domains more challenging; see, e.g.,  
\cite{greenbaum1993laplace,HELSING2008,MIKHLINdirichlet,Nasser_2011,Oudet_2021,Trefethen_2018,WEGMANN200836}. 
\end{remark} 
While several approaches consider multiply-connected regions in free space, extending such methods to doubly-periodic domains requires care due to the lack of periodicity in the standard logarithm~\cite{Barnett_2018,greenbaum1993laplace,Greengard_2004,Greengard_1994,Helsing_1995,kaoharmonic2023}. Instead, we associate a logarithm-like term with each hole via the Green’s function on a torus, which addresses both periodicity and multiple connectivity within a single, consistent framework.

\subsubsection*{Green's function on a torus}
The Green's function on $\mathbb{T}_\tau$ is defined by
\begin{equation}
    G(z) :=-\frac{1}{2\pi}\log|\vartheta_1(z)|+\frac{1}{2b}\Im(z)^2, 
    \qquad z \in \mathbb{T}_\tau,
    \label{eqn: green}
\end{equation}
where $ b = \Im \tau>0$~\cite{lin2010elliptic}. 
The function $G$ is doubly-periodic and satisfies 
\begin{align}
\Delta G(z) &=\left[\partial^2_{x}+\partial^2_{y}\right]G(z) 
= \frac{1}{b}-\delta(z), 
\qquad z \in \mathbb{T}_\tau,
\label{eqn: not harmonic G}
\end{align} 
where $\delta(z)$ is the Dirac delta.  Fig.~\ref{fig:fun domain}(right) illustrates the Green's function $G$; the properties of the  Jacobi theta function $\vartheta_1$ and $G$ are reviewed in Sec.~\ref{sec:LayerPotentialProperties}. 
Note in~\eqref{eqn: not harmonic G} that $G$ is not harmonic away from the origin; in fact, there is no doubly-periodic Green's function satisfying \(\Delta G(z) = -\delta(z)\), since integrating both sides over \(\mathbb{T}_\tau\) and applying Green's theorem yields a contradiction~\cite{Gemmrich2008}. The constant term \(\frac{1}{b}\) on the right-hand side of~\eqref{eqn: not harmonic G} ensures that both sides of the equation integrate to zero. Using the Green's function $G$, we define layer potentials to solve \eqref{eqn: dirichlet}--\eqref{eqn: steklov}.

\subsubsection*{Main results: layer potential methods} 
Layer potential methods express solutions to Laplace BVPs as the convolution of Green's function with an unknown boundary density \(\phi\). Excellent general references on layer potential methods can be found in~\cite{folland2020introduction,kress1989linear,MIKHLINdirichlet}. 
We define the single- and double-layer potentials on $\mathbb T_\tau$ using the doubly-periodic Green's function~\eqref{eqn: green}. Properties of these operators are established in Sec.~\ref{sec:LayerPotentialProperties}.

\clearpage
\begin{definition}
\label{def: layer pots}
Let $\Omega$ satisfy Assumption~\ref{assumption1}. 
For \(\phi \in C(\partial \Omega)\), we define the single-layer and double layer potentials by
\begin{align}
\mathcal{S}[\phi](z) &:= \int_{\partial \Omega}  G(z-\xi)\phi(\xi) \, |d\xi|, 
\qquad z \in \mathbb{T}_\tau \quad \text{and} \label{eqn: single layer}\\
\mathcal{D}[\phi](z) &:= \int_{\partial \Omega}  \partial_{\nu_\xi}G(z-\xi)\phi(\xi) \, |d\xi|, 
\qquad z \in \mathbb{T}_\tau \setminus \partial\Omega, \label{eqn: double layer}
\end{align} respectively. 
Here, $\nu$ is the unit normal vector pointing into \( D \) (see Fig.~\ref{fig:fun domain}(left)). Additionally, for \(\phi \in C(\partial \Omega)\), we define the \textit{modified single-layer} potential by
\begin{equation}
\mathcal{S}_0[\phi](z) := \mathcal{S}[\phi -\M[\phi]](z) + \overline{\phi}, 
\qquad z \in \mathbb{T}_\tau, 
\label{eqn: single-layer Kress formulation}
\end{equation}
where 
\begin{equation}\M[\phi] \coloneq \overline{\phi} = \frac{1}{|\partial \Omega|} \int_{\partial \Omega} \phi(\xi) \, |d\xi| \label{eqn: average operator M}
\end{equation}
is the mean value of $\phi$. 
\end{definition}

\begin{definition} 
\label{def: boundary double} 
Let $\Omega$ satisfy Assumption~\ref{assumption1}. Define the boundary operator $\K \colon C(\partial\Omega) \to C(\partial\Omega)$ and its adjoint $\K^{*} \colon C(\partial\Omega) \to C(\partial\Omega)$ as
 \begin{align}
     &\K[\phi](z_0) \coloneq \mathrm{p.v.}\int_{\partial \Omega} \partial_{\nu_\xi}G(z_0-\xi)\phi(\xi)  \, |d\xi|, 
     \qquad z_0 \in \partial \Omega, 
     \label{eqn: boundary double-layer}\\
     &\K^*[\phi](z_0) \coloneq\mathrm{p.v.}\int_{\partial \Omega} \partial_{\nu_z}G(z_0-\xi)\phi(\xi)  \, |d\xi|, 
     \qquad z_0 \in \partial \Omega, 
    \label{eqn: adjoint double-layer}
 \end{align}
where $\mathrm{p.v.}$ denotes  the Cauchy principal value.  Let $\S \colon C_0(\partial \Omega) \to C_0(\partial \Omega)$ and \\ $\S_0 \colon C(\partial \Omega) \to C^1(\partial \Omega)$ 
denote the continuous restrictions of $\mathcal{S}$ and $\mathcal{S}_0$ to the boundary, respectively. 
\end{definition}

The following three theorems are the main results of this paper, providing solutions to the Dirichlet BVP~\eqref{eqn: dirichlet}, Neumann BVP~\eqref{eqn: neumann}, and Steklov EVP~\eqref{eqn: steklov} via layer potentials; proofs are provided in Sec.~\ref{sec: BIE}.

\begin{theorem}{(Solution to the Dirichlet BVP)}
\label{thm: dirichlet}
Let $\Omega$ satisfy Assumption~\ref{assumption1},
$g \in C(\partial \Omega)$, and 
fix points $\beta_j \in D_j$. Then there exist unique \(\phi \in C(\partial \Omega)\) and $A_j \in \mathbb{R}$ for $j=1,2,\ldots, M$ that satisfy the following system:
\begin{subequations}
\label{e:DirIntEq}
\begin{align}\label{e:DirIntEqA}
    & -\onehalf  \phi(z_0) +\K[\phi](z_0) + \sum_{j =1}^M A_j G(z_0 - \beta_j) = g(z_0), \quad  \, z_0 \in \partial \Omega, \\
     \label{e:DirIntEqB}& \int_{\partial D_j} \phi(z_0) \, |dz_0| = 0,  \quad \forall \, j =1,2,\ldots, M - 1, 
    \quad \textrm{and} \quad 
    \sum_{j=1}^M A_j = 0.
\end{align}
\end{subequations}
The solution to the Dirichlet BVP~\eqref{eqn: dirichlet} is represented as
\begin{equation}
    u(z) = \mathcal{D}[\phi](z) + \sum_{j=1}^M A_j G(z - \beta_j), \quad \quad z\in \Omega. \label{e:uRep}
\end{equation}
Here, the quantity $A_j$ can be interpreted as the flux across $\partial D_j$.  
\end{theorem}
\begin{rmk}
In particular, for \(M = 1\), the condition $\sum A_j = 0$ implies  \(A_1 = 0\). In this case,~\eqref{e:DirIntEq} reduces to 
\begin{equation}
-\onehalf \phi(z_0) + \K[\phi](z_0)
= g(z_0), 
\qquad z_0 \in \partial \Omega
\label{eqn: DirIntEQM1}
\end{equation}
and $u(z) = \mathcal{D}[\phi](z) $ is the solution to the Dirichlet BVP~\eqref{eqn: dirichlet}. 
\label{rmk: m1 dirichlet}
\end{rmk}

\begin{theorem}{(Solution to the Neumann BVP)}
\label{thm: neumann}
Let $\Omega$ satisfy Assumption~\ref{assumption1} and \(g \in C_{0}(\partial \Omega)\). 
There exists a unique $\phi \in C_0(\partial \Omega)$ that satisfies  
\begin{align}
\label{eqn: neumann bdry eqn}
&\onehalf\phi(z_0) +\K^*[\phi](z_0) = g(z_0),
 \qquad z_0\in \partial\Omega. 
\end{align}
Here, a solution to the Neumann BVP~\eqref{eqn: neumann} is represented as
\begin{align}
\label{e:uRep, neu}
&u(z) = \mathcal{S}[\phi](z)+ C, 
\qquad  z\in \Omega 
\end{align}
for any $C \in \mathbb R$. The quantity $\int_{\partial D_j} g(\xi)\, |d\xi| = \int_{\partial D_j}\phi(\xi)\, |d\xi| $
can be interpreted as the flux across $\partial D_j$,  for $j=1,2,\ldots,M$.  
\end{theorem}

Note that in Theorem~\ref{thm: dirichlet}, the solution to the Dirichlet BVP~\eqref{eqn: dirichlet} is represented in~\eqref{e:uRep} as the sum of a double-layer potential and a linear combination of Green's functions, while in Theorem~\ref{thm: neumann} the solution to the Neumann BVP~\eqref{eqn: neumann}  is represented in~\eqref{e:uRep, neu} as a single-layer potential alone. 
The inclusion of the sum in~\eqref{e:uRep} can mathematically be understood from the viewpoints discussed in Remark~\ref{rem:MultConnDom}.  Physically,  in electrostatics, a single-layer potential represents the potential induced by a charge $\phi$ distributed on $\partial\Omega$, which can generate a net charge. In contrast, the double-layer potential represents the potential induced by a dipole distribution on $\partial\Omega$, which does not generate net charge; thus, it requires point charges introduced by $A_j$ \cite{folland2020introduction}.

\smallskip

Finally, we consider the Steklov EVP~\eqref{eqn: steklov};  recent surveys on Steklov eigenvalues can be found in  \cite{girouard2017spectral,kuznetsov2014legacy}. 
The Steklov spectrum is discrete; we enumerate the eigenvalues, counting multiplicity, as:
$0 = \sigma_1  <  \sigma_2  \leq  \sigma_3  \leq \ldots  \to \infty $. 
The eigenspace corresponding to $\sigma_1 =0$ is spanned by the constant function, and the restriction of the Steklov eigenfunctions to the boundary, 
$\{ u_j |_{\partial \Omega} \}_{j=1}^\infty \subset C^\infty(\partial \Omega)$, 
forms a complete orthonormal basis for $L^2(\partial \Omega)$. In addition, the Steklov spectrum coincides with the spectrum of the Dirichlet-to-Neumann operator $\Gamma \colon H^{\frac1 2}(\partial \Omega) \to H^{- \frac 1 2}(\partial \Omega)$, defined as $\Gamma w =  \partial_\nu ( \mathcal H w)$, where $\mathcal H w$ denotes the unique harmonic extension of $ w \in H^{\frac1 2}(\partial \Omega)$ to $\Omega$.

\begin{theorem}{(Solution to the Steklov EVP)}
\label{thm: steklov}
Let $\Omega$ satisfy Assumption~\ref{assumption1}. 
The eigenpair $(\sigma_k, u_k)$  is a solution to the Steklov EVP~\eqref{eqn: steklov} if and only if  $(\sigma_k, \phi_k)$ solves
\begin{align}
     (\K^* + \onehalf I )(I-\M)[\phi_k](z_0) = \sigma_k \S_0[\phi_k](z_0),\qquad z_0\in \partial\Omega,
     \label{e: steklov bie}
\end{align} 
where  $\M$ is defined in~\eqref{eqn: average operator M}. The eigenfunction $u_k$ is represented as 
\begin{equation}
\label{e:uRep, stek}
u_k(z) = \mathcal{S}_0[\phi_k](z), 
\qquad  z\in \Omega. 
\end{equation}
The flux across $\partial D_j$ for $j= 1,2,\ldots,M$ is given by $\int_{\partial D_j} \partial_\nu u_k(\xi)\, |d\xi| = \int_{\partial D_j} (\K^{*}+\onehalf I )(I-\M) [\phi_k] (\xi)\, |d\xi| $. 
\end{theorem}

Theorems~\ref{thm: dirichlet},~\ref{thm: neumann}, and~\ref{thm: steklov} serve as the basis for numerically approximating solutions to the Dirichlet BVP, Neumann BVP, and Steklov EVP. We outline the structure of the paper as follows. In Sec.~\ref{sec:LayerPotentialProperties}, we review the properties of the doubly-periodic Green's function and establish the fundamental properties of the associated layer potentials. In Sec.~\ref{sec: BIE}, we provide proofs for Theorems~\ref{thm: dirichlet}--\ref{thm: steklov}. Proofs of the lemmas in Sec.~\ref{sec:LayerPotentialProperties}--\ref{sec: BIE} are given in Appendix~\ref{sec: proofs of lem thms-properties}.
In Sec.~\ref{sec: Numerical Experiments}, we present the results of six numerical experiments that demonstrate the effectiveness of our approach for each of these problems. To numerically approximate the boundary integrals, we implement Nystr\"om method in MATLAB, which is available on GitHub~\cite{bohyunBIEgithub}. Finally, we conclude in Sec.~\ref{s:disc} with a brief discussion.

\subsection{Notation}
\label{sec: notation}
\begin{itemize}
    \item For \(z_1, z_2\in \mathbb{C}\), define $z_1 \circ z_2 := \Re(z_1 \overline{z_2})$.
    \item \(\int_\Gamma f(z) \, |dz|\) denotes the integral of \(f(z)\) with respect to the arc length along the curve \(\Gamma\). 
    For a bounded region \(R\), we write \(\int_R f(z) \, dA\) to denote the double integral, where \(dA = dx \, dy\). 

    \item For a function $u$ with a jump continuity at $\partial \Omega$, we denote the limiting values of \(u\) approaching \(z_0 \in \partial\Omega\) from \(\Omega\) and \(D\) by 
    \begin{align*}
        u(z_0^{+}) = \lim_{\Omega \ni z \to z_0} u(z) 
        \qquad \textrm{and} \qquad 
       u(z_0^{-}) = \lim_{D \ni z \to z_0} u(z).
    \end{align*}
    Similarly, for a function $u$ with a jump continuity in its gradient at $\partial \Omega$, 
        \begin{align*}
        \partial_\nu u(z_0^{+}) =
        \lim_{\Omega \ni z \to z_0} \nu(z_0) \circ \nabla u(z)
        \quad \textrm{and} \quad 
        \partial_\nu u(z_0^{-}) = \lim_{D \ni z \to z_0} \nu(z_0) \circ \nabla u(z). 
    \end{align*}
    \item For $j=1,2, \ldots, M$, let $\mathbbm{1}_j(z)= \begin{cases}
1 & \text{if } z \in \partial D_j, \\
0 & \text{if } z \in \cup_{i\neq j} \partial D_i,
\end{cases}$ be the indicator function on $\partial D_j$,  let  $\mathbbm{1}$ be the the indicator function on $\partial\Omega$, and let $\mathbbm{1}_{D_j}(z)= \begin{cases}
1 & \text{if } z \in D_j, \\
0 & \text{if } z\in \cup_{i\neq j} D_i
\end{cases}$ be the indicator function on $D_j$.
\end{itemize}

\section{Properties of Doubly-Periodic Layer Potentials}
\label{sec:LayerPotentialProperties} 
We first review several properties of the doubly-periodic Green's function $G(z)$ defined in~\eqref{eqn: green}.  Excellent references include \cite{Bergweiler2016,lin2010elliptic} as well as \cite{chandrasekharan2012elliptic}. The definition of  $G(z)$ involves the Jacobi theta function as defined in \cite[Eq. (7.1)]{lin2010elliptic}, 
\begin{equation}
    \vartheta_1 (z) \coloneq  \vartheta(\pi z; \tau) 
    = -i \sum_{n\in \mathbb Z } (-1)^n q^{(n+\frac{1}{2})^2}e^{(2n+1)\pi i z}   
    \label{eqn: jacobi theta}
\end{equation}
with the nome $q = e^{\pi i \tau}$. The Jacobi theta function $\vartheta_1 (z)$ is an entire function  and the series representation \eqref{eqn: jacobi theta} is exponentially convergent. Moreover, $G(z)$ is doubly-periodic with periods 1 and $\tau$, since $\vartheta_1 (z)$ satisfies $\vartheta_1(z+1) = -\vartheta_1(z)$ and $\vartheta_1(z+\tau) = -q^{-1}e^{-2\pi i z}\vartheta_1(z)$. To see the $\tau$ periodicity,  for $b= \Im\tau$ and  $y = \Im z$, we compute 
\[ G(z+\tau) - G(z) = -\frac{1}{2\pi}\log e^{\pi b+2\pi y} + \frac{1}{2b}((y+b)^2-y^2) = 0.\] Furthermore, $G(-z) = G(z)$ since $\vartheta_1 (z)$ is an odd function. Finally, $G(z)$ has singularities only at the lattice points since $\vartheta_1 (z)$ only has simple zeros occurring at the lattice points.

Observe that the doubly-periodic single~\eqref{eqn: single layer} and double-layer~\eqref{eqn: double layer} potentials utilize $G(z-\xi)$ and $\partial_{\nu_\xi} G(z-\xi)$.
As $z\to \xi$, their asymptotic expansions satisfy (as proven in Appendix A) 
\begin{align}
    &G(z - \xi) = -\frac{1}{2\pi}\log|z - \xi| + \frac{ \log|\vartheta_1'(0)|}{2\pi} + O(|z-\xi|^2) \qquad \textrm{and} 
    \label{eqn: asymptotic G}\\ 
    &\partial_{\nu_\xi} G(z - \xi) = \frac{(z - \xi) 
 \circ \nu(\xi)}{2\pi |z - \xi|^2} + O(|z - \xi|). 
 \label{eqn: asymptotic Gnu}
\end{align}
Thus, standard analysis techniques for the free-space Green's function, \(G_{\textrm{free}}(z) = -\frac{1}{2\pi}\log|z|\), can be used to analyze the doubly-periodic single- and double-layer potentials. In Lemmas~\ref{lem: exact double layer}--\ref{lem: properties single layer}, we provide properties that are analogous to those of \(G_{\textrm{free}}(z)\), with proofs provided in Appendix~\ref{sec: proofs of lem thms-properties}.

\begin{lemma}[Gauss' Lemma for finitely-connected tori]
\label{lem: exact double layer}
Let $\Omega$ satisfy Assumption~\ref{assumption1}. Let $|D_j|$ denote the area of the hole $D_j$ for each $j = 1,2,\ldots, M$. Then,
\begin{equation}
 \int_{\partial D_j} \partial_{\nu_\xi} G(z-\xi)\, |d\xi| =  
    \begin{dcases}
      -\frac{|D_j|}{b} & \text{if } z \in \mathbb{T}_{\tau} \setminus \overline{D_j}, \\
      1 - \frac{|D_j|}{b} & \text{if } z \in D_j, \\ 
      \frac{1}{2} - \frac{|D_j|}{b} & \text{if } z \in \partial D_j,
    \end{dcases}
    \label{eqn: gen gauss sublemma}
\end{equation}
and, denoting $|D| =  \sum_{j}|D_j|$, 
\begin{equation}
   \int_{\partial \Omega} \partial_{\nu_\xi}  G(z-\xi)\, |d\xi| =  
    \begin{dcases}
      -\frac{|D|}{b} & \text{if } z \in \Omega,\\
      1 - \frac{|D|}{b} & \text{if } z \in D, \\ 
      \frac{1}{2} - \frac{|D|}{b} & \text{if } z \in \partial \Omega.
    \end{dcases}
    \label{eqn: gen gauss}
\end{equation} When $z\in \partial D_j$ or $z\in \partial\Omega$, the integrals are interpreted in the $\mathrm{p.v.}$ sense. 
\end{lemma}
Lemma~\ref{lem: exact double layer} is analogous to the classical Gauss' Lemma for $G_{\textrm{free}}(x) = -\frac{1}{2\pi} \log|x|$ (see \cite[Ex.~6.17]{kress1989linear} or \cite[Prop. 3.19]{folland2020introduction}). It differs only by the term $|D|/b$, which arises from the periodicity of $G$~\eqref{eqn: not harmonic G}. Gauss' Lemma is used to prove the following properties of double-layer potentials.

\begin{lemma}[Properties of Double-Layer Potentials]
\label{lem: properties double layer}
Let $\Omega$ satisfy Assumption~\ref{assumption1}. 
The double-layer potential $\mathcal{D}[\phi]$ in~\eqref{eqn: double layer} satisfies the following properties:
\begin{enumerate}
    \item  For $\phi \in C(\partial \Omega)$, $\mathcal{D}[\phi]$ is doubly-periodic, is smooth on $\mathbb{T}_\tau\setminus  \partial\Omega$, and satisfies $\Delta \mathcal{D}[\phi](z) = 0$ for $z \in \mathbb{T}_\tau\setminus  \partial\Omega$. 
 
    \item For $\phi \in C(\partial \Omega)$, the double layer potential satisfies the ``jump relations'',
    \[
    \mathcal{D}[\phi](z_0^{\pm}) = \K[\phi](z_0) \mp \onehalf\phi(z_0), \qquad z_0 \in \partial \Omega,
    \]
    where $\K$ is defined in~\eqref{eqn: boundary double-layer},
    \item The kernels of $\K$ and $\K^*$ in Definition \ref{def: boundary double} satisfy
\[
\lim_{z\to \xi} \ \partial_{\nu_\xi}G(z-\xi) =  -\frac{\kappa(\xi)}{4\pi} 
\quad \textrm{and} \quad 
\lim_{z\to \xi} \  \partial_{\nu_z}G(z-\xi) = -\frac{\kappa(\xi)}{4\pi}, 
\]
where $\kappa$ is the signed curvature of $\partial \Omega$. Since these kernels can be continuously extended to $z = \xi \in \partial \Omega $, $\K$ and $\K^*$ are compact linear operators. 
\item The null space of $\K - \onehalf I $, denoted by $N(\K -\onehalf I )$,  has dimension $M-1$ and has basis $\{ \phi_{jM}\}_{j=1}^{M-1}$, where $b = |\mathbb{T}_{\tau}|$ and 
\begin{align}
\phi_{jM}(z):= \frac{b}{|D_j|}\mathbbm{1}_j(z) - \frac{b}{|D_M|}\mathbbm{1}_M(z), \qquad   z\in\partial\Omega.
\label{eqn: characteristic phi}
\end{align}
Here, $\mathbbm{1}_j(z)$ is the indicator function on $\partial D_j$ (see Sec.~\ref{sec: notation}).
\item Similarly, $N(\K^{*}-\onehalf I )$ has dimension $M-1$. If  $\{ \psi_{k} \}^{M-1}_{k=1}$ is a basis for  $N(\K^*-\onehalf I )$, then each $\psi_{k}$ satisfies, for some constants $s_{kj}\in\mathbb{R}$:
\begin{subequations}
\label{eqn: Spsi span}
\begin{align}
\label{eqn: Spsi span A}
    & \mathcal{S}[\psi_{k}](z) =  \sum_{j=1}^M s_{kj} \mathbbm{1}_{D_{j}}(z), \quad  z\in D, \quad \text{and} \\
    \label{eqn: Spsi span B}
    & \S[\psi_{k}](z^{-}_0) =  s_{kj}, \quad z_0\in \partial D_{j}. 
\end{align}
\end{subequations} 
Furthermore, letting $\mathbf{s}_{k} = [s_{k1}, s_{k2}, \ldots, s_{kM}]^{T} \in \mathbb R^M$, $\{\mathbf{s}_1,\mathbf{s}_2, \ldots, \mathbf{s}_{M-1}, \mathbf{1}_M \}$ forms a basis for $\mathbb{R}^M$. 
\item $\dim N(\K^*+\onehalf I ) =\dim N(\K+\onehalf I )= 0$. 
    \end{enumerate}
\end{lemma}

We note that the limit of the double-layer potential given in Lemma~\ref{lem: properties double layer}(2) matches the limit of the classical double-layer potential (see \cite[Thm. 6.18]{kress1989linear}), with the difference in sign arising from our choice of the normal direction (see Fig.~\ref{fig:fun domain}(left)). 

\begin{figure}[t!]
\begin{center}
\includegraphics[height=0.25\textwidth]{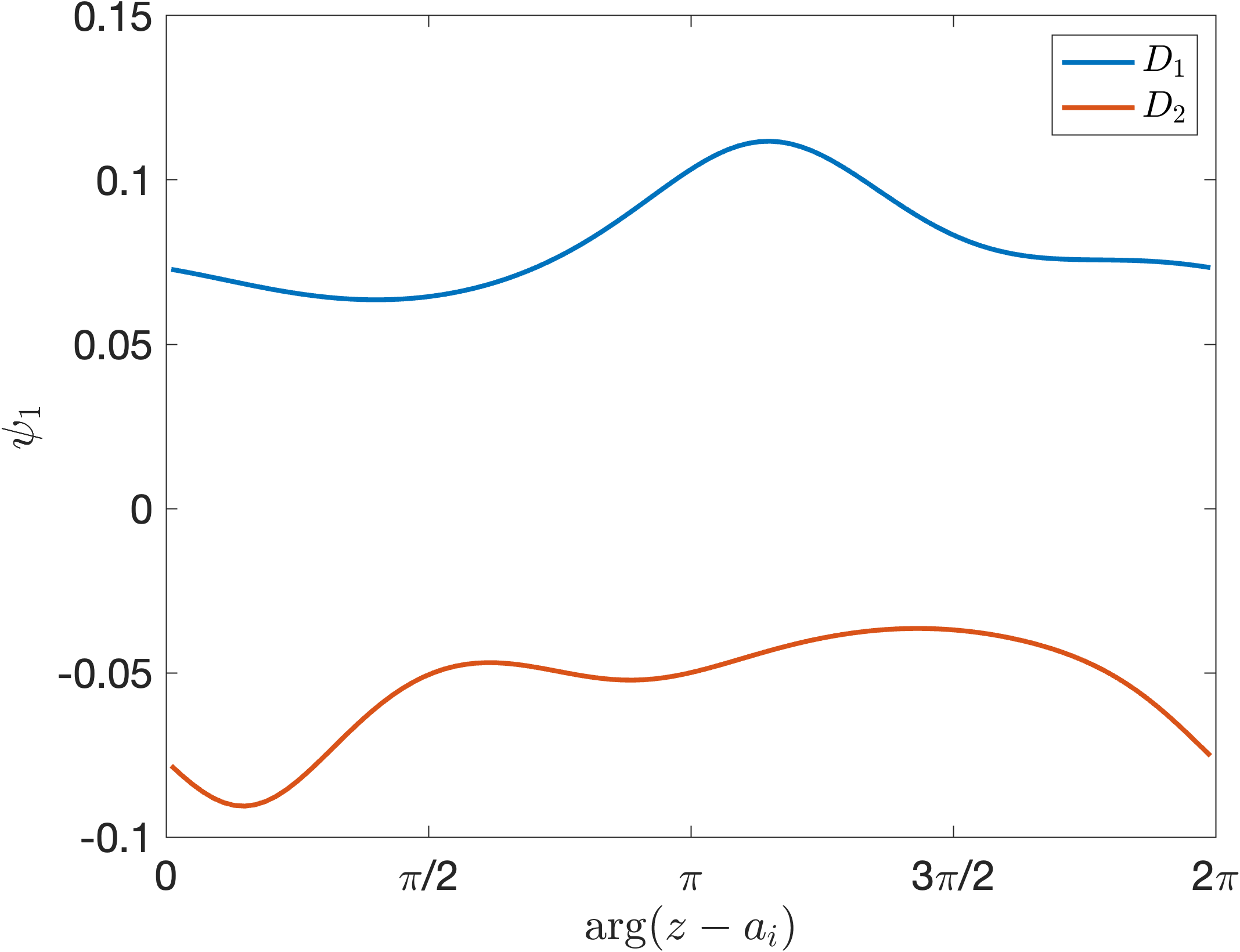}
\qquad 
\includegraphics[height=0.25\textwidth]{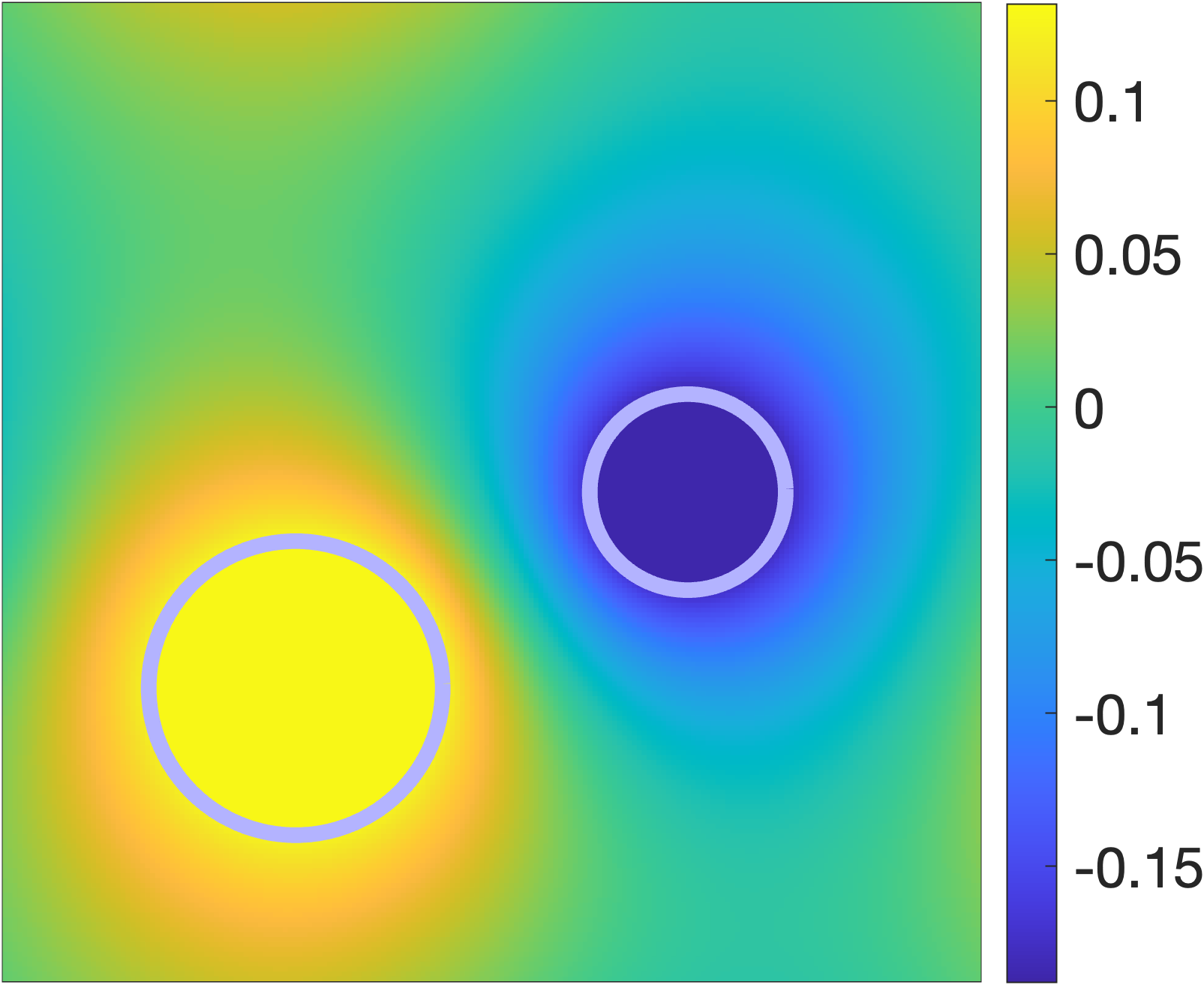}
\end{center}
\caption{For a square torus with $M=2$ circular holes, we plot 
$\psi_1 \in N(\K^{*}-\onehalf I )$ {\bf~(left)} 
and $\mathcal{S}[\psi_1](z)$, $z \in \mathbb T_{\tau}$ {\bf (right)}.  
As proven in Lemma~\ref{lem: properties double layer} (5), 
$\mathcal{S}[\psi_1]$ is constant on each hole; see  Remark~\ref{r:NullSpace}.}
\label{fig: Dstar kernel}
\end{figure}

\begin{rmk}
\label{r:NullSpace}
In Lemma~\ref{lem: properties double layer}(5), a  basis $\{\psi_k\}_{k =1}^{M-1}$ of $N(\K^*-\onehalf I )$ corresponds to functions $\mathcal{S}[\psi_{k}]$ that are constant on each hole $D_j$. In Fig.~\ref{fig: Dstar kernel}, for a square torus $\mathbb T_{\tau}$ ($\tau=i$) with $M=2$ circular holes, we plot $\psi_1$ on $\partial \Omega$~(left) and  $\mathcal{S}[\psi_{1}](z)$, $z \in \mathbb T_{\tau}$~(right). 
In this case, $\mathbf{s}_1 = [-0.1856, 0.1290]^{T}$. The holes have centers $a_1 = 0.7+0.5i$ and $a_2 =0.3+ 0.3i$ with radii $r_1 = 0.1$ and $r_2 = 0.15$. 
See Sec.~\ref{sec: Numerical Experiments} for a discussion on the numerical methods used to generate this figure. 
\end{rmk}

\begin{lemma}[Properties of Single-Layer Potential]
\label{lem: properties single layer}
Let $\Omega$ satisfy Assumption~\ref{assumption1}. The single-layer potential $\mathcal{S}[\phi]$ given in~\eqref{eqn: single layer} satisfies the following properties:  
\begin{enumerate}
\item For  $\phi\in C_0(\partial \Omega)$,  $\mathcal{S}[\phi]$ is doubly-periodic, continuous on $\mathbb{T}_\tau$,  smooth on $\mathbb{T}_\tau \setminus \partial \Omega$, and satisfies  $\Delta \mathcal{S}[\phi](z) = 0$ for $z \in \mathbb{T}_\tau\setminus\partial\Omega$.

\item For $\phi \in C_0(\partial\Omega)$,the normal derivative of the single layer potential satisfies the ``jump relations'', 
    \begin{align*}
    \partial_\nu \mathcal{S}[\phi](z^{\pm}_0) =  \K^{*}[\phi](z_0)\pm \onehalf\phi(z_0), 
    \qquad z_0 \in \partial \Omega,
    \end{align*}
    where $\K^*$ is defined in~\eqref{eqn: adjoint double-layer}.
    \item \(\S\colon C_0(\partial \Omega) \to C_0(\partial \Omega) \) is a compact linear operator. 
\end{enumerate}
\end{lemma}

\section{Layer Potential Methods}
\label{sec: BIE}
In this section, we prove 
Theorems~\ref{thm: dirichlet}--\ref{thm: steklov}, 
which establish layer potential representations for the solutions to the Dirichlet BVP \eqref{eqn: dirichlet}, Neumann BVP~\eqref{eqn: neumann}, and Steklov EVP~\eqref{eqn: steklov}. Throughout the proofs, we employ the Fredholm alternative~\cite[Appendix D]{evanspartial} and adapt arguments for multiply-connected Euclidean domains: for the Dirichlet BVP, we follow~\cite[Sec. 29]{MIKHLINdirichlet} and \cite[Sec. 3.1]{greenbaum1993laplace}; for the Neumann BVP, we refer to~\cite[Sec. 37, Lemma 2]{mikhlin2020linear}; and for the Steklov EVP, we utilize~\cite[Theorem 7.36]{kress1989linear}.

\subsection{The Dirichlet BVP  and Proof of Theorem~\ref{thm: dirichlet}}
\label{sec: dirichlet proof}
Lemma~\ref{lem: properties double layer}(4) shows that the null space of $\K - \onehalf I $ is trivial for $M = 1$ and nontrivial for $M \geq 2$. Accordingly, we introduce the following integral operator to handle both cases simultaneously. 
\begin{definition}[Characteristic operator]
\label{def: characteristic}
Define the characteristic operator as the linear operator  $\X \colon C(\partial \Omega) \to C(\partial \Omega)$ for $M=1$ by  $\X[\phi] \equiv 0$, and for $M\geq2$ by
\begin{align}
&\X[\phi](z_0) := \int_{\partial\Omega} \chi(z_0,\xi) \phi(\xi) \,|d\xi|, 
\qquad 
z_0\in\partial\Omega, 
\label{eqn: characteristic operator}
\end{align}
where 
$\chi(z_0,\xi)\coloneq \sum_{j=1}^{M-1}\mathbbm{1}_j(z_0) \mathbbm{1}_j(\xi)$. 
\end{definition}
The operator $\X$ integrates $\phi$ only over the boundary component containing $z_0$, excluding the $M$th component. The following lemma establishes preliminary properties of the characteristic operator. 
\begin{lemma}[Properties of the Characteristic Operator]
\label{lem: characteristic}
Let \(\Omega\) satisfy Assumption~\ref{assumption1}. 
The characteristic operator in \eqref{eqn: characteristic operator} satisfies the following properties: 
\begin{enumerate}
\item  $\X$ is a compact linear operator. For any $\phi \in C(\partial \Omega)$, \(\X[\phi] \)   is a constant function on each $\partial D_j$,  for $j=1,2,\ldots, M$, and vanishes on $\partial D_M$.
\item $\K + \X - \onehalf I \colon C(\partial \Omega) \to C(\partial \Omega)$ is an injective compact linear operator; consequently, by the Fredholm alternative, it is also surjective. 
\end{enumerate}
\end{lemma}
Lemma~\ref{lem: characteristic} is used in the proof of Theorem~\ref{thm: dirichlet} and is proven in Appendix~\ref{sec: proofs of lem thms-properties}. The injectivity of $\K+\X -\onehalf I$ is established using the properties of the complex double-layer potential in Lemma~\ref{lem: Complex layer}.

\begin{proof}[Proof of Theorem~\ref{thm: dirichlet}] 
Our goal is to express system~\eqref{e:DirIntEq} using $\K + \X - \onehalf I$. By applying Lemma~\ref{lem: characteristic}(2), we establish the existence and uniqueness of the system.

\medskip 
\noindent \underline{$M=1$.}
Since $\X\equiv 0$ by Definition~\ref{def: characteristic}, \eqref{e:DirIntEq} reduces to $\left(\K - \onehalf I\right)[\phi] = g$ (see Remark~\ref{rmk: m1 dirichlet}). Therefore, for any $g\in C(\partial \Omega)$, there exists a unique $\phi\in C(\partial \Omega)$ satisfying~\eqref{e:DirIntEq}. This proves the theorem for the case $M=1$. 

\medskip 
\noindent \underline{$M\geq 2$.} 
Assume $g\in C(\partial\Omega)$ and fix $\beta_j \in D_j$. We show that there exist unique $\phi_A$ and $A_j \in \mathbb{R}$ for $j=1,2,\ldots,M$ that solve the system~\eqref{e:DirIntEq}. By Lemma~\ref{lem: characteristic}(2), there exists a unique $\phi^j\in C(\partial\Omega)$ corresponding to $g^j\in C(\partial\Omega)$ for $j=0,1,\ldots,M$, solving 
\begin{align}
&(\K+\X-\onehalf I )[\phi^j](z_0) = g^j(z_0), \qquad z_0\in \partial\Omega, \quad \text{where}
    \label{eqn: phij gj}\\
\nonumber    &g^j(z_0)  = 
\begin{cases}
g(z_0) & \text{if}\quad j =0,\\
    G(z_0-\beta_j) & \text{if} \quad j=1,2,\ldots, M.
\end{cases} 
\end{align}

First, we show  that there exists a unique set $\{A_j\}_{j=1}^M$ satisfying
\begin{equation} \label{e:AkSystem}
\sum_{j=1}^M A_j \X[\phi^j](z_0) = \X[\phi^0](z_0), \qquad   z_0\in\partial\Omega 
\qquad \text{and} \qquad 
\sum_{j=1}^M A_j = 0.
\end{equation}Since $\X[\phi^j]$ for $j= 0, 1,\ldots, M$ are piecewise constant on each $\partial D_k$, and vanish on $\partial D_M$, we may view~\eqref{e:AkSystem} as an $M\times M$ linear system. Showing that this system is nonsingular implies the uniqueness of $A_j$. To show it is nonsingular, assume $\sum_{j=1}^{M}A_j\X[\phi^j](z_0) = 0$ on $ \partial\Omega$ and $\sum_{j=1}^M A_j = 0$. Substituting the first condition into~\eqref{eqn: phij gj} and applying the Fredholm alternative, we obtain 
\begin{align*}
    &\sum_{j=1}^{M}A_j\,G(z_0-\beta_j) 
    \ \in \ 
    R(\K-\onehalf I) 
    \ = \ 
    N(\K^*-\onehalf I )^\perp.
\end{align*}
Let $\{\psi_k\}_{k=1}^{M-1}$ be the basis for $N(\K^*-\onehalf I)$ as in Lemma~\ref{lem: properties double layer}(5). By~\eqref{eqn: Spsi span},
\begin{align*}
    0 &= \left\langle \sum_{j=1}^M A_j\,G( \cdot-\beta_j), \psi_{k} \right\rangle 
 = \sum_{j=1}^M A_j\mathcal{S}[\psi_{k}](\beta_j)= \sum_{j=1}^{M} A_j s_{kj} = \begin{bmatrix}A_1 \, A_2 \,\ldots \, A_M\end{bmatrix} \mathbf{s_k}, 
\end{align*} for each $k=1,2,\ldots,M-1$. In addition to these $M-1$ conditions, $\sum_{j=1}^M A_j = 0$ in \eqref{e:AkSystem} provides $[A_1, \ldots, A_M] \mathbf{1}_M  = 0$,
where $\mathbf{1}_M \in \mathbb{R}^M$ is the vector of all ones.  Finally, since $\{ \mathbf{s}_1,\mathbf{s}_2, \ldots, \mathbf{s}_{M-1}, \mathbf{1}_M \}$ forms a basis for $\mathbb{R}^M$ (Lemma~\ref{lem: properties double layer}(5)), we conclude that $A_j=0$ for all $ j=1,2,\ldots,M$.

Next, let
\begin{align}
    &\phi_A(z_0) := \phi^0(z_0) - \sum_{j=1}^M A_j \phi^j(z_0), \qquad  z_0 \in \partial\Omega,\\
     &u_A(z) := \mathcal{D}[\phi_A](z) + \sum_{j=1}^{M} A_j G(z-\beta_j), \qquad  z \in \Omega.  
\end{align}
Recall that the $\phi^j$ are already uniquely determined by~\eqref{eqn: phij gj} and the $A_j$ uniquely satisfy~\eqref{e:AkSystem}. Taking the interior limit, we get
\begin{align*}
u_A(z^+_0) & = (\K-\onehalf I )[\phi_A](z_0) + \sum_{j=1}^{M}A_j G(z_0-\beta_j)\\
&= g(z_0) - \left(\X[\phi^0] (z_0)- \sum_{j=1}^M A_j \X[\phi^j] (z_0)\right)= g(z_0),
\end{align*} and conclude that $u_A$ is the solution to the Dirichlet BVP.  

Finally, the flux across $\partial D_j$ for $j= 1,2,\ldots, M$ is given by
\begin{align*}
   \int_{\partial D_j} \partial_\nu u(z^+)\,|dz| & = -\int_{D_j} \Delta \mathcal{D}[\phi](z)\, dA + \sum_{k=1}^M A_k\int_{D_j} \partial_\nu G(z-\beta_k)\, |dz| \\
   & =  A_j +\sum_{k=1}^M A_k\left(-\frac{|D_j|}{b}\right)  = A_j.
\end{align*}
Here, we used Lemma~\ref{lem: exact double layer} and the condition $\sum_{k=1}^M A_k = 0$.
\end{proof}

\subsection{The Neumann BVP and  Proof of Theorem~\ref{thm: neumann}}

\begin{proof}[Proof of Theorem~\ref{thm: neumann}] 
Let  $u(z) = \mathcal{S}[\phi](z)$ as in~\eqref{e:uRep, neu}. 
Calculating $\partial_\nu u(z^+_0)$ by using Lemma~\ref{lem: properties single layer}(2), we obtain 
\begin{equation*}
\K^*[\phi](z_0)+\frac{\phi(z_0)}{2} = g(z_0), \qquad z_0 \in \partial \Omega.
\end{equation*}By Lemma~\ref{lem: properties double layer}(3) and Lemma~\ref{lem: properties double layer}(6), $\K^{*} \colon C_0(\partial \Omega) \to C_0(\partial \Omega)$ is compact and $\K^{*} + \onehalf I$ is injective (here, $C_0(\partial\Omega)$ is a closed subspace of $C(\partial\Omega)$). Hence, by the Fredholm alternative, $\K^{*} + \onehalf I$ is surjective. Therefore, we conclude that for any $g\in C_0(\Omega)$, there is a unique $\phi\in C_0(\partial \Omega)$ such that $u(z)$ is the solution to the Neumann BVP~\eqref{eqn: neumann}, up to an additive constant. 
Finally, the flux across $\partial D_j$ is given by
{\begin{align*}
    \int_{\partial D_j} \partial_\nu u(z^+)\,|dz|& = \int_{\partial D_j} \partial_\nu \mathcal{S}(z^+)\,|dz| = \left\langle (\K^* + \onehalf I )[\phi], \mathbbm{1}_j \right\rangle = \left\langle \phi, (\K + \onehalf I )[\mathbbm{1}_j] \right\rangle \\
    &=  \left\langle \phi,  \frac{\mathbbm{1}_j}{2}-\frac{|D_j|}{b}\mathbbm{1} +\frac{\mathbbm{1}_j}{2} \right\rangle  =  \int_{\partial D_j} \phi (\xi)\, |d\xi|,
\end{align*}}
using Lemma~\ref{lem: exact double layer}, Lemma~\ref{lem: properties single layer}(2), and the fact that $\phi\in C_0(\partial\Omega)$.
\end{proof}

\subsection{The Steklov EVP and  Proof of Theorem~\ref{thm: steklov}}
To prove Theorem~\ref{thm: steklov}, we use the modified single-layer potential in~\eqref{eqn: single-layer Kress formulation} to find a density $\phi \in C(\partial\Omega)$ instead of $\phi \in C_0(\partial\Omega)$.  Second, we match the boundary data of $\sigma u(z)$ with $\partial_\nu u(z)$.
To proceed, we define the Neumann-to-Dirichlet map.
\begin{prop}
\label{l:NtD}
The \emph{Neumann-to-Dirichlet map}, $\Lambda \colon C_0(\partial \Omega) \to  C_0^1(\partial \Omega)$, can be expressed in terms of layer potentials as
\begin{equation}
\Lambda \coloneq (I-\M ) \S \left(\K^*+ \onehalf I \right)^{-1}, 
\label{eqn: ntd}
\end{equation}
which, given  Neumann boundary data, outputs zero mean  Dirichlet boundary data. 
\end{prop}
The term $(I - \M)$ enforces the Dirichlet data to have zero mean on the boundary. We adopt this convention to ensure the uniqueness of the Dirichlet solution, since Neumann boundary data determines a harmonic function only up to an additive constant.

\begin{proof}[Proof of Proposition~\ref{l:NtD}]
Let $g \in C_0 (\partial\Omega)$ be Neumann data. By Theorem~\ref{thm: neumann}, there exists a unique $\psi\in C_0(\partial\Omega)$ satisfying~\eqref{eqn: neumann bdry eqn}. 
By Theorem~\ref{thm: neumann}, 
$\mathcal{S}[\psi]+ C$ solves the Neumann BVP~\eqref{eqn: neumann} with boundary data $g$, for any $C\in \mathbb{R}$. We fix the additive constant by applying $(I-\M)$; that is, $C$ is the average value of $-\S[\psi]$ on $\partial\Omega$.  
\end{proof}

\begin{proof}[Proof of Theorem~\ref{thm: steklov}] 
Suppose $(\sigma_k, u_k)$ is an eigenpair of the Steklov EVP \eqref{eqn: steklov} with $k \geq 1$ such that $\sigma_k > 0$. Applying Proposition~\ref{l:NtD}, there exists a unique $\psi_k \in C_0(\partial\Omega)$ such that $(I-\M)\S[\psi_k]$ is a solution to the Neumann BVP with data $\partial_\nu u_k$. Define $\phi_k \in C(\partial \Omega)$ by $\phi_k(z_0) = \psi_k(z_0) - \M \S[\psi_k]$ so that $\psi_k = (I-\M) \phi_k$ and $\overline{\phi_k} = - \M\S [\psi_k]$. 
By construction,
\begin{align*}
(\K^* + \onehalf)(I- \M)[\phi_k]= \partial_\nu u_k = \sigma_k u_k = \sigma_k (I-\M)\S[\psi_k] =  \sigma_k\S_0[\phi_k],
\end{align*}
proving the existence of a $\phi_k$  such that $(\sigma_k, \phi_k)$ satisfies \eqref{e: steklov bie}.

Conversely, suppose $(\sigma_k, \phi_k)$ with $\phi_k \in C(\partial\Omega)$ satisfies~\eqref{e: steklov bie}. Define $u_k(z) = \mathcal{S}_0[\phi_k] (z),$
which is harmonic on $\Omega$. $\partial_\nu u_k$ is  continuous on $\partial\Omega$, since $\partial_\nu u_k = (\K^{*}+\onehalf I )(I-\M)\phi_k$ by Lemma~\ref{lem: properties single layer}. Thus, $u_k$ solves the Steklov EVP~\eqref{eqn: steklov} with eigenvalue $\sigma_k$. 
\end{proof}

\section{Numerical Experiments}
\label{sec: Numerical Experiments}
In this section, we present numerical experiments for approximating solutions of the  Dirichlet BVP~\eqref{eqn: dirichlet}, Neumann BVP~\eqref{eqn: neumann}, and Steklov EVP~\eqref{eqn: steklov}. In each of the examples,  we consider tori with periods $(1, \tau)$ with $\tau = i$ (square torus) and $\tau = \frac{1}{2} + \frac{\sqrt{3}}{2} i$ (equilateral torus). Our implementation utilizes a standard trapezoid Nystr{\"o}m discretization~\cite{hao2011high}\cite[Ch. 12]{kress1989linear} based on the layer potential methods developed in Theorems~\ref{thm: dirichlet}--\ref{thm: steklov}. All numerical results have been obtained using Matlab direct solvers on fully discretized matrices. We utilize the backslash operator (\texttt{$\backslash$}) for the Dirichlet~\eqref{e:DirIntEq} and Neumann BVPs~\eqref{eqn: neumann bdry eqn}, and the 
\texttt{eig} function for the Steklov EVP~\eqref{e: steklov bie}. The complete implementation of our computational methods is available on GitHub~\cite{bohyunBIEgithub}. The runtime for solving the linear systems to compute the density $\phi$ in Examples 1--5 is between 10 and 100 ms on a 2020 MacBook Pro M1. The runtime for Example 6  is $19.6$ minutes.

\subsection{Dirichlet boundary value problem} 
\label{sec: dirichlet num}
In the following examples, we solve the Dirichlet BVP~\eqref{eqn: dirichlet} on a multiply-connected torus, $\Omega =  \mathbb{T}_\tau \setminus \cup_{j=1}^M \overline{D_j}$, by varying the number of holes $M$, the period $\tau$, and the Dirichlet boundary data $g$.

\begin{figure}[t!]
\begin{center}
\includegraphics[width=0.28\textwidth]{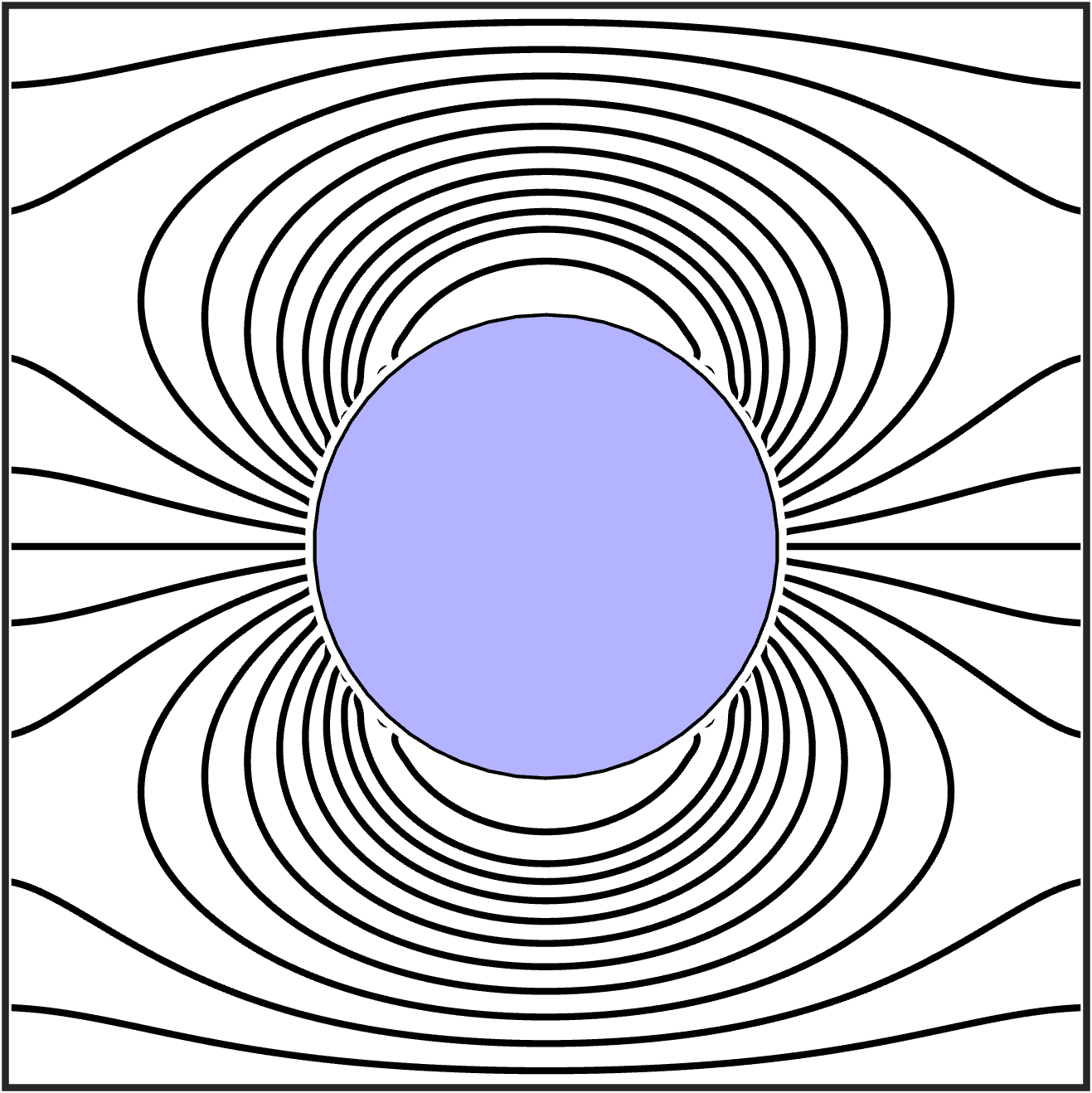}
\includegraphics[width=0.48\textwidth]{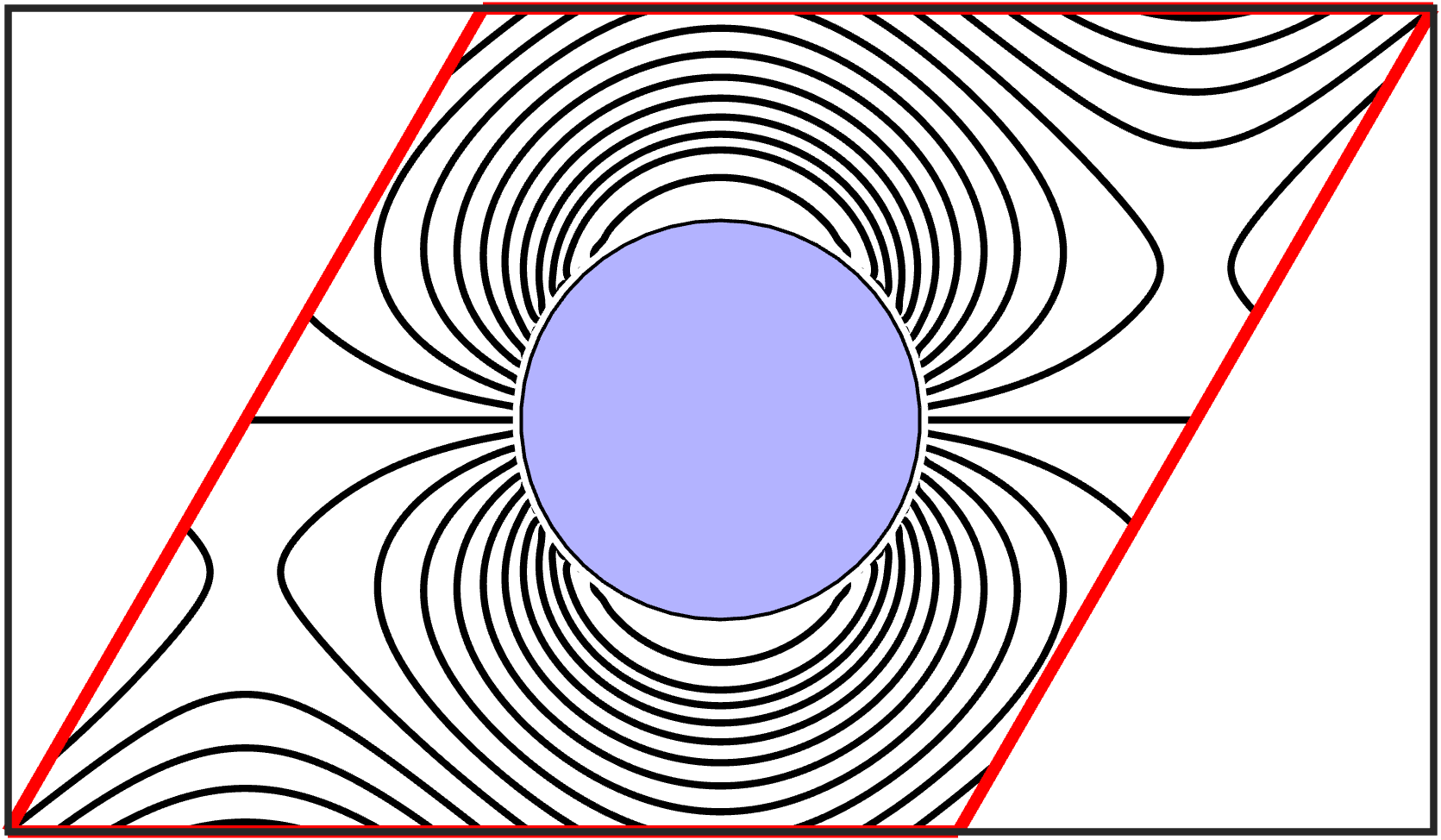}
\end{center}
\caption{Approximate Dirichlet BVP solutions on a square torus {\bf~(left)} and an equilateral torus {\bf (right)} with a single circular hole. See Sec.~\ref{sec: dirichlet num}, Example 1.}
\label{fig: dirchlet simp circ}
\end{figure}

\vspace{3mm} 
\noindent {\bf Example 1. Dirichlet BVP with a single ($M=1$) hole.} 
We solve the Dirichlet BVP on tori with a single circular hole, as shown in Fig.~\ref{fig: dirchlet simp circ}.
Let $D = B(a_1, r)$ be the disk of radius $r = 0.2$ centered at $a_1 = 0.5 + 0.5i$. 
The boundary data is $g = -\S[\psi ]$ with $\psi(\xi) = \sin \arg (\xi - a_1)$ for $\xi \in \partial D$. By construction, the exact solution is given by 
$u_{\text{exact}}(z) = -\mathcal{S}[\psi](z)$ for $z \in \Omega$. Note that the kernel of the single-layer potential exhibits a logarithmic singularity near the boundary. To address this, we utilize a splitting technique similar to those in \cite[Sec. 12.3]{kress1989linear} and \cite{kress1991boundary}. We decompose $G$ into a logarithmically singular part, $-\frac{1}{4\pi}\log(4\sin^2((s-t)/2))$, and a remainder that is continuous as $s\to t$. The implementation details are available on GitHub~\cite{bohyunBIEgithub}.

To approximate the solution, consistent with Theorem~\ref{thm: dirichlet}, we use  the representation \( u(z) = \mathcal{D}[\phi](z) \). 
We use $N=50$ points to discretize the boundary integral equation~\eqref{e:DirIntEq}.
The evaluation of the double-layer potential near the boundary $\partial \Omega$ is relatively inaccurate due to discretization error. 
Therefore, to evaluate the error in the numerical computation, we define $E = \| u - u_{\text{exact}}\|_{L^{\infty}(\gamma)}$, where $\gamma$ is the boundary of an enlarged disk, $\partial B(a_1,r+0.15)$. With $N=50$ grid points, the errors obtained are 
$E = 1.604\times 10^{-13}$ (square torus) and 
$E = 1.605 \times 10^{-13}$ (equilateral torus).

\begin{figure}[t!]
\begin{center}
\includegraphics[width=0.28\textwidth]{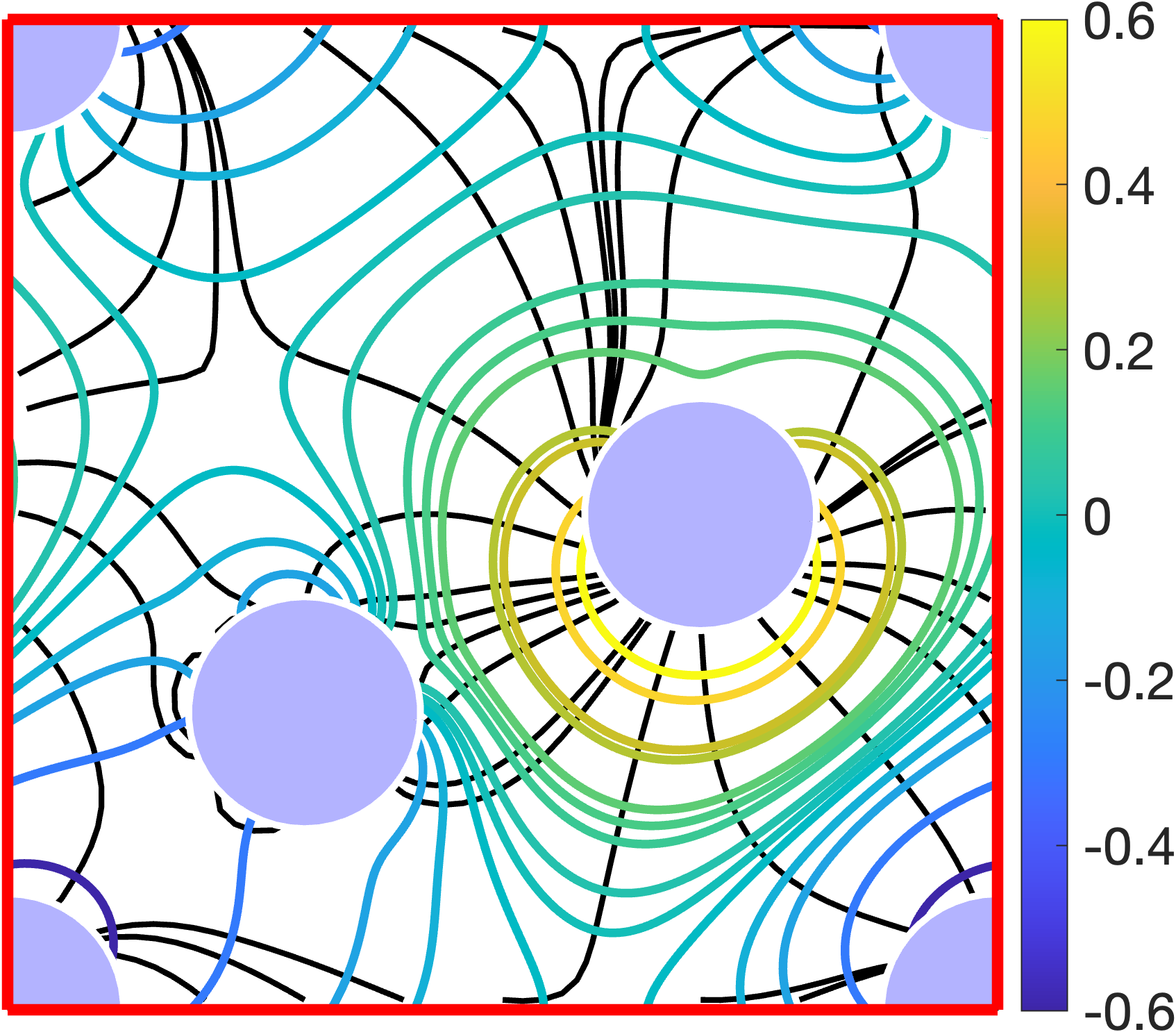}
\includegraphics[width=0.48\textwidth]{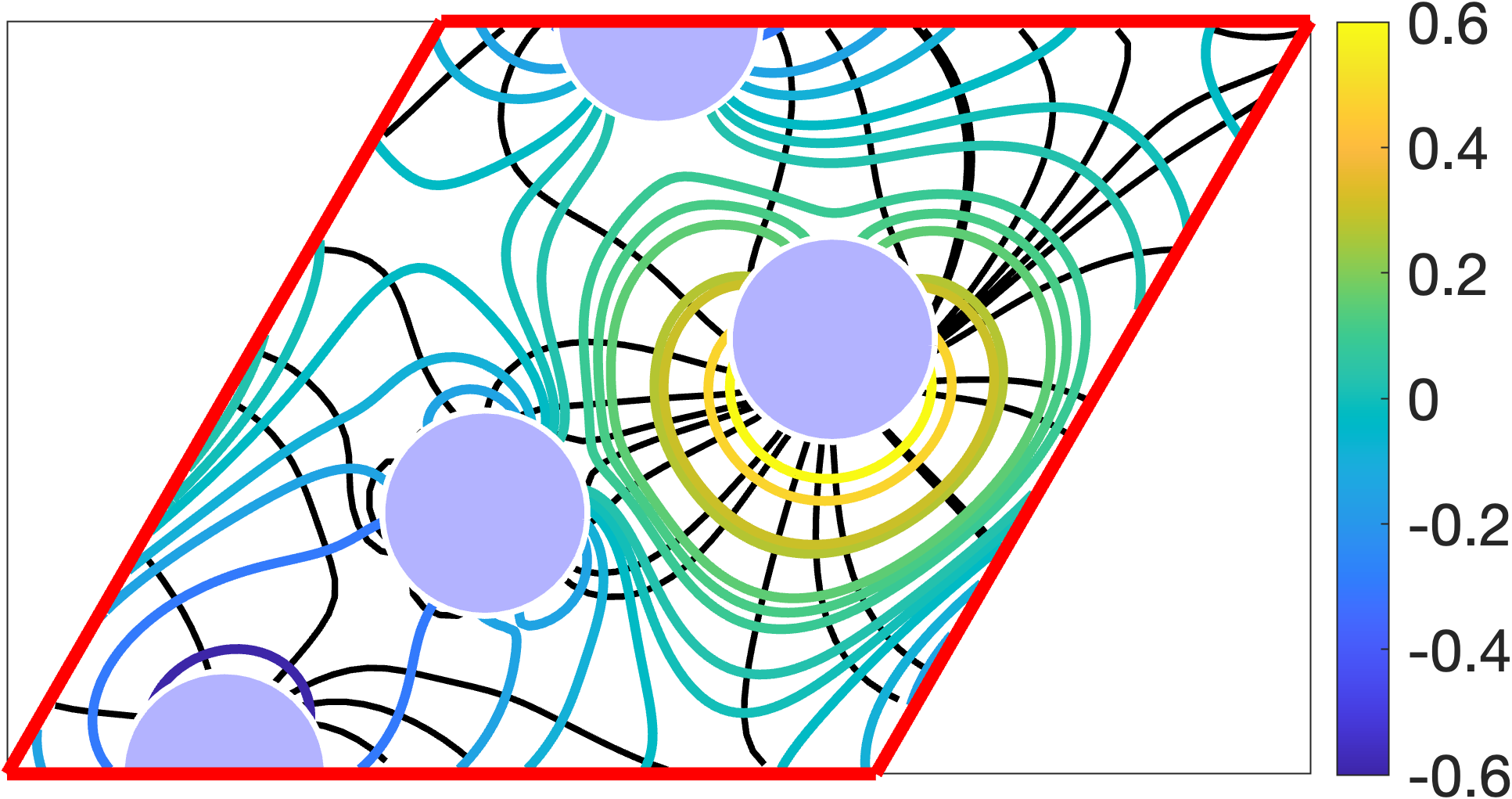} \\ 
\includegraphics[width=0.28\textwidth]{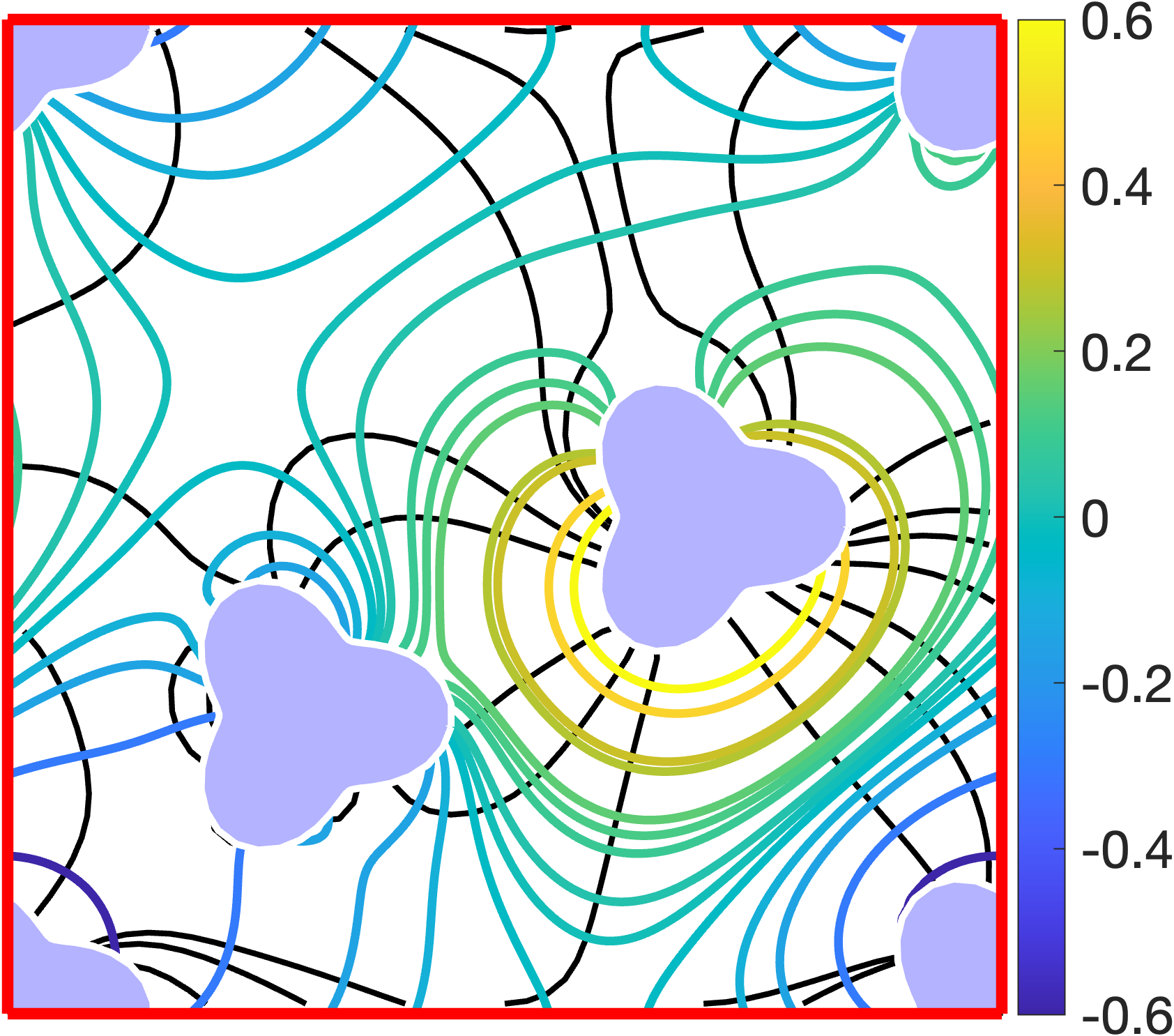}
\includegraphics[width=0.48\textwidth]{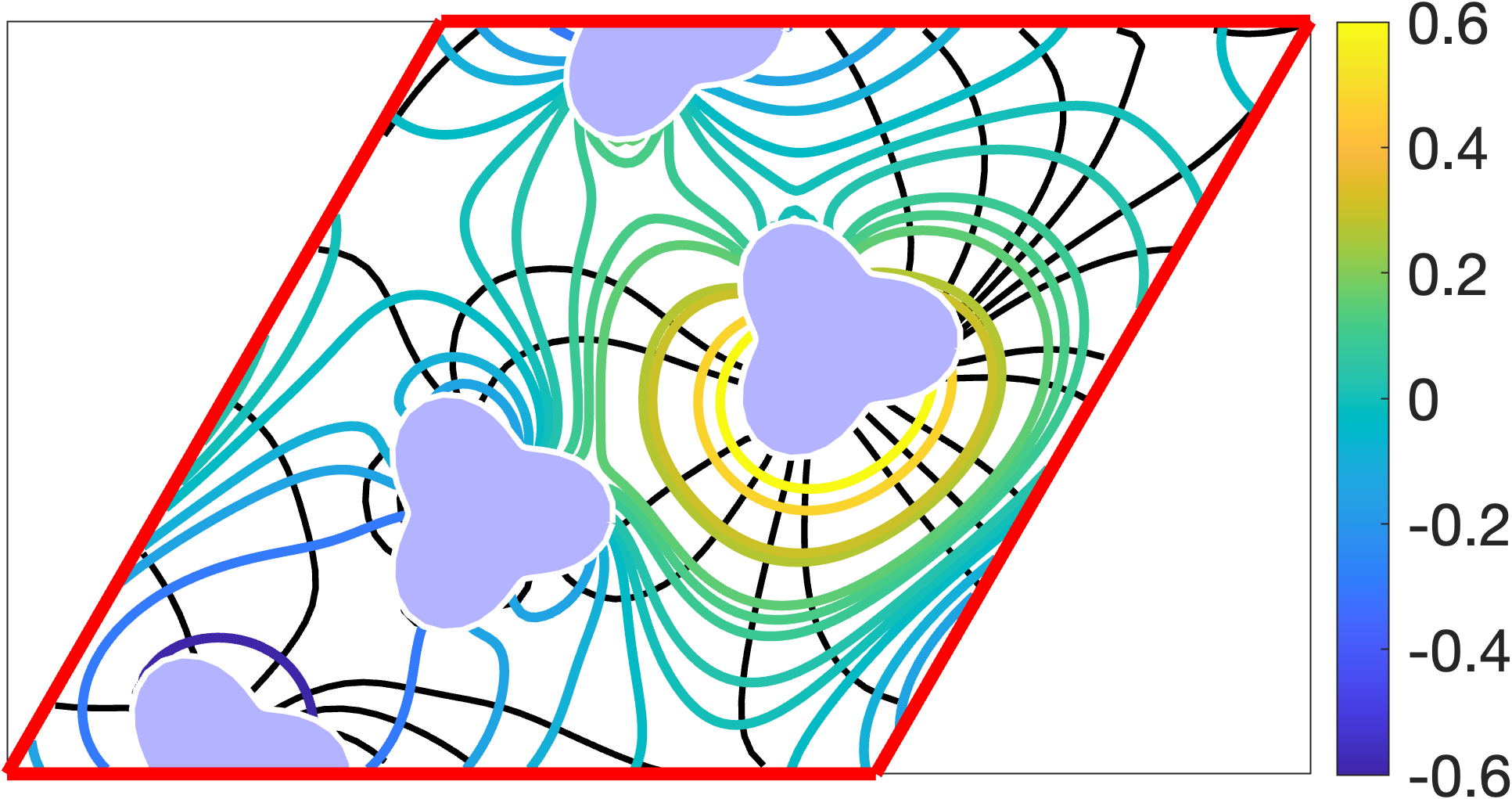} \\ 
\includegraphics[width=0.382\textwidth]{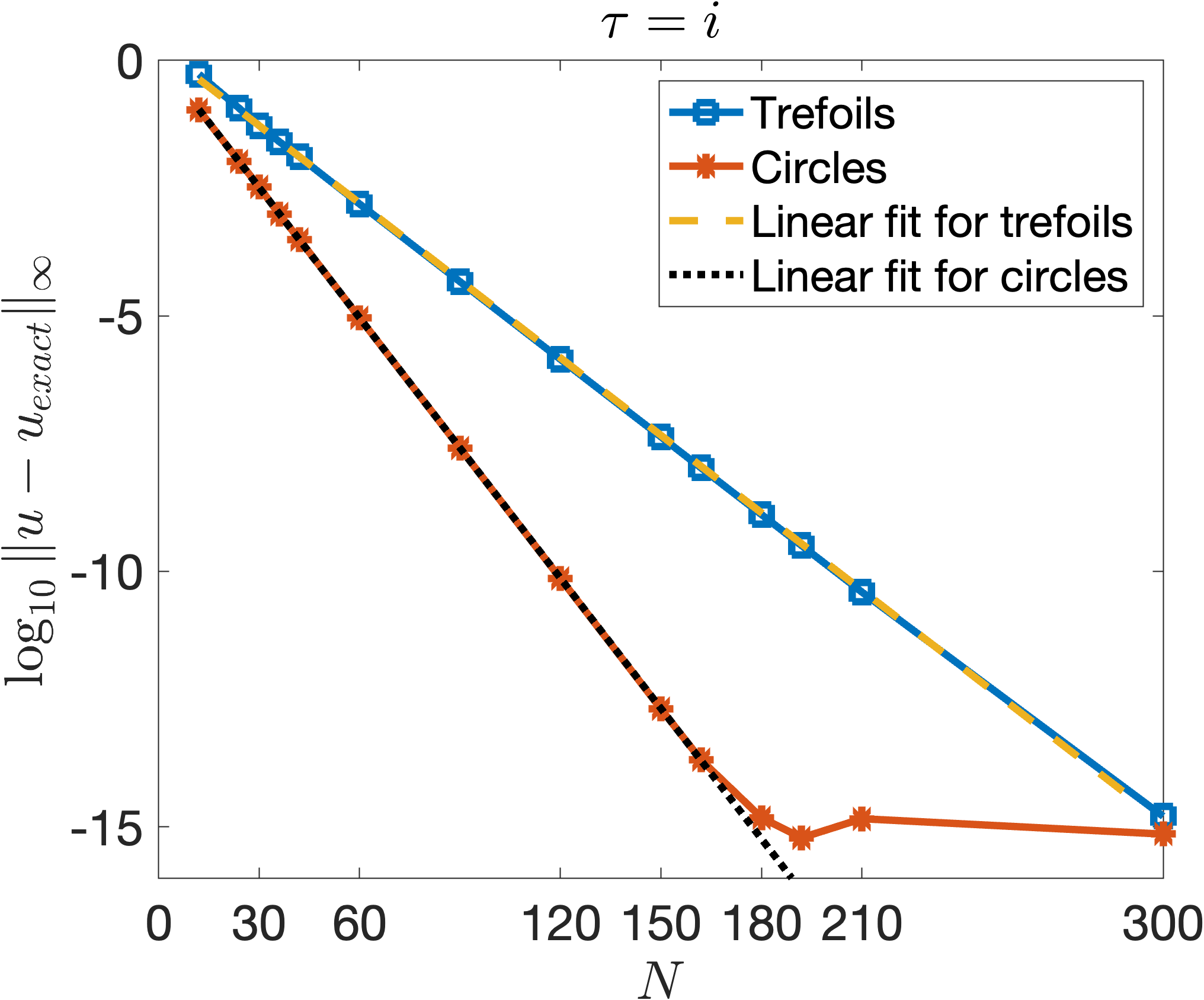}
\includegraphics[width=0.382\textwidth]{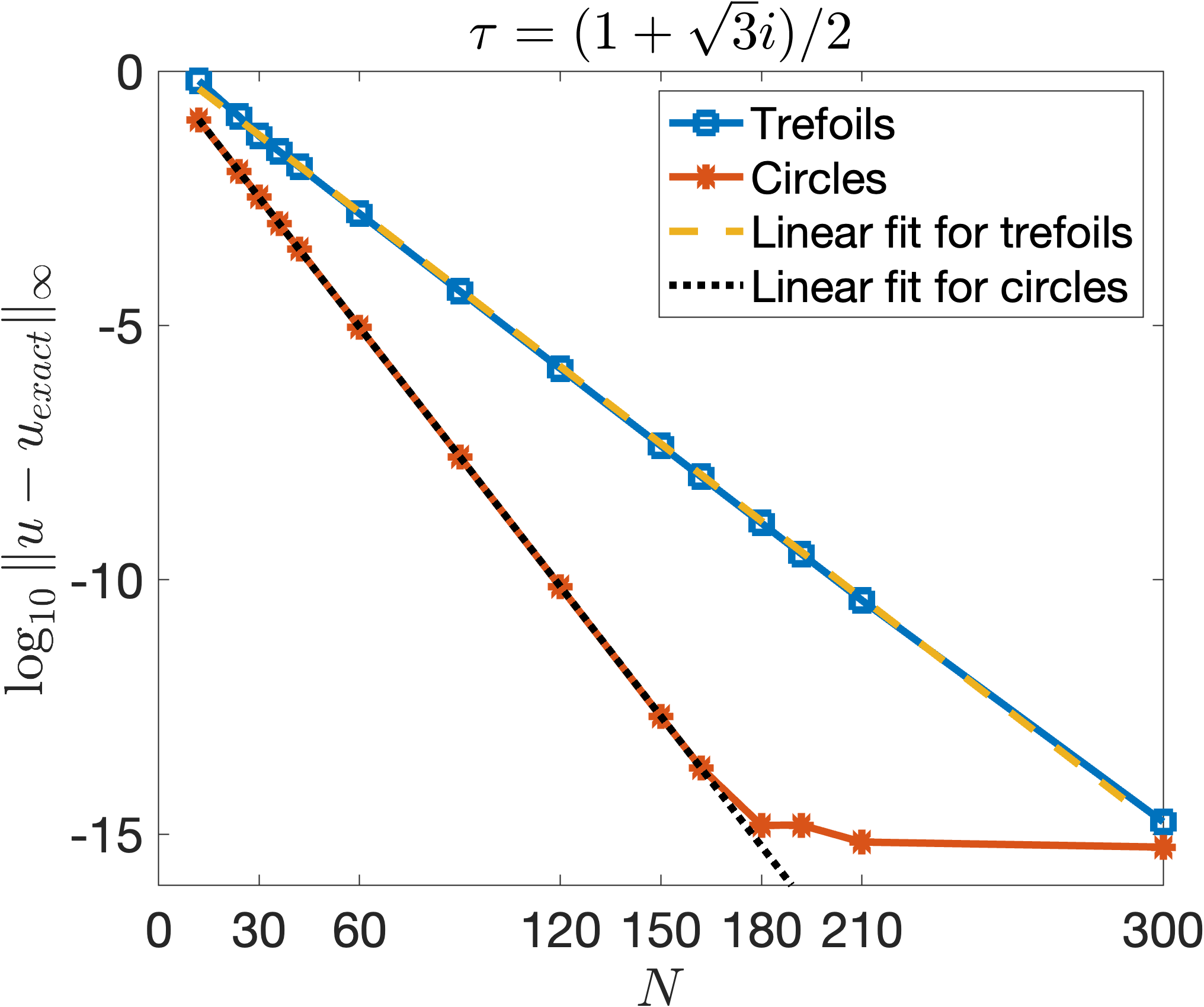}
\end{center}
\caption{Approximate Dirichlet BVP solutions on the square torus {\bf~(left)} and equilateral torus {\bf~(right)} with $M=3$ circular holes {\bf (top)}  and trefoil-shaped holes {\bf (center)}. 
{\bf (bottom)} The convergence plots illustrate that the numerical methods converge spectrally for both the three circular holes and the three trefoil-shaped holes on square tori {\bf~(left)} and equilateral tori {\bf~(right)}. The convergence rates for the two shaped holes are similar, but slightly better for the circles. See Sec.~\ref{sec: dirichlet num}, Example 2.}
\label{fig: dirchlet mult holes}
\end{figure}

\vspace{3mm} 
\noindent {\bf Example 2. Dirichlet BVP with $M=3$ holes.}
We solve the Dirichlet BVP on tori with 
$M=3$ circular and ``trefoil''-shaped holes; see Fig.~\ref{fig: dirchlet mult holes}. Both types of holes are centered at $a_1 = 0.7+0.5i$, $a_2 = 0.3+0.3i$, and $a_3 = 0$ with radius $r = 0.1$. The trefoils are parameterized by 
\begin{equation}
\label{e:tref}
 f(z,r,a_i) \coloneq r \left( 1 + 0.3\cos \left(3 \arg(z-a_i) \right) \right), 
 \qquad z \in \mathbb{T}_\tau.
\end{equation}
For both the square and equilateral tori, define 
\begin{align*}
    g(z) & = \mathcal{S}[\psi ](z) + 3 G(z-a_1)-G(z-a_2)-2G(z-a_3), \quad z\in \Omega, \\
    \text{where} \quad \psi(\xi) &=\begin{cases} -10\sin \arg (\xi - a_1) & \text{if } \xi \in \partial D_1,\\
10 \sin (3\arg (\xi - a_2)) & \text{if } \xi \in \partial D_2,\\
-10 \sin \arg (\xi - a_3) & \text{if } \xi \in \partial D_3.
\end{cases}
\end{align*} 
Using $g\big|_{\partial \Omega}$ as boundary Dirichlet data, by construction, the exact solution is given by 
$u_{\text{exact}}(z) = g(z)$ for $z \in \Omega$.

To approximate the solution, we represent the solution as in~\eqref{e:uRep}. We use $N_i=50$ points per boundary component ($N=150$ total) to discretize the boundary integral equation~\eqref{e:DirIntEqA}, plus three additional equations to satisfy the conditions in~\eqref{e:DirIntEqB}. Consequently, system~\eqref{e:DirIntEq} is represented by a $153\times153$ matrix equation.  
As above, to evaluate the error, we define 
$E = \| u - u_{\text{exact}}\|_{L^{\infty}(\gamma)}$, where $u_{\text{exact}}$ is calculated on a fine grid with $N = 3750$. For the circular holes in Fig.~\ref{fig: dirchlet mult holes}, we use 
$\gamma = {\cup_{i=1}^3}\partial B(a_i,r+0.08)$ and obtain  $E =2.032\times 10^{-13}$~(left) and  $E =  2.028\times 10^{-13}$(right). For the trefoil holes in Fig.~\ref{fig: dirchlet mult holes}, we use 
$\gamma$ as in~\eqref{e:tref} with $r=0.18$ 
and obtain $E = 4.300\times 10^{-8}$~(left) and $E = 4.299\times 10^{-8}$ (right).

In the bottom plots of Fig.~\ref{fig: dirchlet mult holes}, we plot the convergence of the error for increasing $N=4,10,20,\ldots, 250$. 
These plots illustrate that the numerical methods converge spectrally for both the circular and trefoil-shaped holes. 
The convergence rate for the trefoil-shaped holes is not significantly degraded from the convergence rate for the circular holes; a linear fit for the slope in the log plot yields 
$-0.0505$ (trefoils) and $-0.0850$ (circles) for the square torus, and 
$-0.0506$ (trefoils) and $-0.0850$ (circles) for the equilateral torus.  
For irregularly shaped holes, our layer potential approach exhibits a better convergence rate than the method of particular solutions (MPS) given in \cite{kaoharmonic2023}. In particular,  the $\log_2$ error plot in \cite[Fig. 4]{kaoharmonic2023} shows that  MPS requires approximately 500 degrees of freedom (dof) to achieve an error that is approximately $3\times10^{-6}$ on a square torus with two trefoil-shaped holes. Our method achieves comparable error with $N=120$ dof on a square torus with three trefoil-shaped holes.

In Table~\ref{tabl: dirc flux}, we report the values of the fluxes $A_i$ for $i=1,2,3$ for both the trefoil and circular holes on the square~(left) and equilateral torus (right). For both cases, the accuracy of the flux is lower for the trefoil holes than for the circular holes. Because the contribution to flux from $\S[\psi]$ is zero around each hole, the estimated fluxes are very close to the coefficient of $G(z-a_i)$ in $g(z)$.

\begin{table}[t!]
    \centering
    {\small
  \begin{tabular}{|c|c|c|c|}
    \hline
     & Flux, $A_i$, & Flux errors & Flux errors   \\
    i & for $D_i$ & (trefoil) &(circle) \\ 
    \hline
    1 & 3  & $6.161 \times 10^{-11}$  & $4.400 \times 10^{-16}$ \\
    \hline
    2 & -1 & $4.411 \times 10^{-11}$  & $2.200 \times 10^{-16}$ \\
    \hline
    3 & -2 & $1.057 \times 10^{-10}$ & $2.200 \times 10^{-16}$ \\
    \hline
\end{tabular}
\begin{tabular}{|c|c|c|}
    \hline
     & Flux errors & Flux errors   \\
    i & (trefoil)  &(circle) \\ 
    \hline
    1 & $7.849 \times 10^{-11}$  & $1.330 \times 10^{-15}$ \\
    \hline
    2 & $7.463 \times 10^{-11}$  & $8.900 \times 10^{-16}$ \\
    \hline
    3 & $3.860 \times 10^{-12}$  & $2.220 \times 10^{-15}$ \\
    \hline
\end{tabular}
    \vspace{0.1em}
    }
   \caption{Errors of fluxes  for the square {\bf~(left)} and equilateral {\bf~(right)} tori in Example 2 (Dirichlet BVP with $M = 3$ holes). See Sec.~\ref{sec: dirichlet num}  and Fig.~\ref{fig: dirchlet mult holes}.}
    \label{tabl: dirc flux}
\end{table}

\begin{figure}[t!]
\begin{center}
\includegraphics[width=0.28\textwidth]{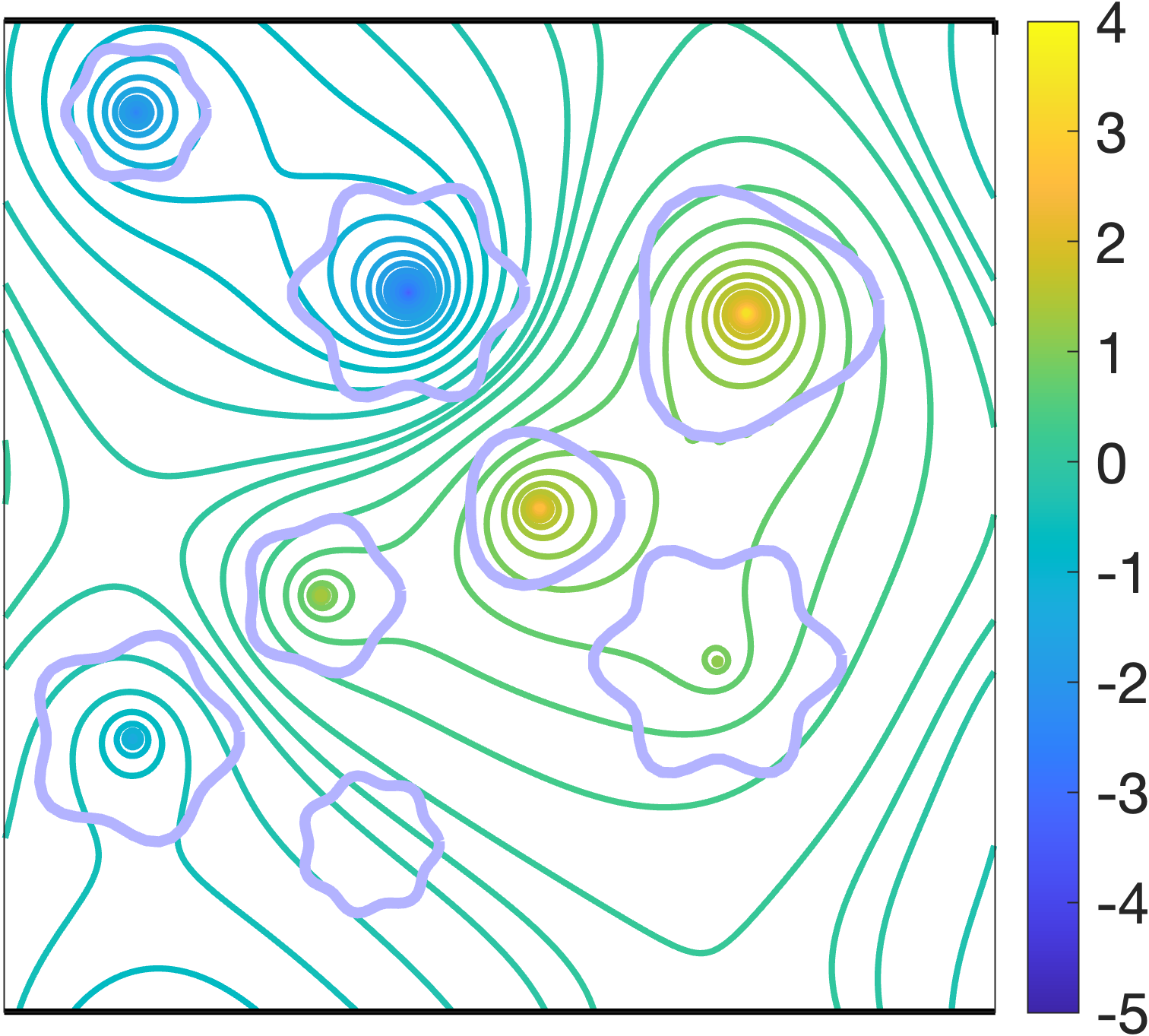}
\includegraphics[width=0.48\textwidth]{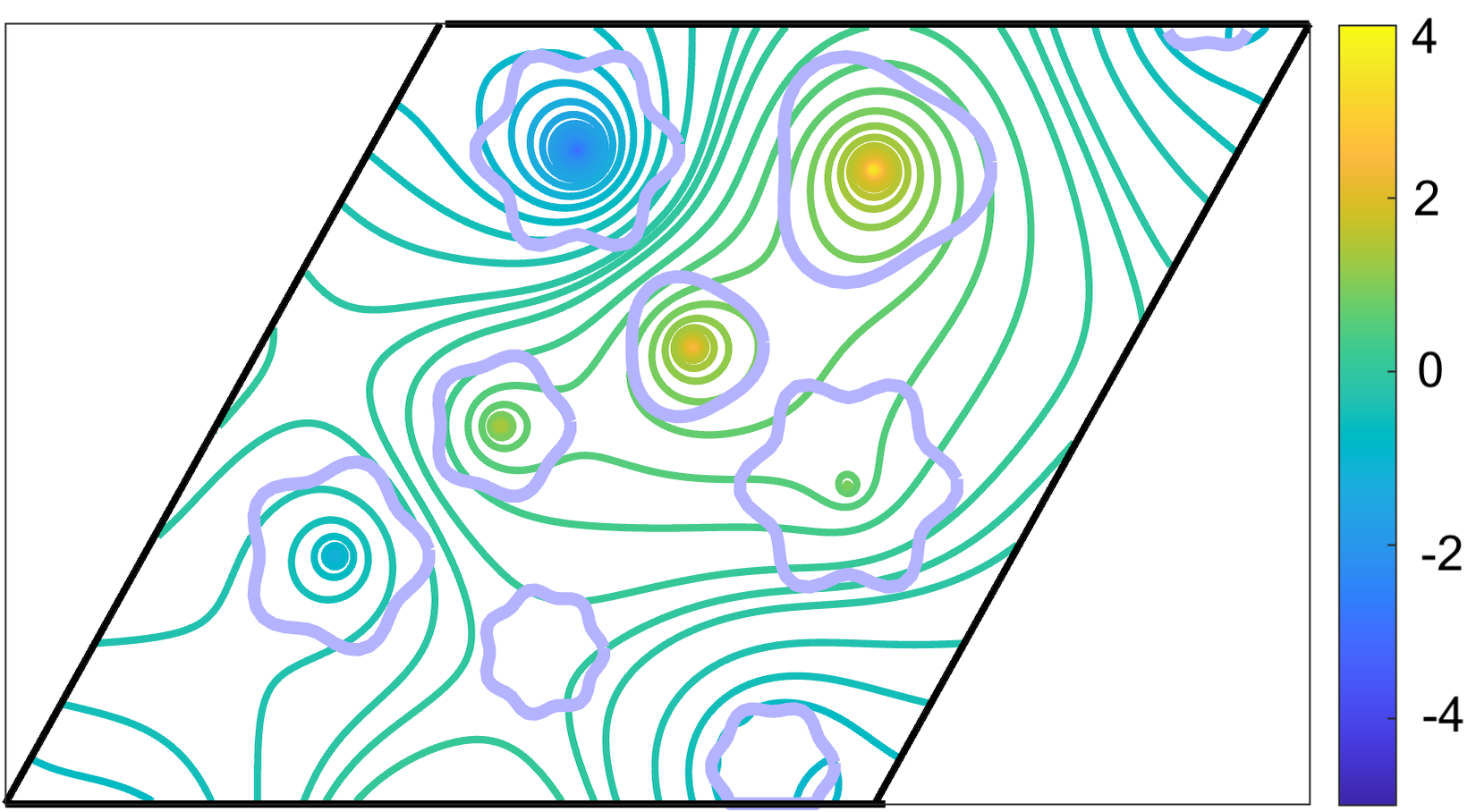}
\end{center}
\caption{Approximate solutions to the Neumann BVP on the square torus {\bf~(left)} and the equilateral torus {\bf (right)} with multiple holes. The boundary $\partial \Omega$ is drawn in purple. Since $\mathcal{S}[\phi]$ is continuous across the boundary, we have also plotted $\mathcal{S}[\phi]$ inside $D$. See Sec.~\ref{s: neumann bvp}, Example 3 and Table~\ref{tabl: neumann flux}.}
\label{fig: neumman mult}
\end{figure}

\subsection{Neumann boundary value problem} 
\label{s: neumann bvp}
\ \\ 
\noindent {\bf Example 3. Neumann BVP with $M=8$ holes.} 
We solve the Neumann BVP on tori with $M=8$ holes; the results are shown in Fig.~\ref{fig: neumman mult}. To generate the holes, we randomly chose each center $a_i$ for $i=1,2,\ldots,8$, an oscillation factor $\omega_i$ with $3\leq \omega_i\leq7$, and the maximum radius of the hole \( r_i \) (given in Table~\ref{tabl: neumann flux} (left)). The boundary of the hole \( \partial D_i \) is parametrized by 
\begin{equation}
 f(z,r_i,\omega_i, a_i) \coloneq  \frac{r_i}{r_i+1} (1+r_i\cos\omega_i\arg (z - a_i)), \qquad  z\in \mathbb{T}_\tau.
    \label{eqn: neumann bdry shape}
\end{equation}
The boundary data is chosen as $g(z_0) = \sum_{i=1}^8 A_i \partial_\nu G(z_0-a_i)$ for $z_0\in\partial\Omega$, where the value of $A_i$
for each hole is reported in Table~\ref{tabl: neumann flux}(right).  By construction, the exact solution is given by 
$u_{\text{exact}}(z) =\sum_{i=1}^8 A_i G(z-a_i)+C$ for $z \in \Omega$ and $C\in\mathbb{R}$ is an arbitrary constant.

To approximate the solution, we represent the solution as in~\eqref{e:uRep, neu}. We use $N_i=50$ points per boundary component ($N=400$ total) to discretize the boundary integral equation~\eqref{eqn: neumann bdry eqn}. The approximate solutions are plotted in Fig.~\ref{fig: neumman mult}. 
Since the single-layer potential is continuous across $\partial \Omega$, we plot the solution on the entire $\mathbb{T}_\tau$ (the domain $\Omega$ and the removed holes). 

We evaluate the error in the numerical computation at the boundary $\partial\Omega$ as
$
E = \| u - u_{\text{exact}} + C \|_{L^\infty(\partial\Omega)}$, where $C$ is chosen so that $u$ and $u_{\text{exact}}$ have the same mean. We compute $C$ by averaging $u - u_{\text{exact}}$  on $250$ randomly chosen points on $\Omega$. 
As in Example 1, we calculate $\S[\phi]$ by decomposing the kernel $G$ into a logarithmically singular part, $-\frac{1}{4\pi}\log(4\sin^2((s-t)/2))$, and a remainder that is continuous as $s \to t$.
We compute $E = 2.100\times 10^{-7}$ (square torus) and 
$E = 2.100\times 10^{-7}$ (equilateral torus).
In Table~\ref{tabl: neumann flux} (right), we report the estimated flux errors,  $\left|A_i - \int_{\partial D_i} \phi(\xi)\,|d\xi| \right|$ for $i=1,2,\ldots, 8$; see Theorem~\ref{thm: neumann}. For both values of $\tau$, the estimated fluxes are very close to the coefficients $A_i$ in $u_{\text{exact}}$. 

\begin{table}[t!]
    \centering
    {\small
    \begin{tabular}{|c|c|c|c|}
    \hline
i & $\omega_i$ & $r_i$ & $a_i$ \\
& & & \\
\hline
1 & 6 & 0.126 & 0.720 + 0.353i \\
\hline
2 & 5 & 0.081 & 0.320 + 0.420i \\
\hline
3 & 3 & 0.082 & 0.540 + 0.508i \\
\hline
4 & 3 & 0.135 & 0.749 + 0.704i \\
\hline
5 & 6 & 0.118 & 0.408 + 0.725i \\
\hline
6 & 5 & 0.108 & 0.130 + 0.276i \\
\hline
7 & 6 & 0.071 & 0.133 + 0.907i \\
\hline
8 & 7 & 0.071 & 0.369 + 0.169i \\
\hline
    \end{tabular}
    \quad 
       \begin{tabular}{|c|c|c|}
    \hline
     $A_i$ & \text{Flux Errors}  & \text{Flux Errors} \\
      & $\tau = i$ &  $\tau = 1/2+\sqrt{3}i/{2}$\\
    \hline
     1 & $1.284 \times 10^{-10}$ & $1.318 \times 10^{-10}$ \\
    \hline
     2 & $3.720 \times 10^{-12}$ & $4.510 \times 10^{-12}$ \\
    \hline
     3 & $2.212 \times 10^{-11}$ & $2.293 \times 10^{-11}$ \\
    \hline
     4 & $9.690 \times 10^{-12}$ & $1.168 \times 10^{-11}$ \\
    \hline
     -5 & $1.588 \times 10^{-11}$ & $1.744 \times 10^{-11}$ \\
    \hline
     -2 & $3.026 \times 10^{-11}$ & $3.159 \times 10^{-11}$ \\
    \hline
     -3 & $2.990 \times 10^{-12}$ & $3.610 \times 10^{-12}$ \\
    \hline
     0 & $2.985 \times 10^{-12}$ & $3.599 \times 10^{-12}$ \\
    \hline
    \end{tabular}
    \vspace{0.5em}
    }
       \caption{\textbf{~(left)} For the torus in Example 3 (Neumann BVP with $M = 8$ holes), we tabulate oscillations, radii, and centers of each hole $D_i$ in~\eqref{eqn: neumann bdry shape}.
       \textbf{(right)} Error in fluxes for the square and equilateral tori. See Sec.~\ref{s: neumann bvp}, Example 3 and Fig.~\ref{fig: neumman mult}.}
    \label{tabl: neumann flux}
\end{table}

\subsection{Steklov eigenvalue problem}
\label{sec: steklovbvp numerical}
In the following examples, we solve the Steklov EVP~\eqref{eqn: steklov} on a multiply-connected torus, $\Omega =  \mathbb{T}_\tau  \setminus \cup_{j=1}^M \overline{D_j}$, by varying the number of holes $M$ and the period $\tau$. For all examples, we approximate the solution using \eqref{e:uRep, stek} in Theorem~\ref{thm: steklov}. In particular, Examples 4 and 5 are compared to those in~\cite{kaoharmonic2023}; since the domain in our problem is half the size of the one in~\cite{kaoharmonic2023}, we divide the eigenvalues by two before comparing errors.

To compute the eigenvalues \( \sigma_k \) and the corresponding density functions \( \phi_k \), we discretize the boundary integral equation~\eqref{e: steklov bie}. The discretization of the adjoint double-layer potential $\K^*$ is performed similarly to $\K$, where the diagonal limit is determined by Lemma~\ref{lem: properties double layer}(3). For the modified single-layer potential $\S_0 = \S (I-\M) + \M$, we discretize the operators $\S$, $I$, and $\M$ as $N \times N$ matrices. Following the approach in Example 1, we resolve the logarithmic singularity in $\S$ by splitting the kernel into a singular part and a smooth remainder.  Using \texttt{eig} function in MATLAB (specifically \texttt{[V, D] = eig(A, B)}), we solve the generalized EVP. We recall that the first eigenvalue is zero.

\begin{figure}[t!]
\begin{center}
\includegraphics[width=0.28\textwidth]{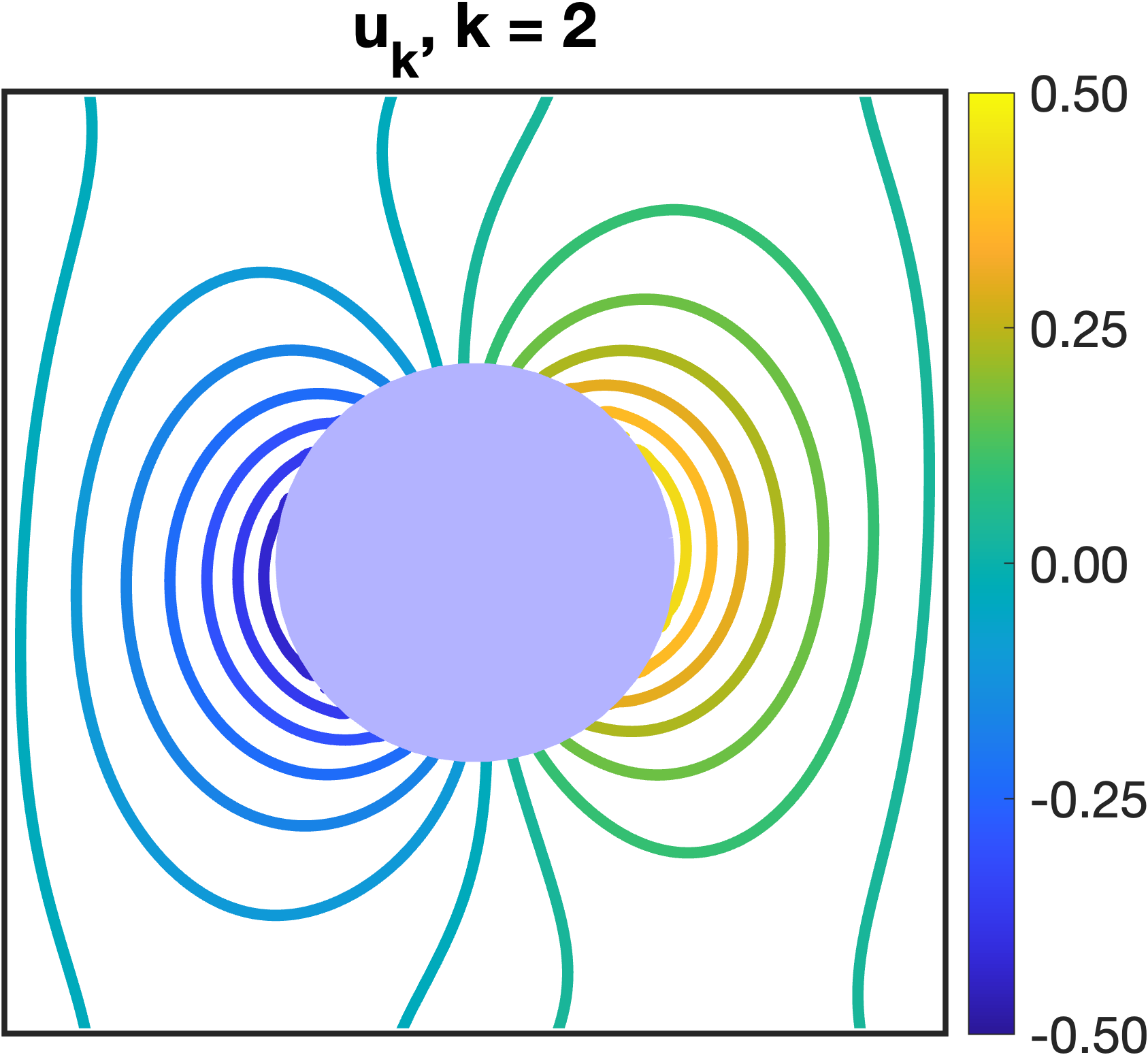}
\includegraphics[width=0.28\textwidth]{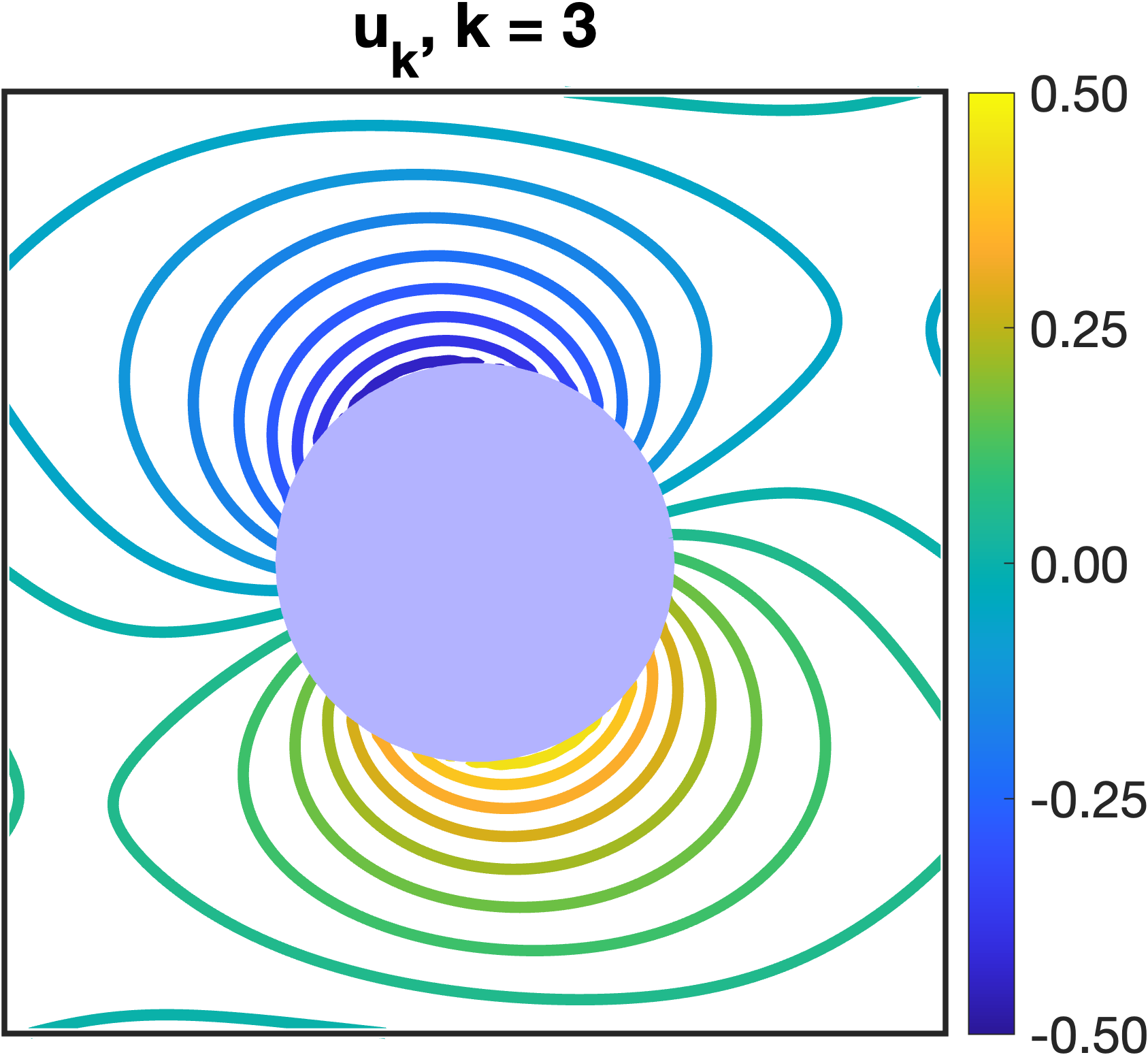}
\includegraphics[width=0.28\textwidth]{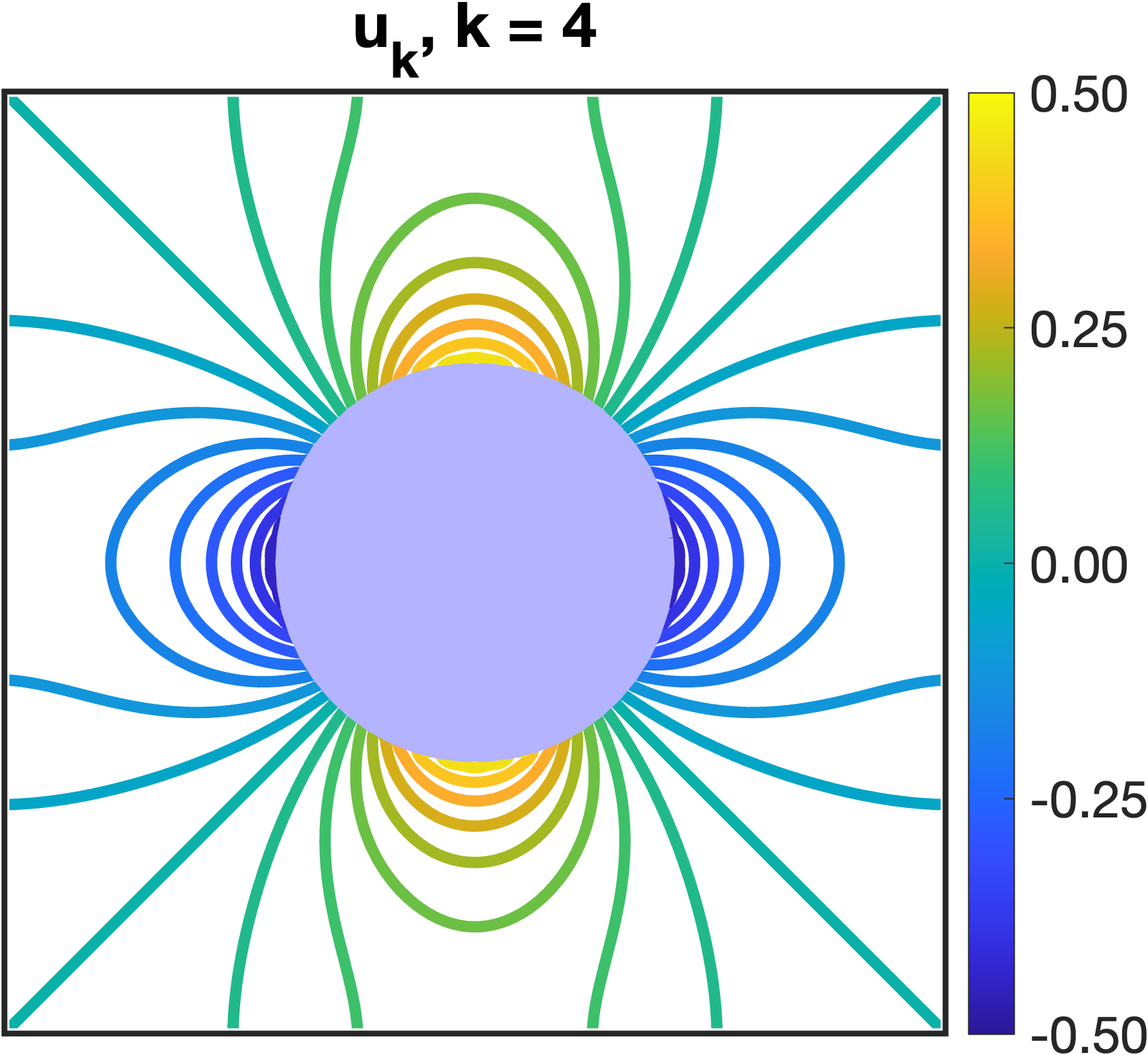} \\ 
\includegraphics[width=0.28\textwidth]{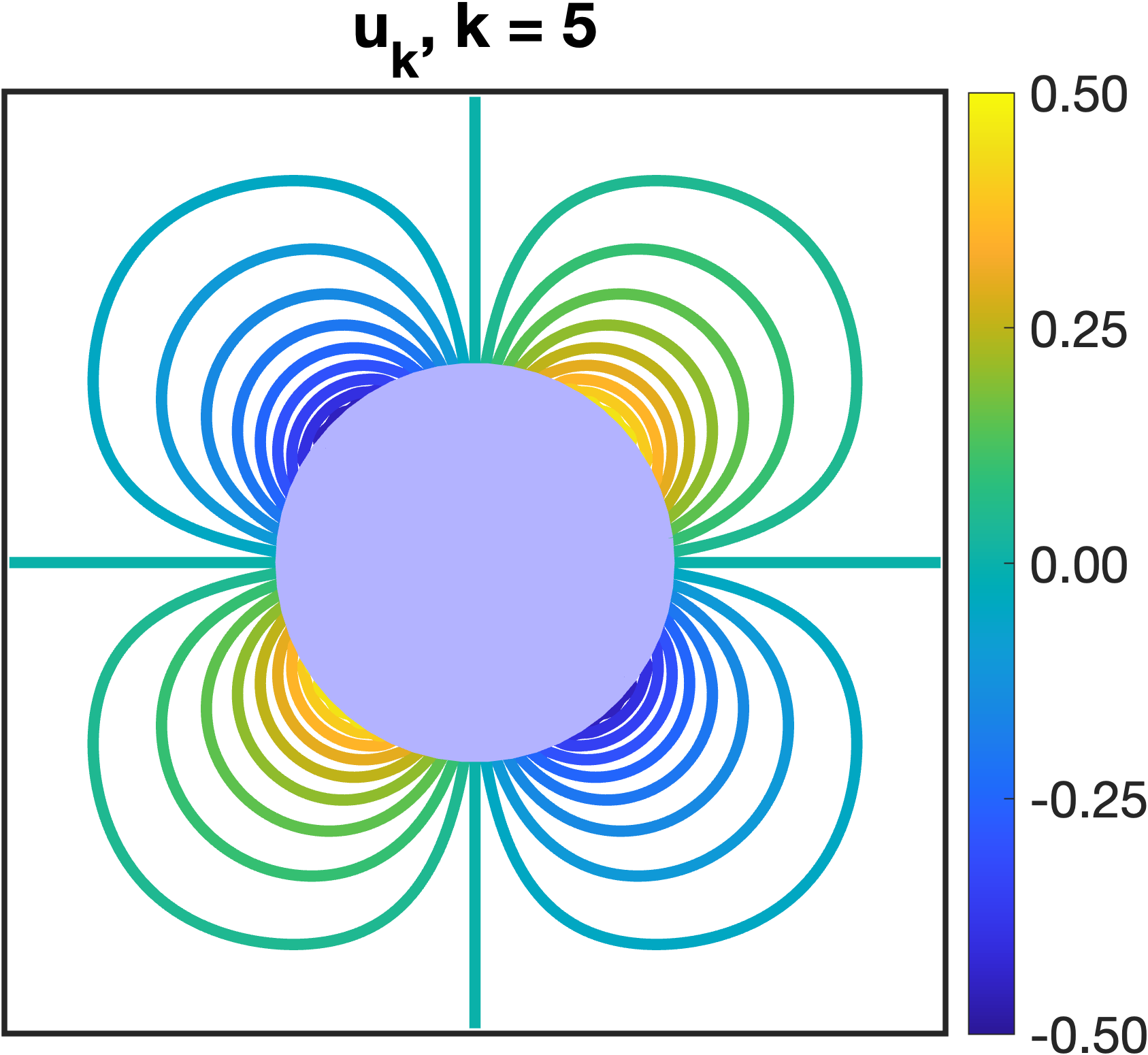}
\includegraphics[width=0.28\textwidth]{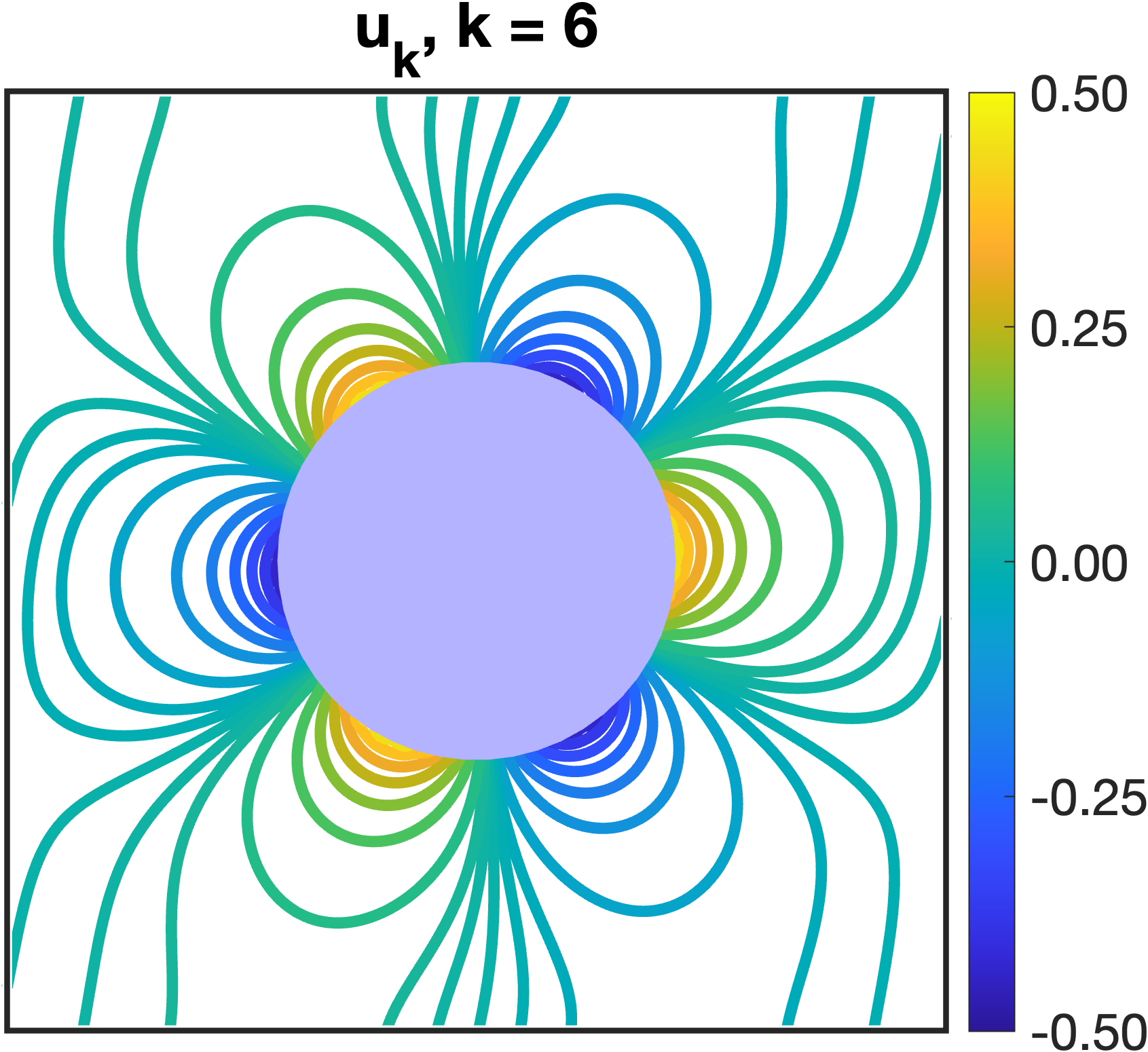}
\includegraphics[width=0.28\textwidth]{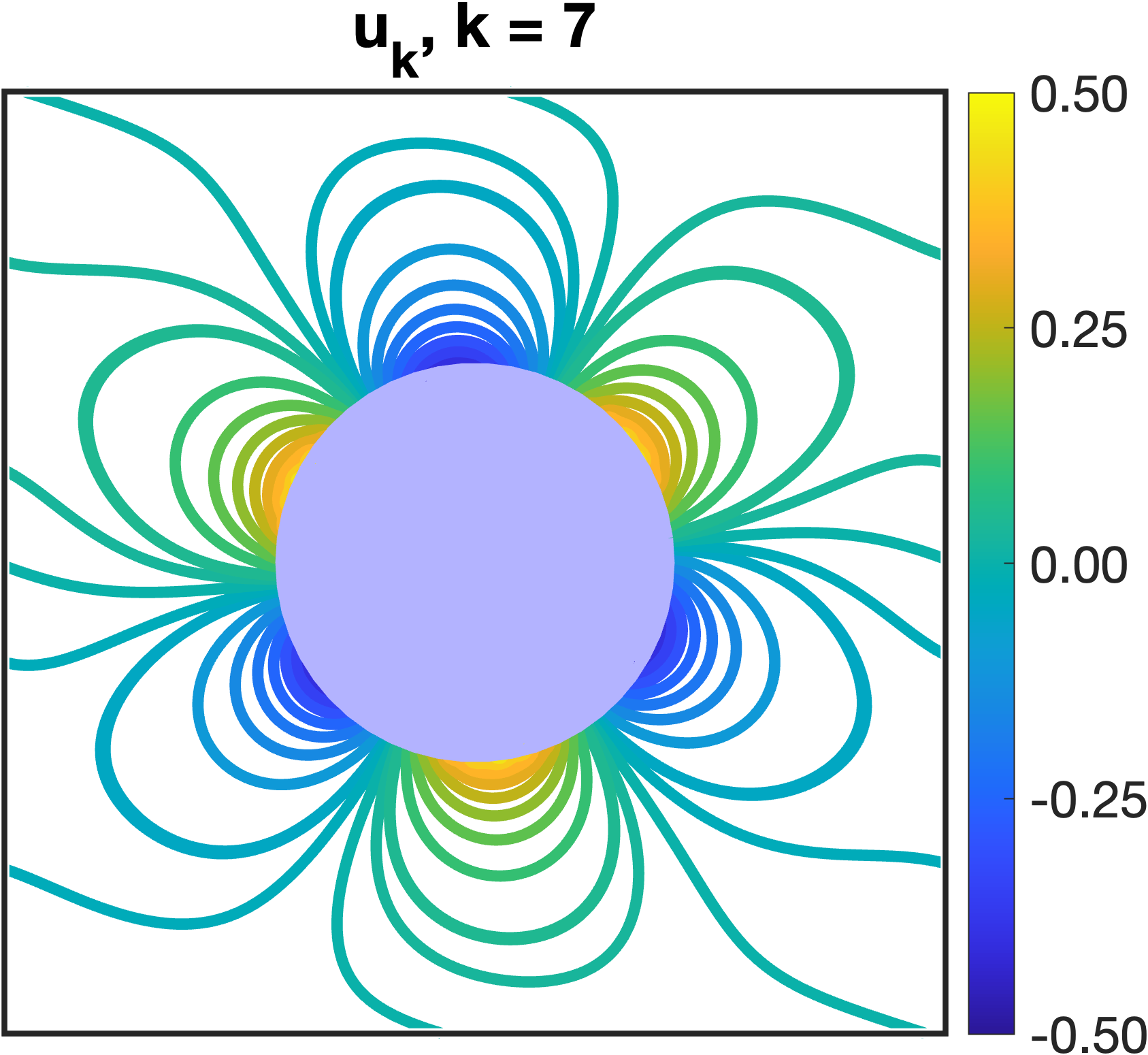} \\ 
\vspace{0.5cm}
\includegraphics[width=0.28\textwidth]{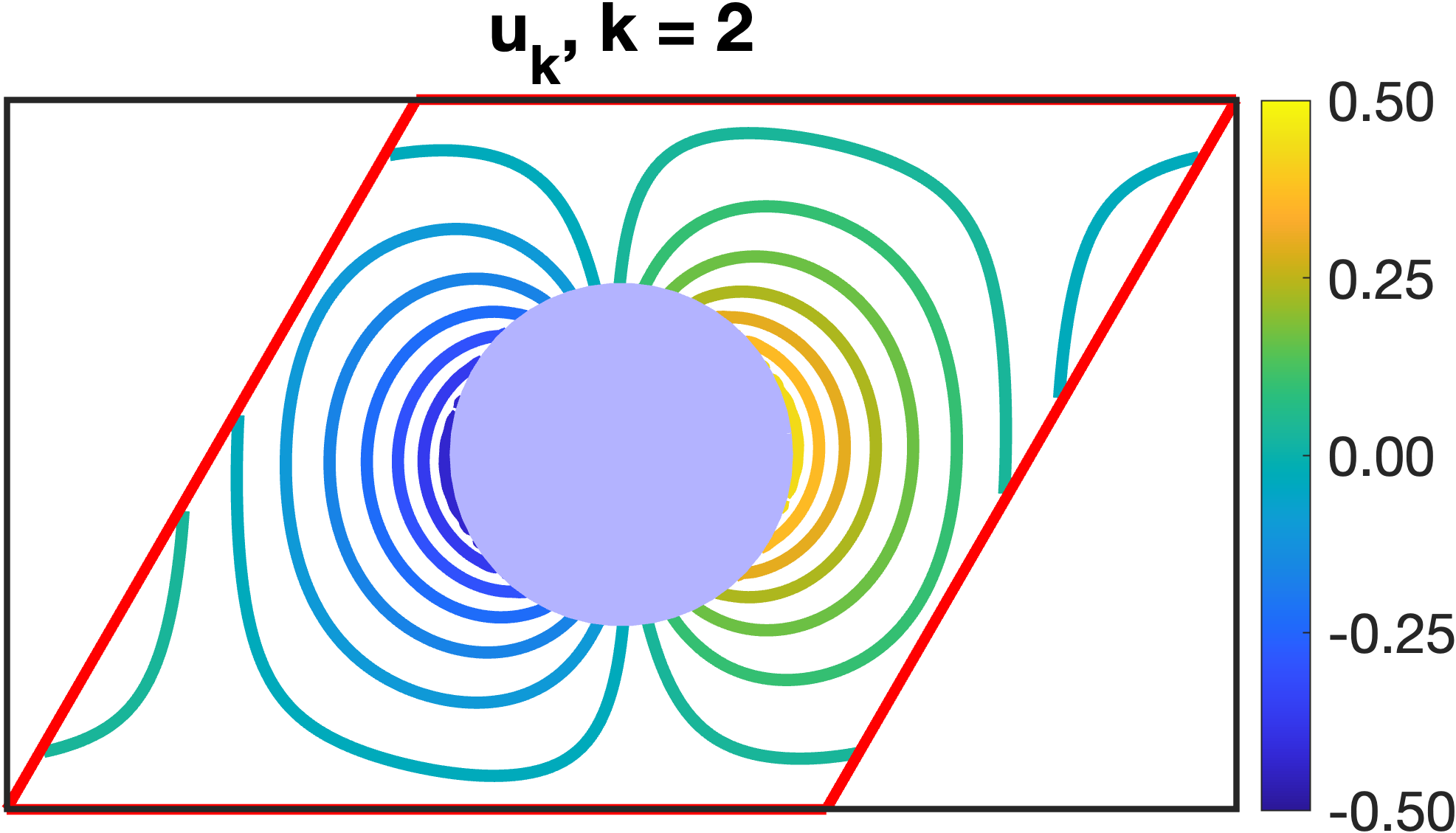}
\includegraphics[width=0.28\textwidth]{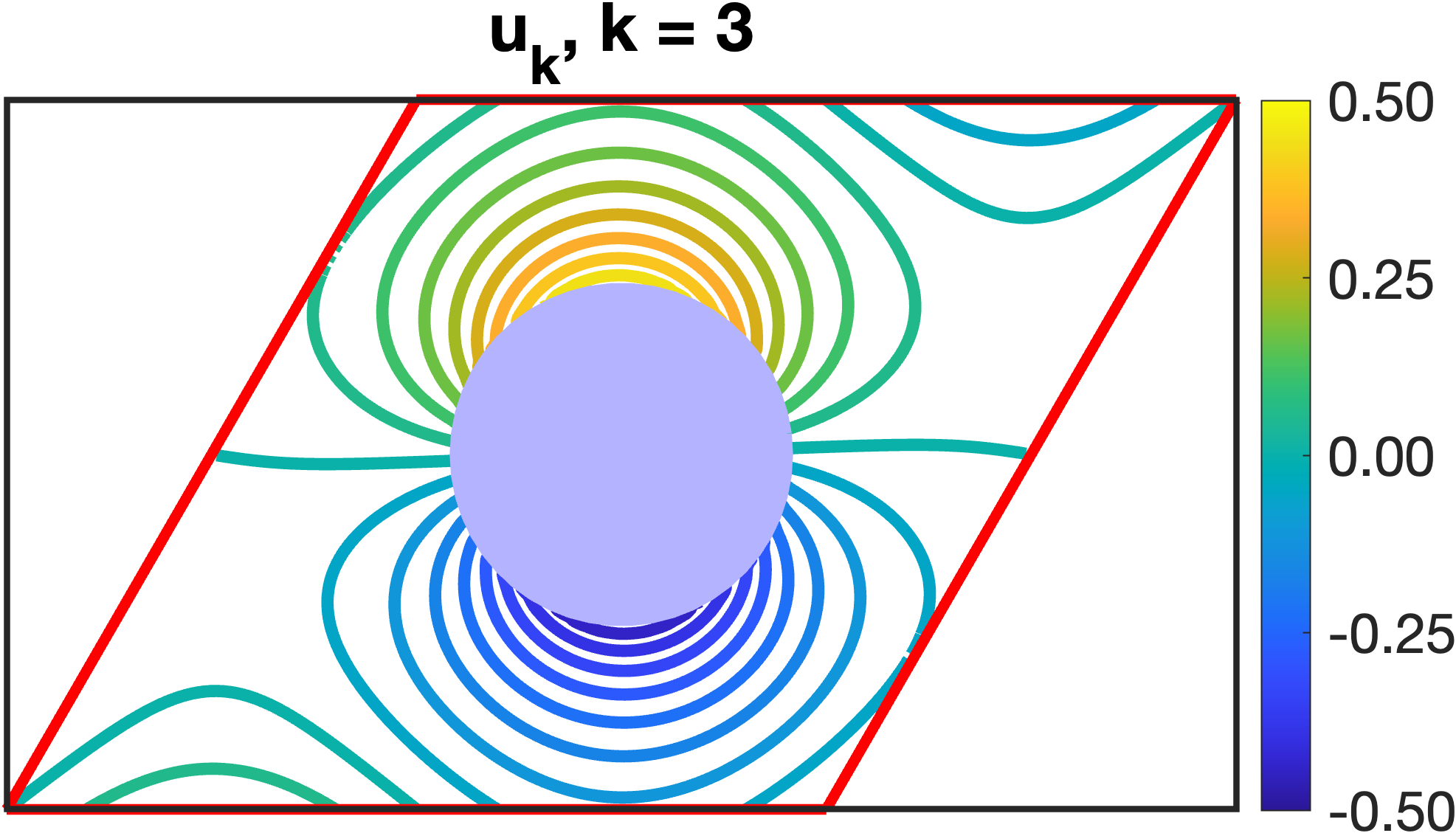}
\includegraphics[width=0.28\textwidth]{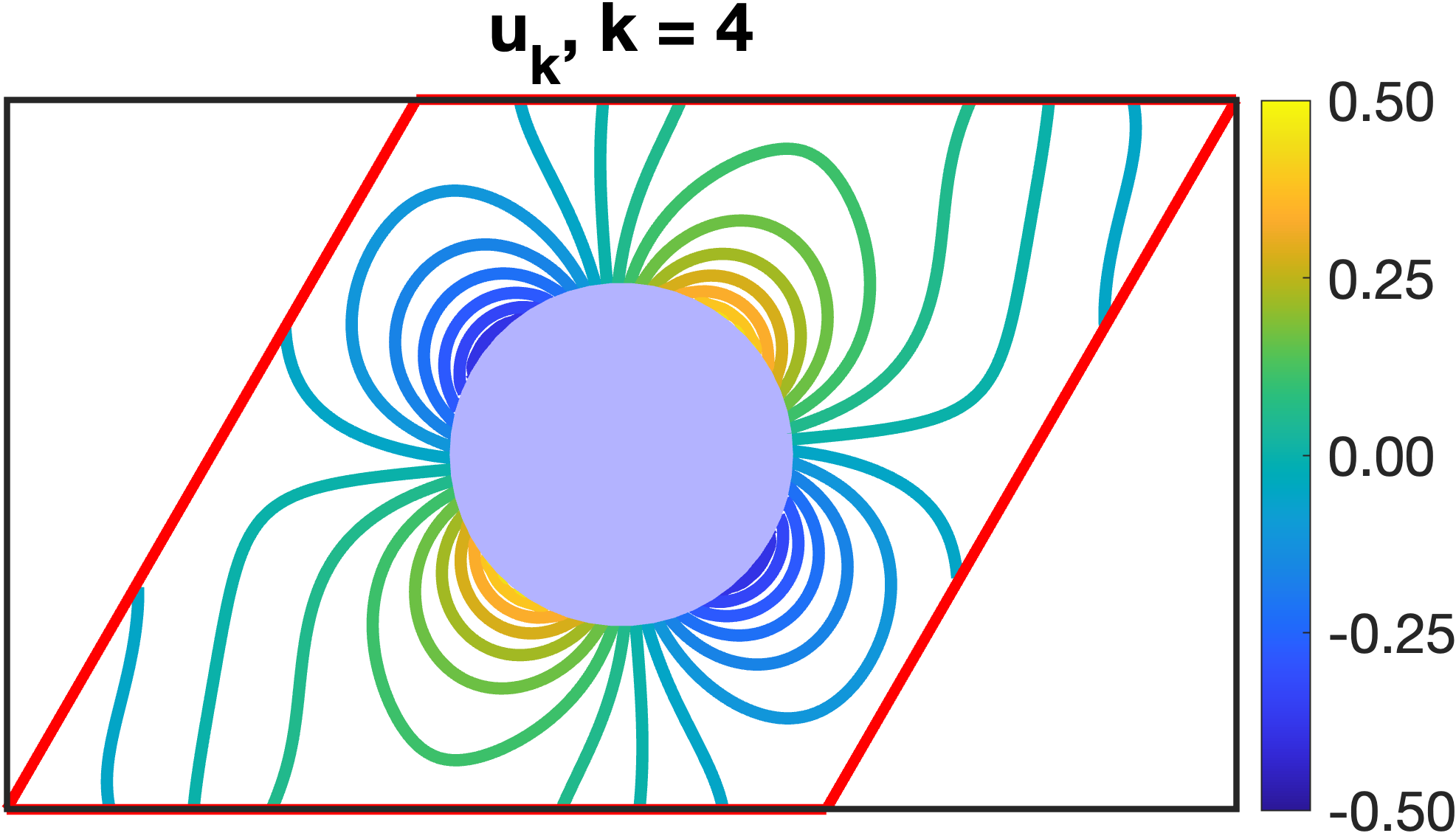} \\ 
\includegraphics[width=0.28\textwidth]{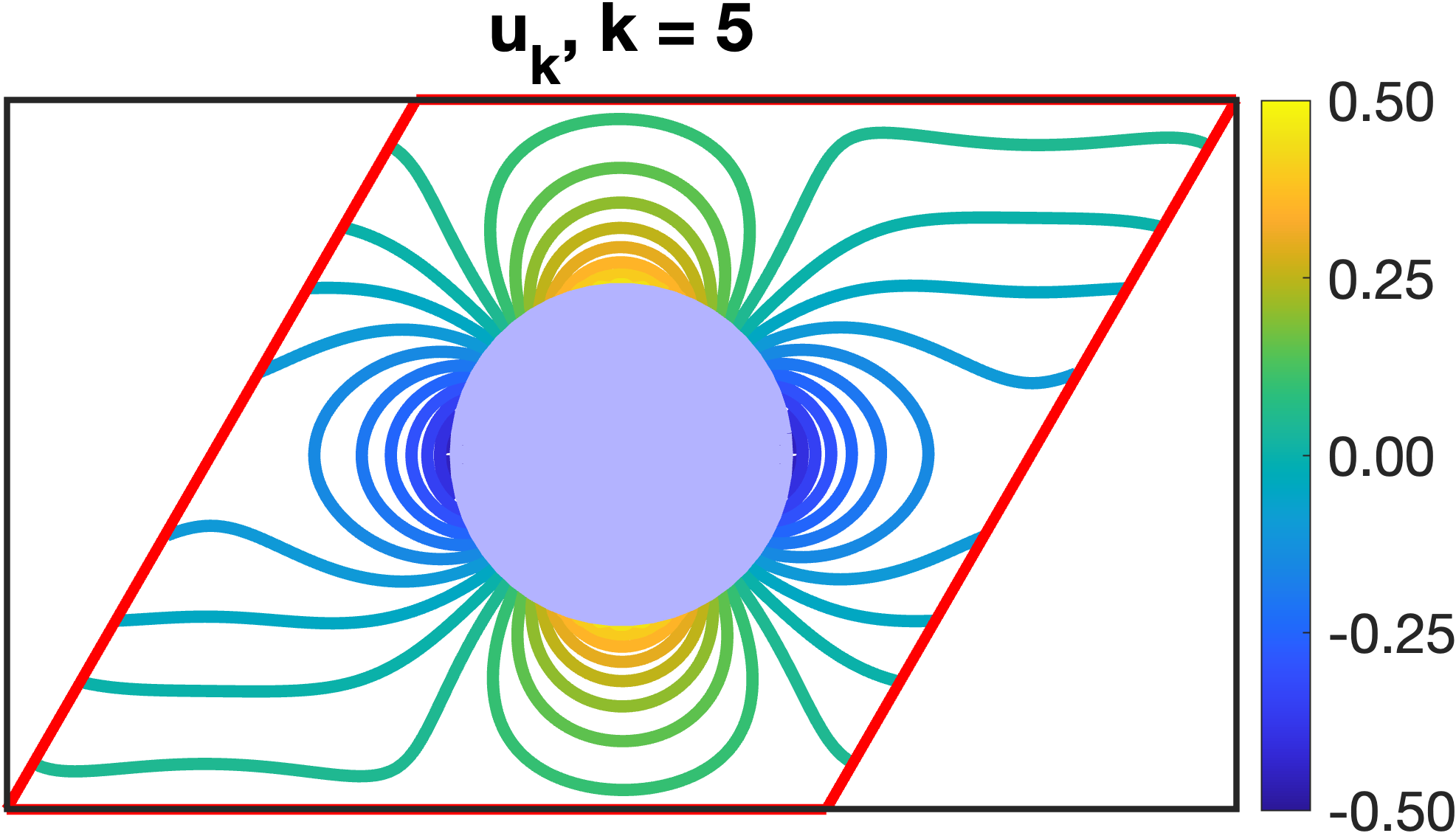}
\includegraphics[width=0.28\textwidth]{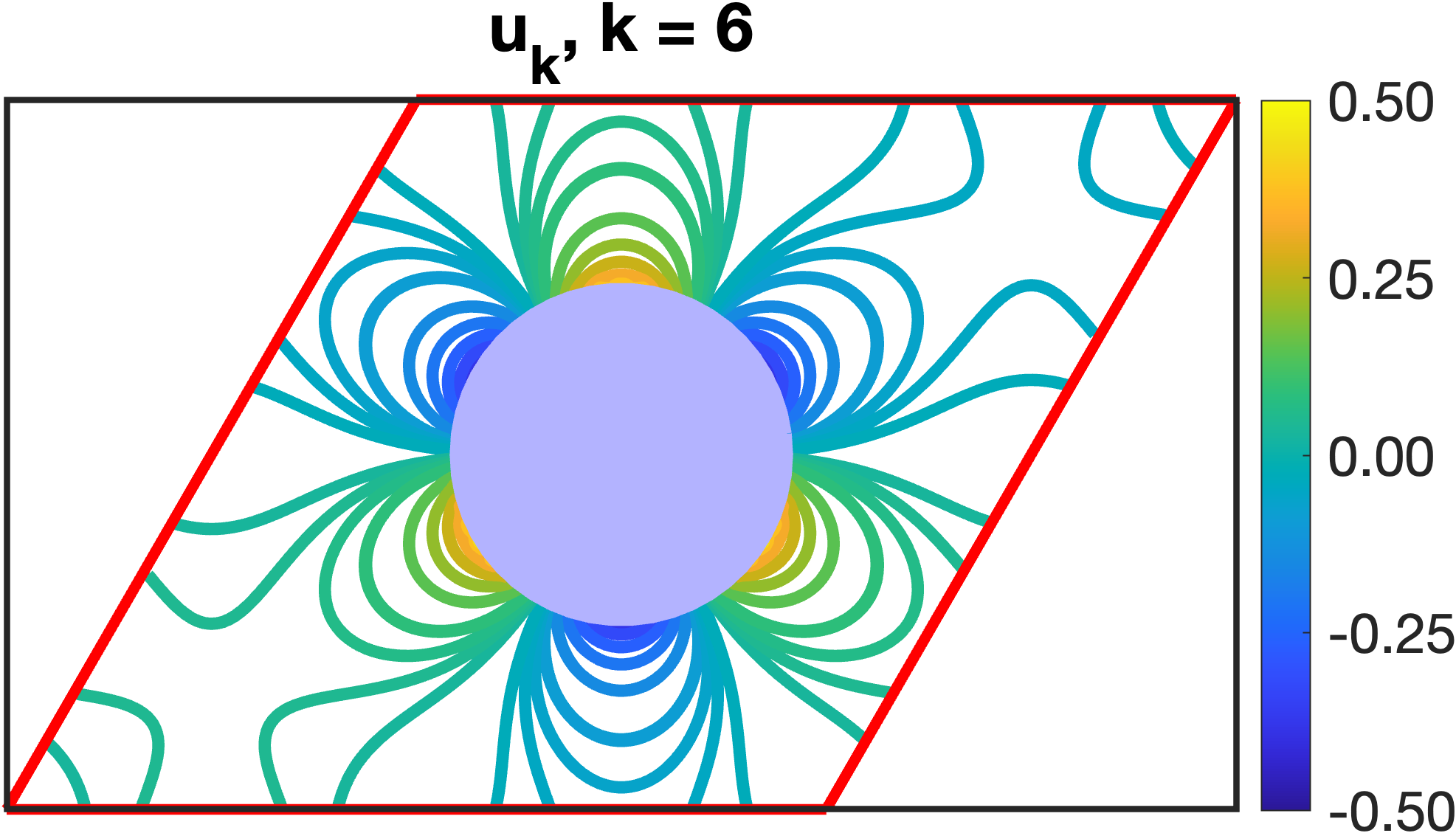}
\includegraphics[width=0.28\textwidth]{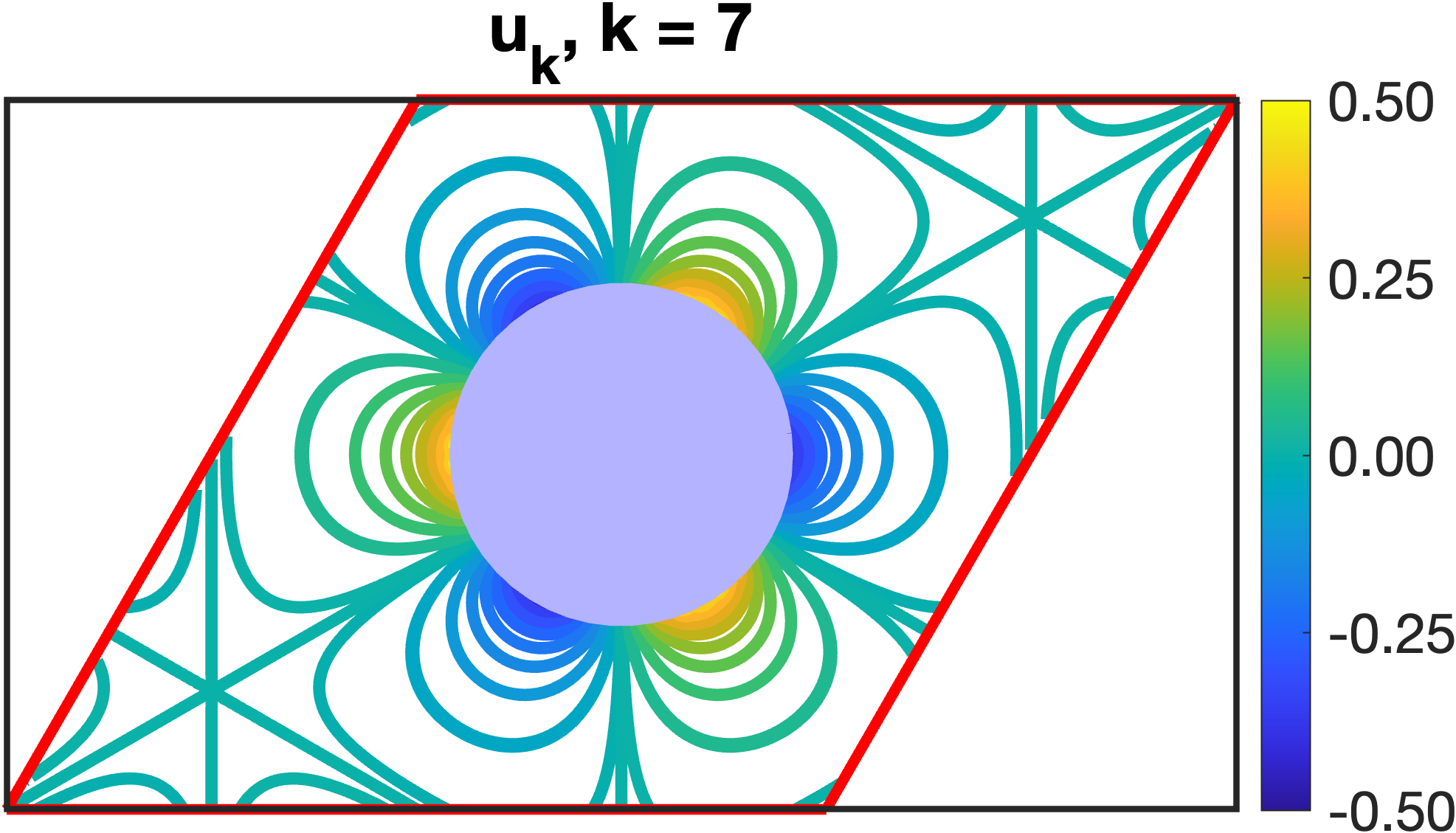}
\end{center}
\caption{Steklov eigenfunctions $u_k$ of the square torus \textbf{ (top two rows)} and equilateral torus \textbf{(bottom two rows)} with a single hole. See Sec.~\ref{sec: steklovbvp numerical}, Example 4, and Table~\ref{tabl: steklov evp err one hole, square}.}
\label{fig: simp steklov}
\end{figure}

\vspace{3mm} 
\noindent {\bf Example 4. Steklov EVP with a single ($M=1$) hole.}
We solve the Steklov EVP on tori with a single circular hole; the results are shown in Fig.~\ref{fig: simp steklov}.
Let $D = B(a_1, r)$ be the disk centered at $a_1 = 0.5 + 0.5i$ with radius $r = 0.2$.  
This example corresponds to Fig.~4 and Fig.~7 in~\cite{kaoharmonic2023}, where the eigenvalues are reported in Tables~1 and 4 of~\cite{kaoharmonic2023}, with an accuracy of approximately 50 digits. 

Using $N=50$ points to discretize the boundary integral equation~\eqref{e: steklov bie}, we calculate 2nd to 7th eigenpairs $(\sigma_k, u_k)$ for $ k=2,3,\ldots,7$ using the representation in~\eqref{e:uRep, stek}; the eigenfunctions are plotted in Fig.~\ref{fig: simp steklov}. 
In Table~\ref{tabl: steklov evp err one hole, square}, we tabulate the estimated 2nd to 7th eigenvalues and their associated errors. 
The errors are computed as the absolute difference between half of our results and those in~\cite{kaoharmonic2023}; we achieve an accuracy of 13 to 16 digits.

\begin{table}[t!]
\centering
{\small
\begin{tabular}{|l|l|l|}
\hline
    $k$ & Estimated $\sigma_k$  & Error \\ \hline
    $2$ & $3.217375$ & $7.550 \times 10^{-15}$ \\ \hline
    $3$ & $3.217375$ & $1.510 \times 10^{-14}$ \\ \hline
    $4$ & $4.850995$ & $1.004 \times 10^{-13}$ \\ \hline
    $5$ & $5.153581$ & $1.155 \times 10^{-13}$ \\ \hline
    $6$ & $7.503050$ & $3.508 \times 10^{-13}$ \\ \hline
    $7$ & $7.503050$ & $3.366 \times 10^{-13}$ \\ \hline
\end{tabular}
\quad 
\begin{tabular}{|l|l|}
\hline
     Estimated $\sigma_k$  & Error  \\ \hline
   3.348656 & $4.440 \times 10^{-15}$ \\ \hline
3.348656 & $1.780 \times 10^{-15}$ \\ \hline
4.999789 & $1.780 \times 10^{-15}$ \\ \hline
4.999789 & $7.100 \times 10^{-15}$ \\ \hline
7.443925 & $1.865 \times 10^{-14}$ \\ \hline
7.556497 & $3.550 \times 10^{-15}$ \\ \hline
\end{tabular}\vspace{0.1em}}
    \caption{For the Steklov EVP on the square torus {\bf~(left)} and 
    equilateral torus  {\bf (right)} with a single hole, we tabulate the estimated eigenvalues $\sigma_k$ and the errors when compared to the values in \cite[App.~B]{kaoharmonic2023}. See Sec.~\ref{sec: steklovbvp numerical}, Example 4 and Fig.~\ref{fig: simp steklov}.}
    \label{tabl: steklov evp err one hole, square}
\end{table}

\begin{figure}[t!]
\begin{center}
\includegraphics[width=0.300\textwidth]{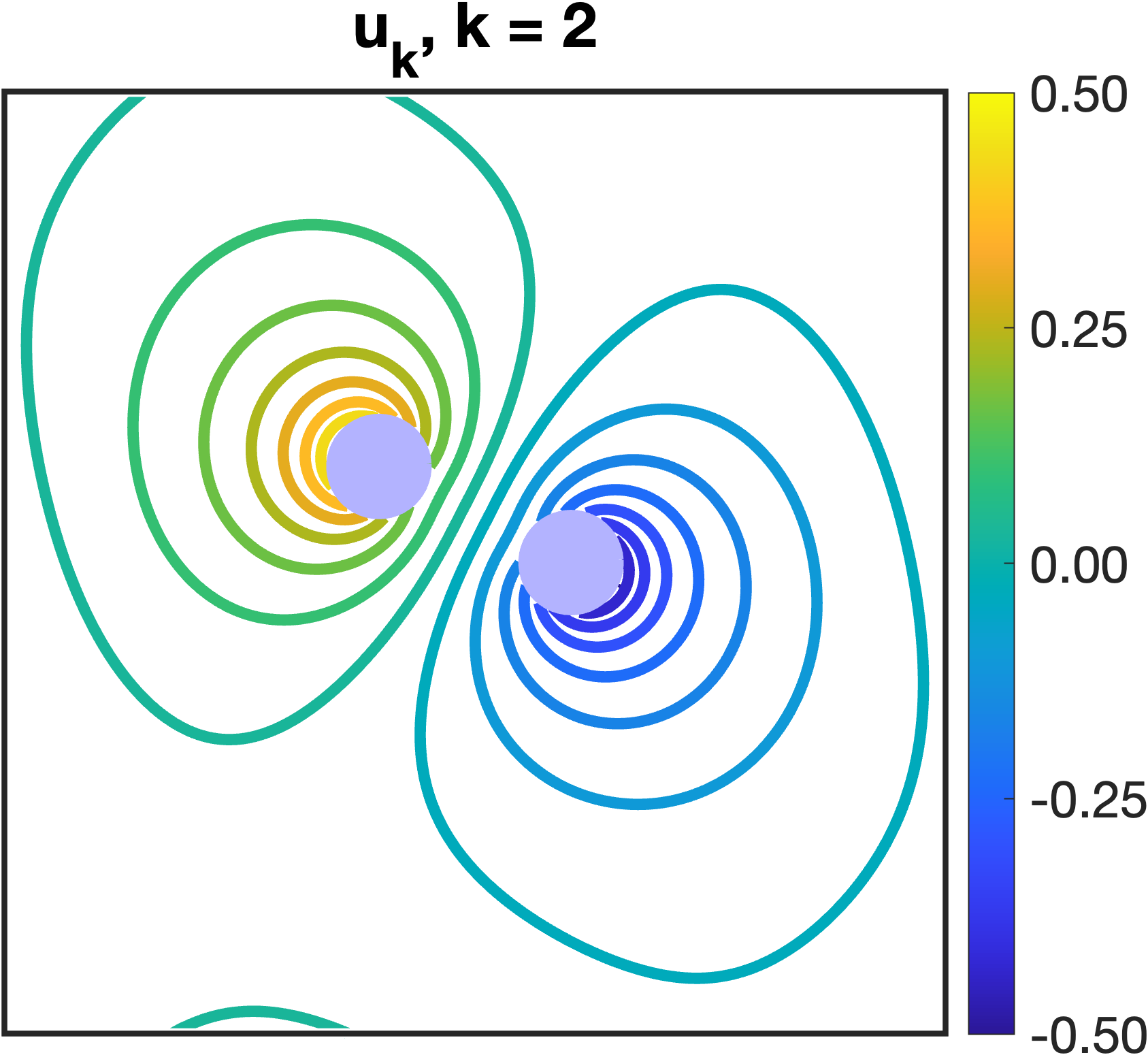}
\includegraphics[width=0.300\textwidth]{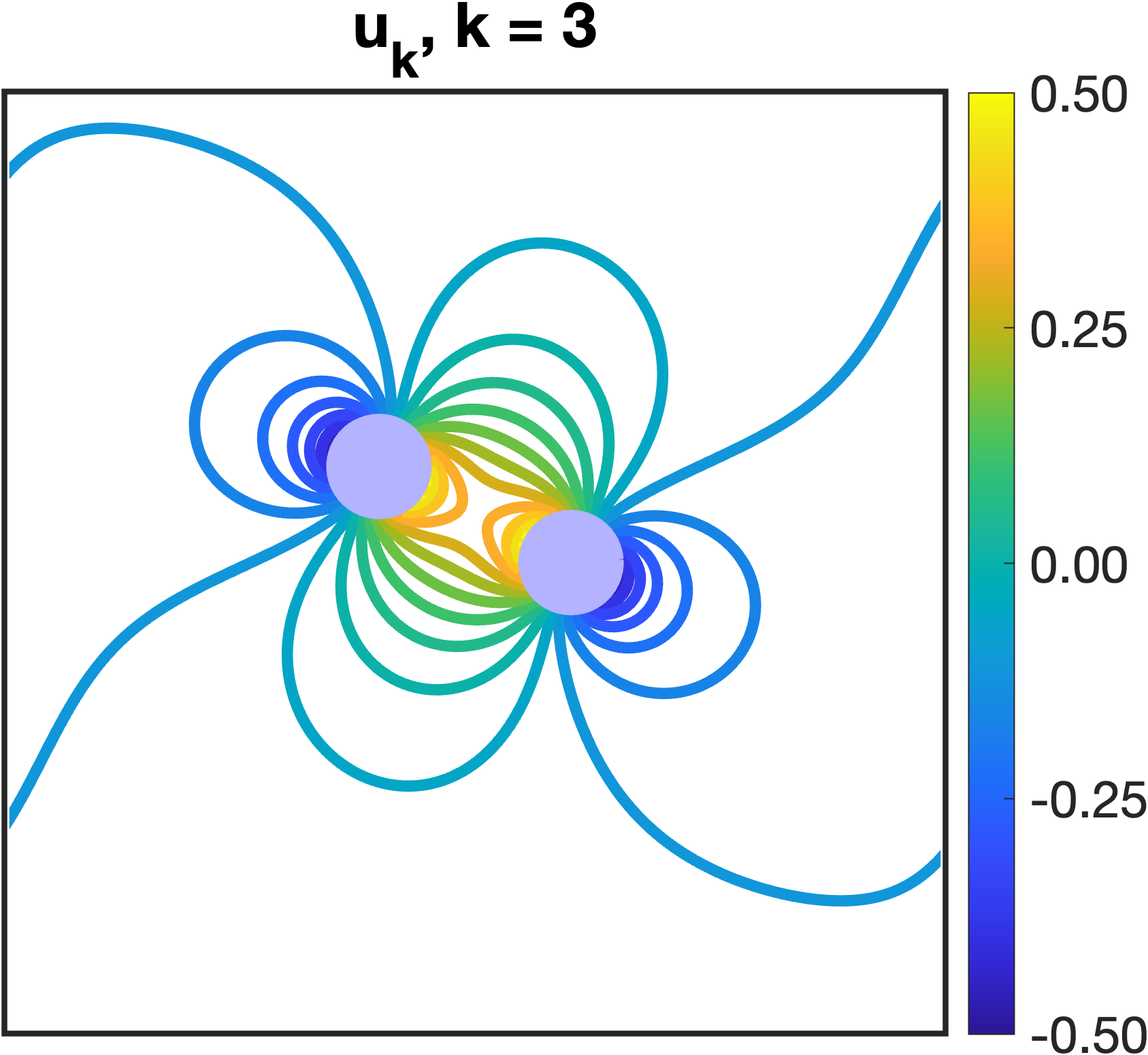}
\includegraphics[width=0.300\textwidth]{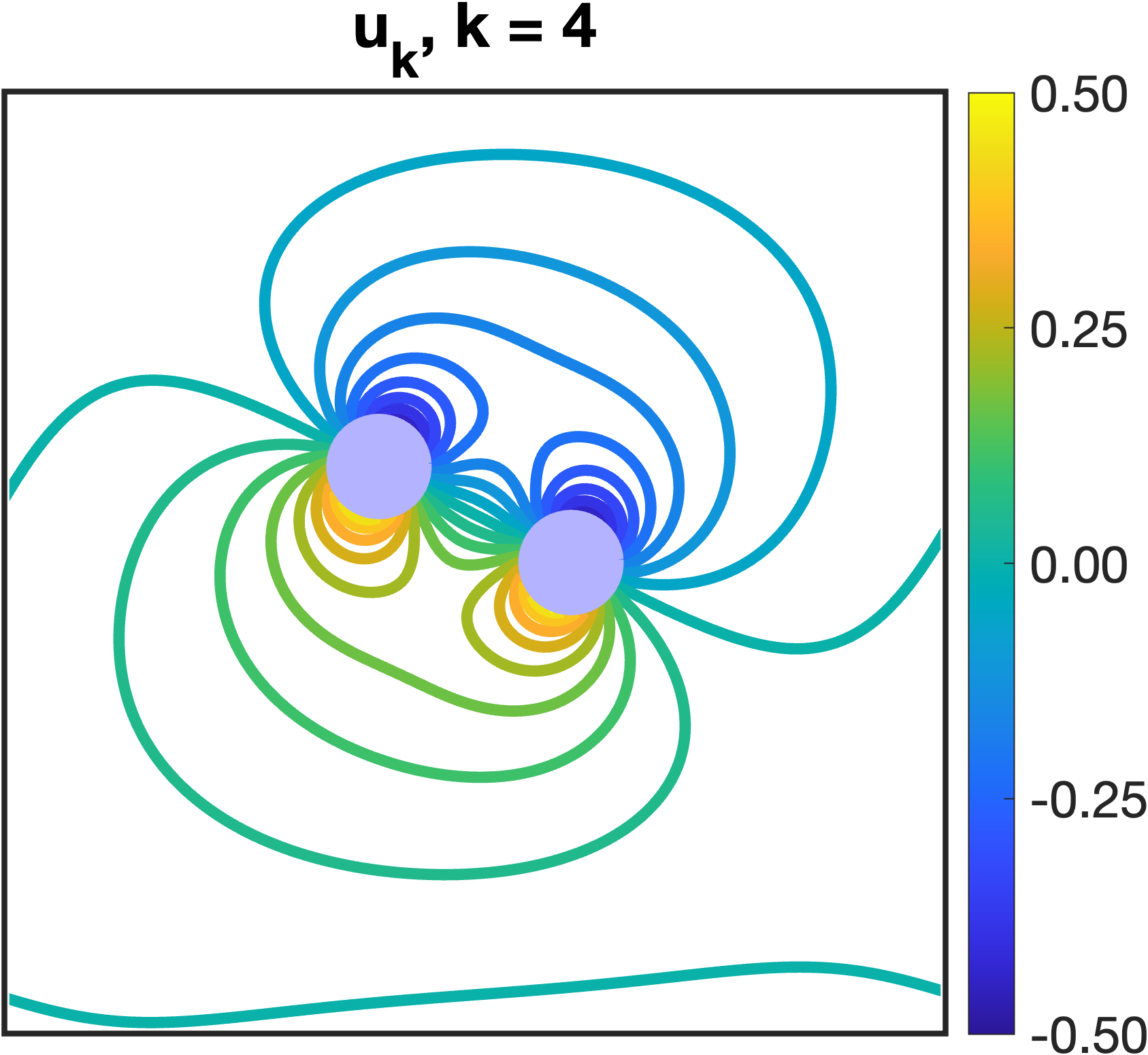}
\includegraphics[width=0.300\textwidth]{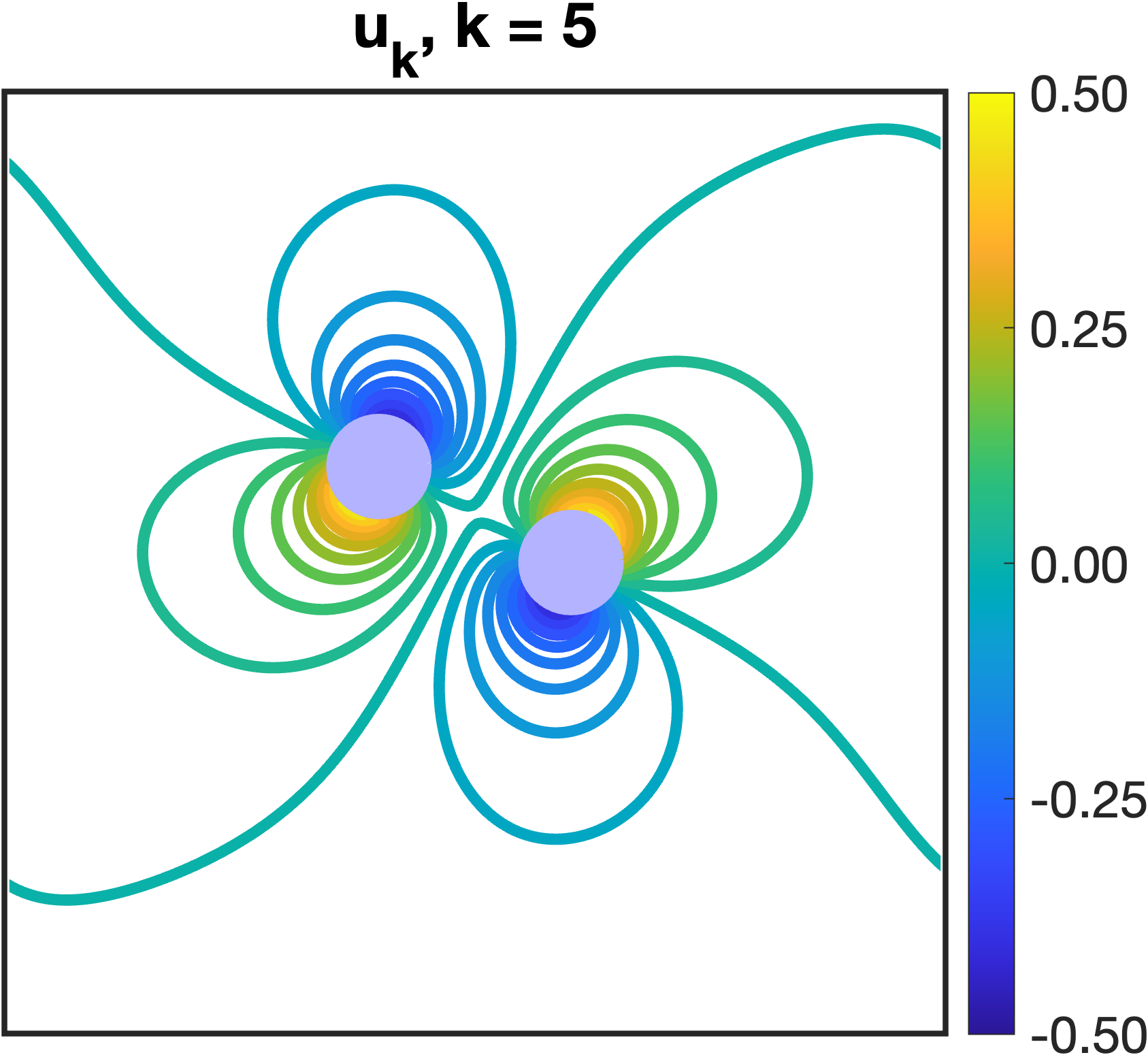}
\includegraphics[width=0.300\textwidth]{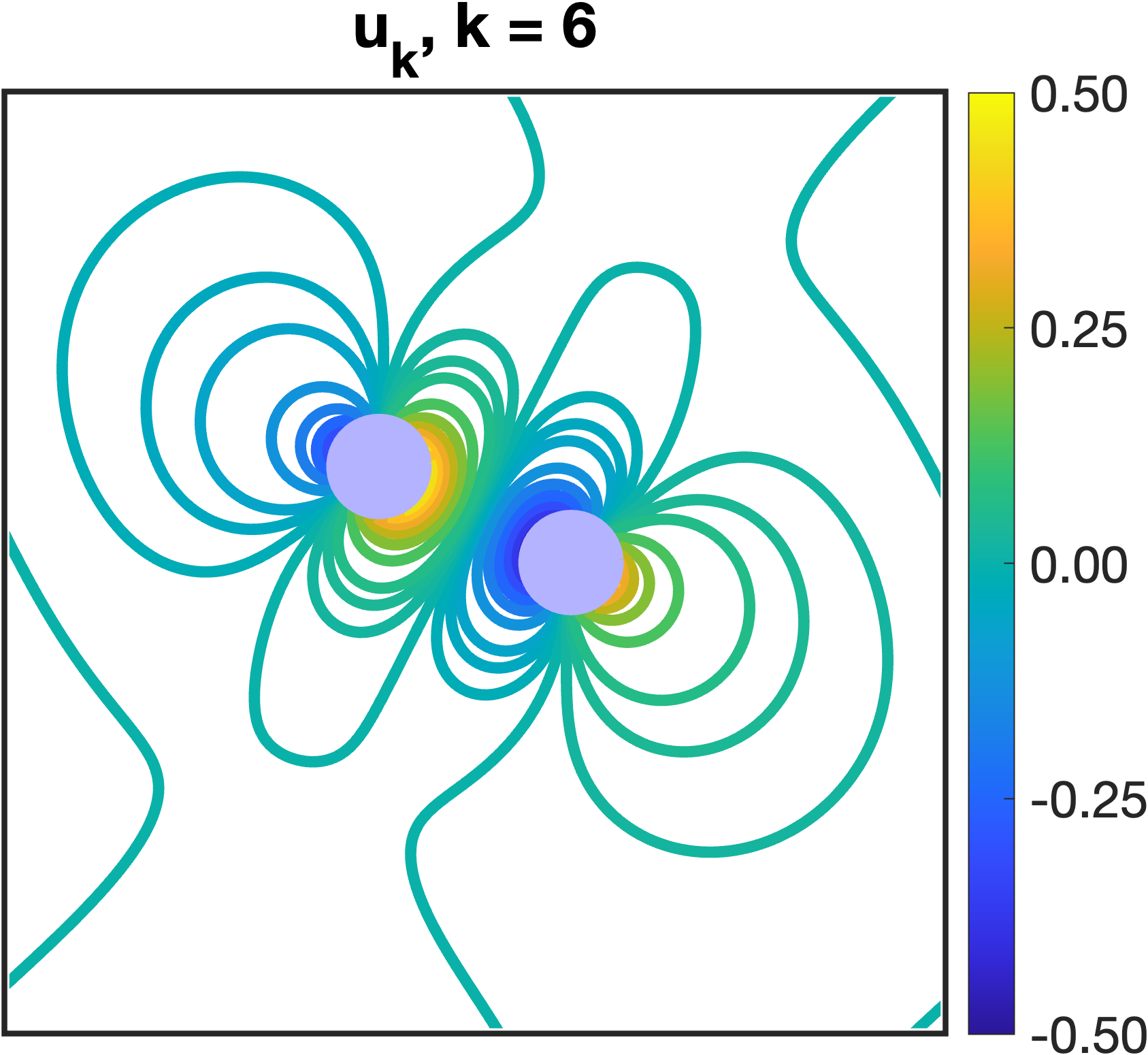}
\includegraphics[width=0.300\textwidth]{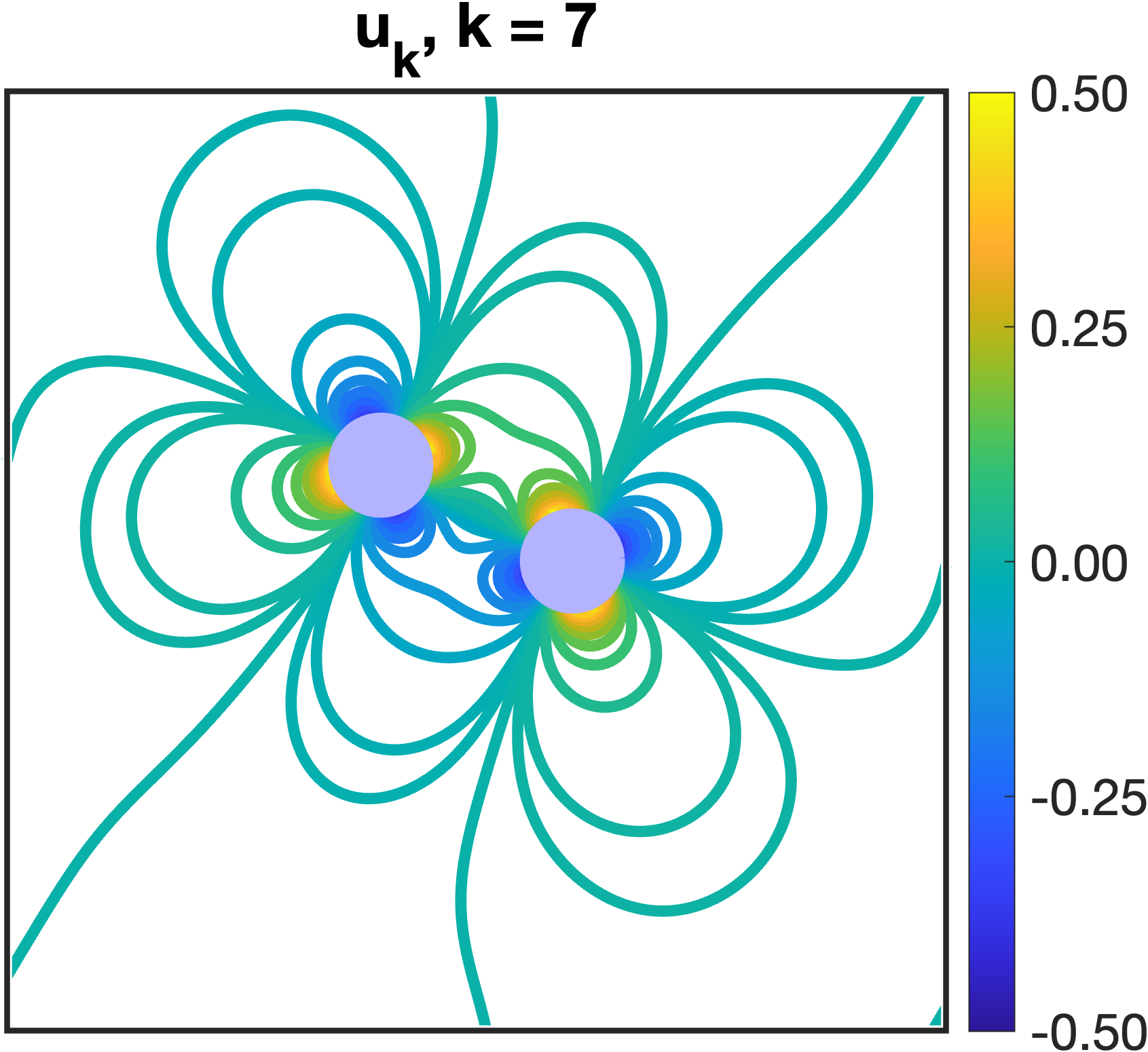}\\ 
\vspace{0.5cm}
\includegraphics[width=0.300\textwidth]{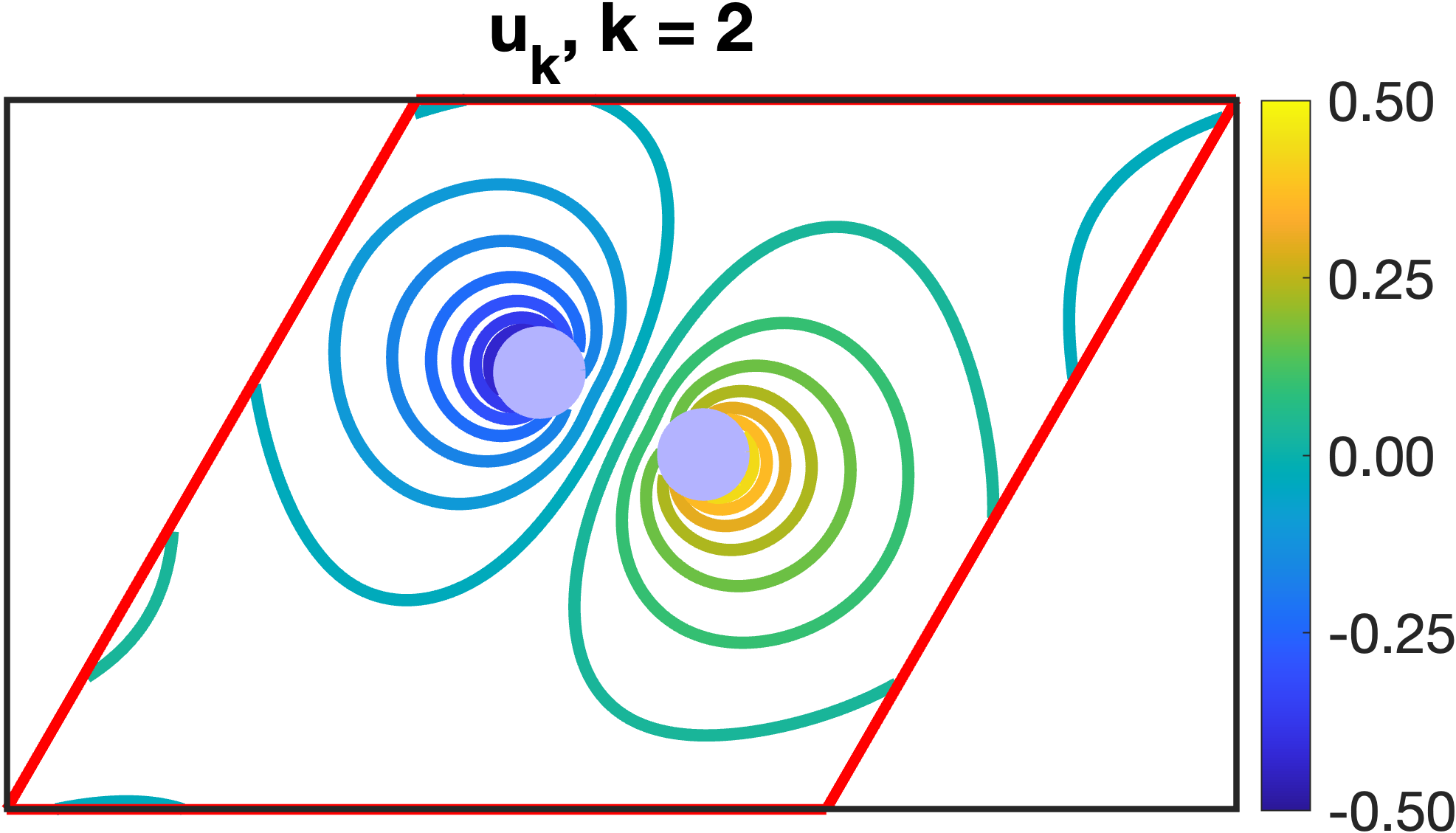}
\includegraphics[width=0.300\textwidth]{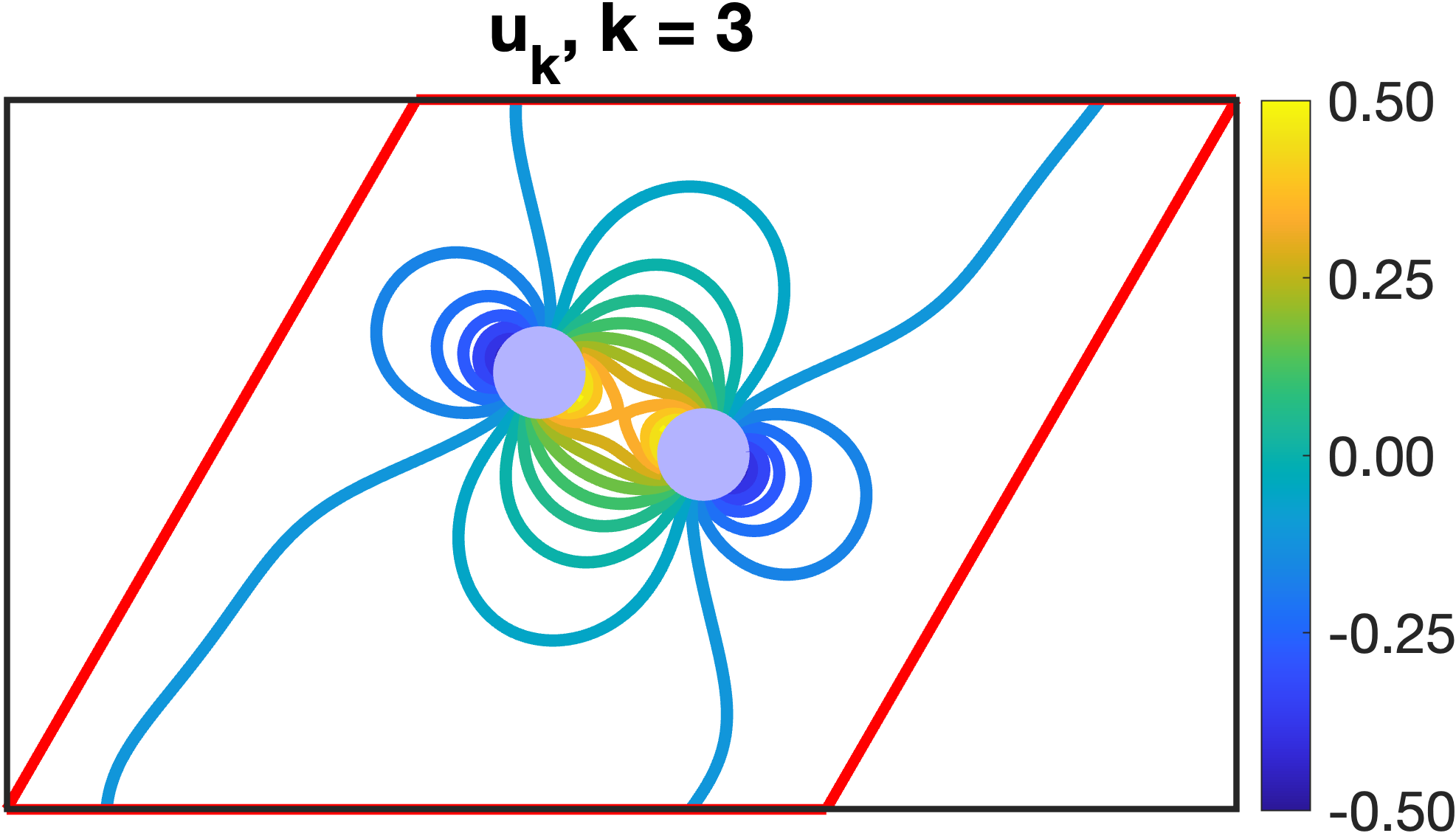}
\includegraphics[width=0.300\textwidth]{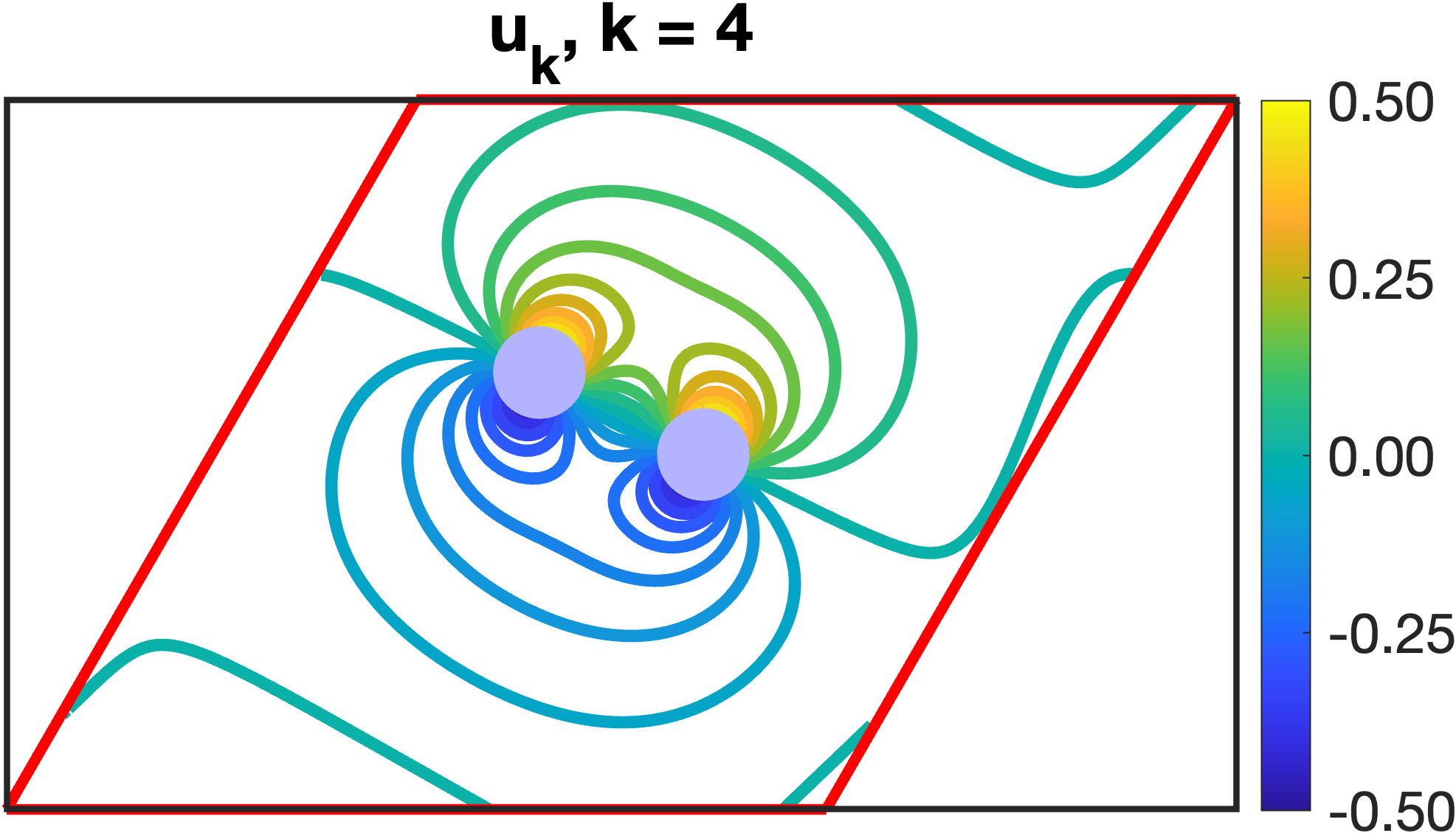}
\includegraphics[width=0.300\textwidth]{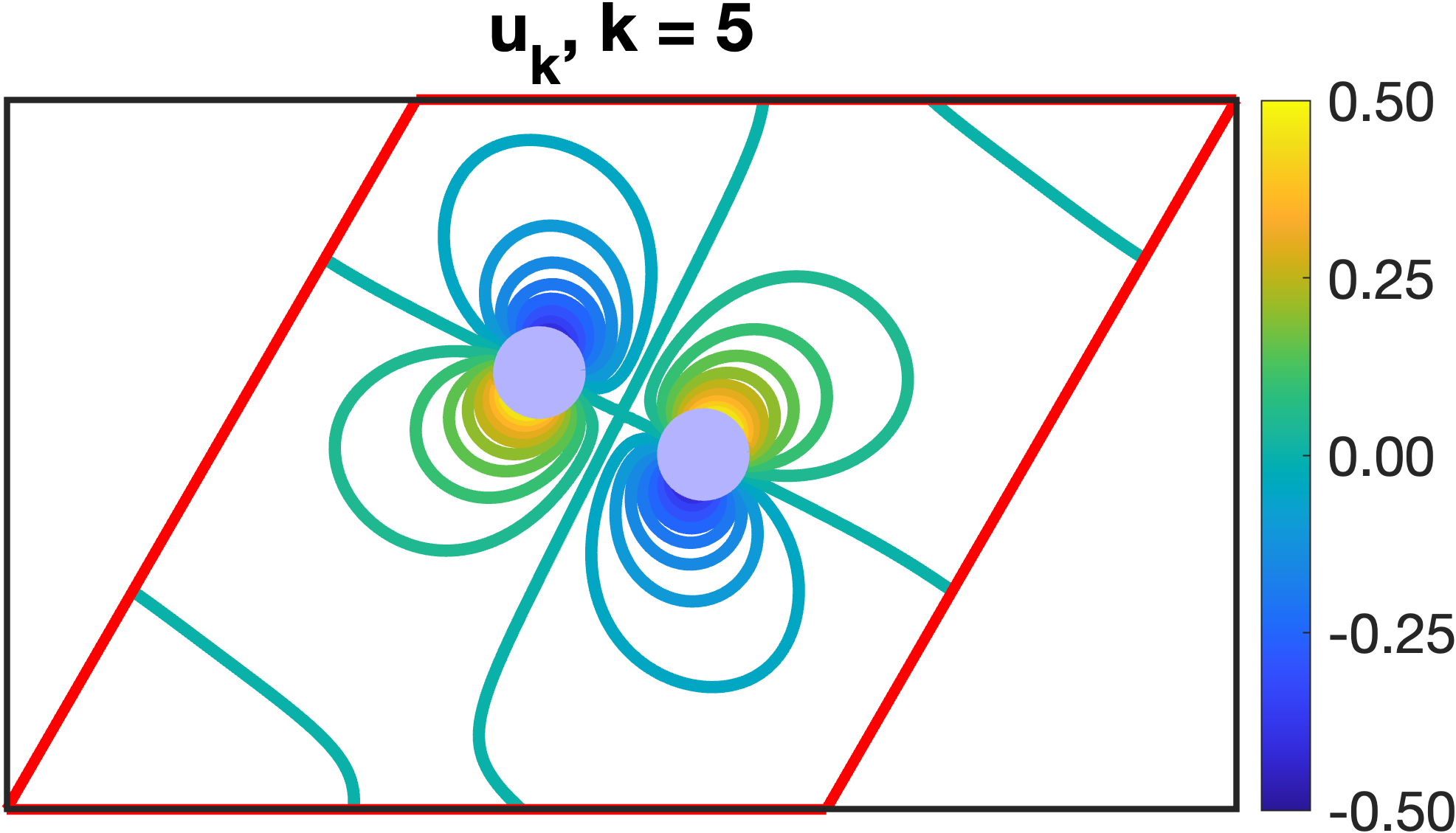}
\includegraphics[width=0.300\textwidth]{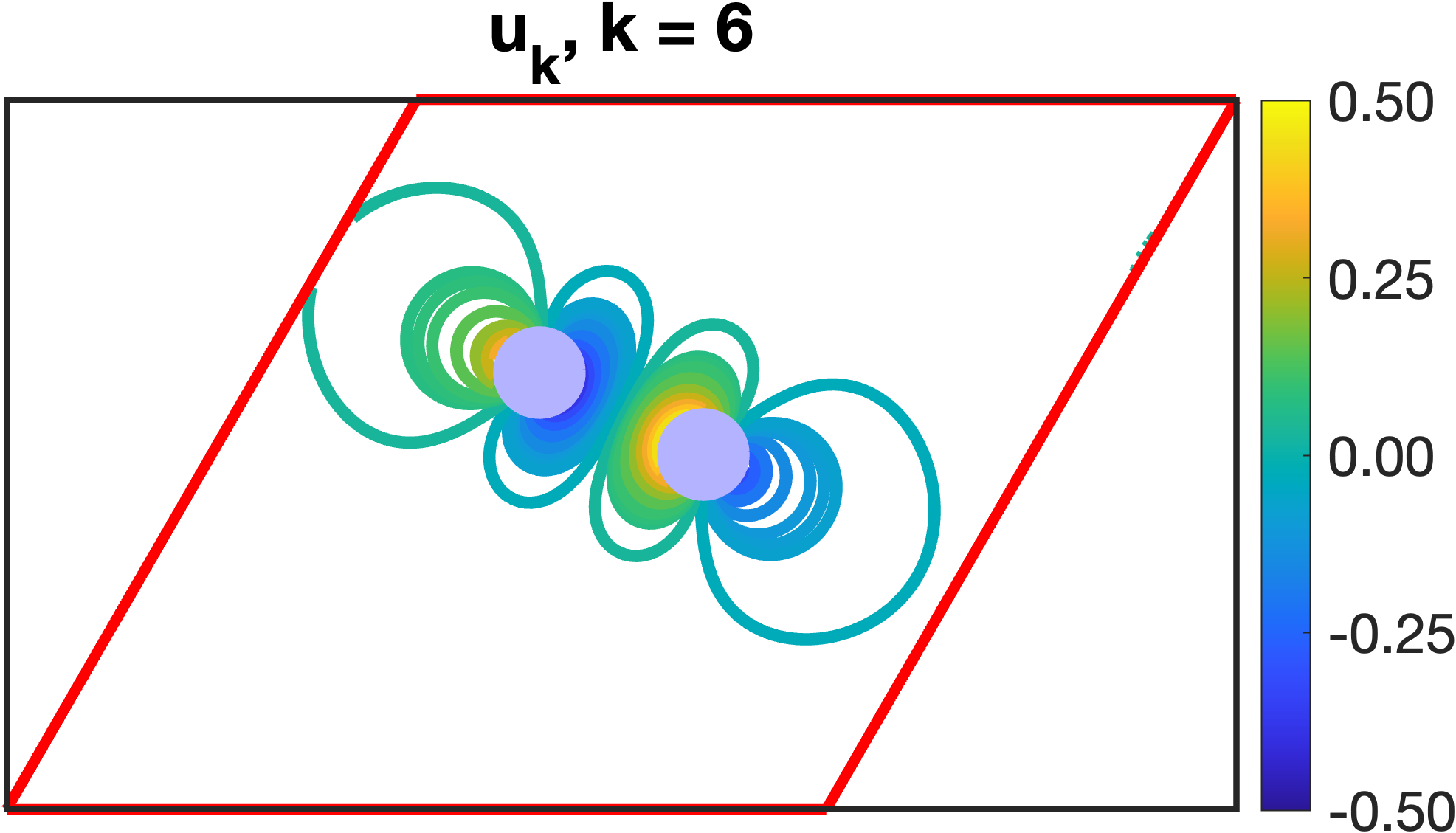}
\includegraphics[width=0.300\textwidth]{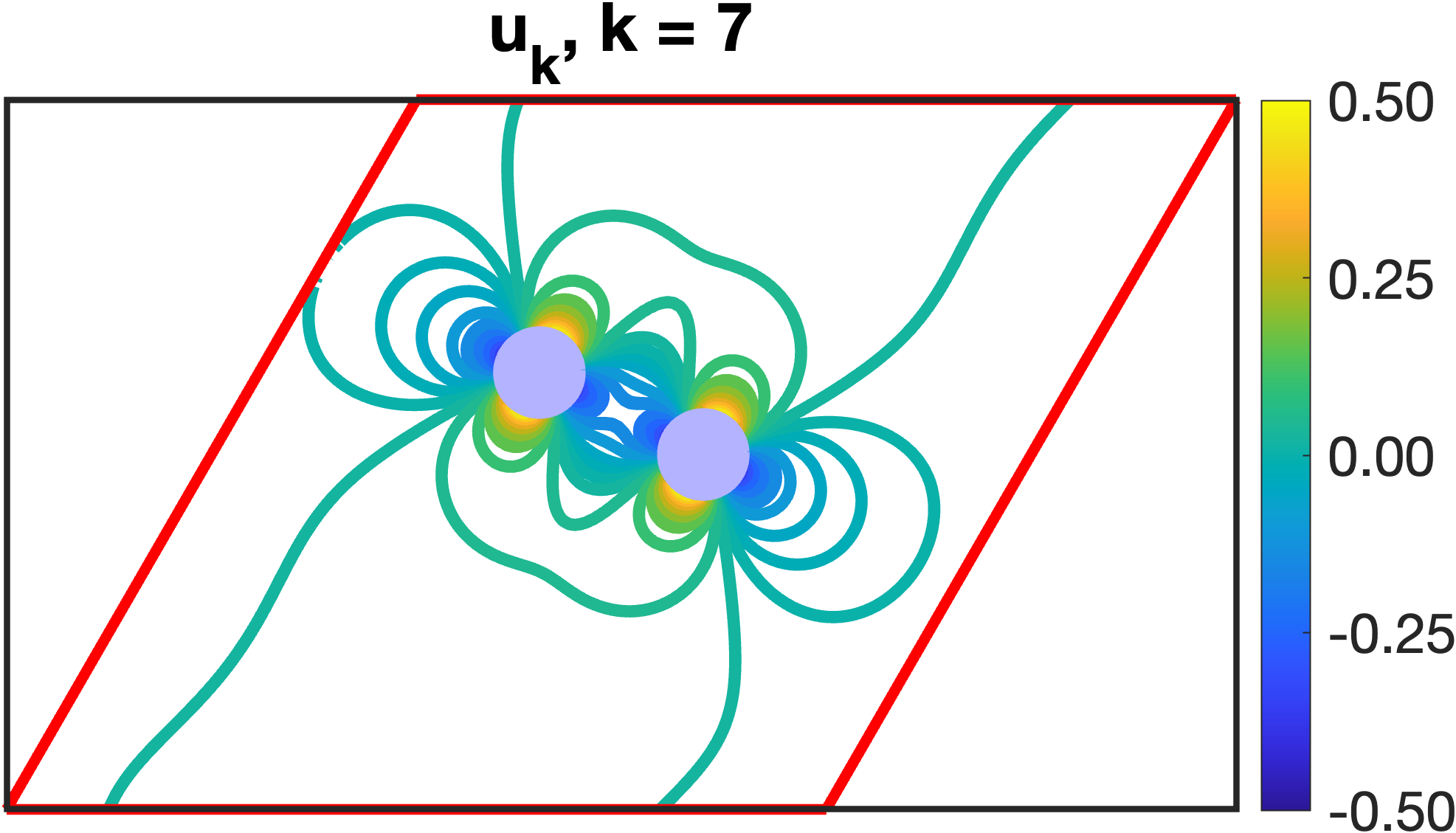}
\end{center}
\caption{Steklov eigenfunctions $u_k$ of the square torus \textbf{ (top two rows)} and equilateral torus \textbf{(bottom two rows)} with $M=2$ holes. See Sec.~\ref{sec: steklovbvp numerical}, Example 5 and Table~\ref{tabl: steklov mult errs}.}
\label{fig: mult steklov}
\end{figure}

\vspace{3mm} 
\noindent {\bf Example 5: Steklov EVP with $M=2$ holes.}
We consider the Steklov EVP on tori with $M=2$ circular holes; the results are shown in Fig.~\ref{fig: mult steklov}.
Let $D = \bigcup_{i=1}^M B(a_i, r)$, the union of disks centered at $a_1 = 0.6+0.5i$ and $a_2 = 0.4+0.6i$ with radii $r = 0.05$. 

This example corresponds to Fig. 5 and Fig. 8 of~\cite{kaoharmonic2023}, where the eigenvalues are reported in Tables 2 and 5 of~\cite{kaoharmonic2023}, which are accurate to approximately 50 digits.
We use $N_i=50$ ($N=100$ total) points to discretize the boundary integral equation~\eqref{eqn: steklov}. 
We calculate eigenpairs $(\sigma_k, u_k)$ for $ k=2,3,\ldots,7$ using the representation in~\eqref{e:uRep, stek}. 
We tabulate the 2nd to 7th eigenvalues $\sigma_k$ in Table~\ref{tabl: steklov mult errs} and plot their corresponding eigenfunctions $u_k$ given in~\eqref{e:uRep, stek} in Fig.~\ref{fig: mult steklov}. We achieve an accuracy of 14 to 16 digits.

We consider the flux of the \( k \)-th eigenfunction through the first hole. Since the \( k \)-th eigenfunction is determined up to a constant multiple, we normalize it with respect to the $L^2(\partial \Omega)$ norm and report the absolute value of the flux, 
$
 |A_1(k)| 
 = \frac{1 }{\| u_k \|_{L^2(\partial \Omega)}}  
 \left| \int_{\partial D_1} \partial_\nu u_k(\xi)\, |d\xi| \right|. 
$
In Table~\ref{tabl: steklov mult errs}, we report the absolute value of the flux since there are only two holes and recall that \( A_2(k) = -A_1(k) \).

 \begin{table}[!t]
 \centering
{\small
\begin{tabular}{|c|c|c|}
    \hline
    $k$ & Error in $\sigma_k$ & $|A_1(k)|$ \\ \hline
    $2$ & $9.800 \times 10^{-15}$ & $4.796$ \\ \hline
$3$ & $7.100 \times 10^{-15}$ & $2.216 \times 10^{-14}$ \\ \hline
$4$ & $7.100 \times 10^{-15}$ & $0.05758$ \\ \hline
$5$ & $7.100 \times 10^{-15}$ & $1.267 \times 10^{-14}$ \\ \hline
$6$ & $3.910 \times 10^{-14}$ & $3.470$ \\ \hline
$7$ & $1.070 \times 10^{-14}$ & $0.01060$ \\ \hline
\end{tabular}
\quad 
\begin{tabular}{|c|c|}
    \hline
    Error in $\sigma_k$ & $|A_1(k)|$ \\  \hline
    $8.000 \times 10^{-15}$ & $4.863$ \\ \hline
    $1.780 \times 10^{-14}$ & $1.412 \times 10^{-14}$ \\ \hline
    $3.200 \times 10^{-14}$ & $0.002723$ \\ \hline
    $1.240 \times 10^{-14}$ & $1.363 \times 10^{-14}$ \\ \hline
    $5.300 \times 10^{-15}$ & $3.434$ \\ \hline
    $1.780 \times 10^{-14}$ & $2.309 \times 10^{-14}$ \\ \hline
\end{tabular}
\vspace{0.5em}
}
\caption{
For the Steklov EVP on a square torus {\bf~(left)} and equilateral torus {\bf~(right)} with
    $M=2$ holes, we tabulate the error for  $\sigma_k$ for $k=2,\ldots,7$ when compared to the values in \cite[App.~B]{kaoharmonic2023}, along with the fluxes $|A_1| = |A_2|$ of the associated eigenfunctions. See Sec.~\ref{sec: steklovbvp numerical}, Example 5 and Fig.~\ref{fig: mult steklov}.}
\label{tabl: steklov mult errs}
\end{table}

\begin{figure}[t!]
\begin{center}
\includegraphics[width=.38\textwidth]{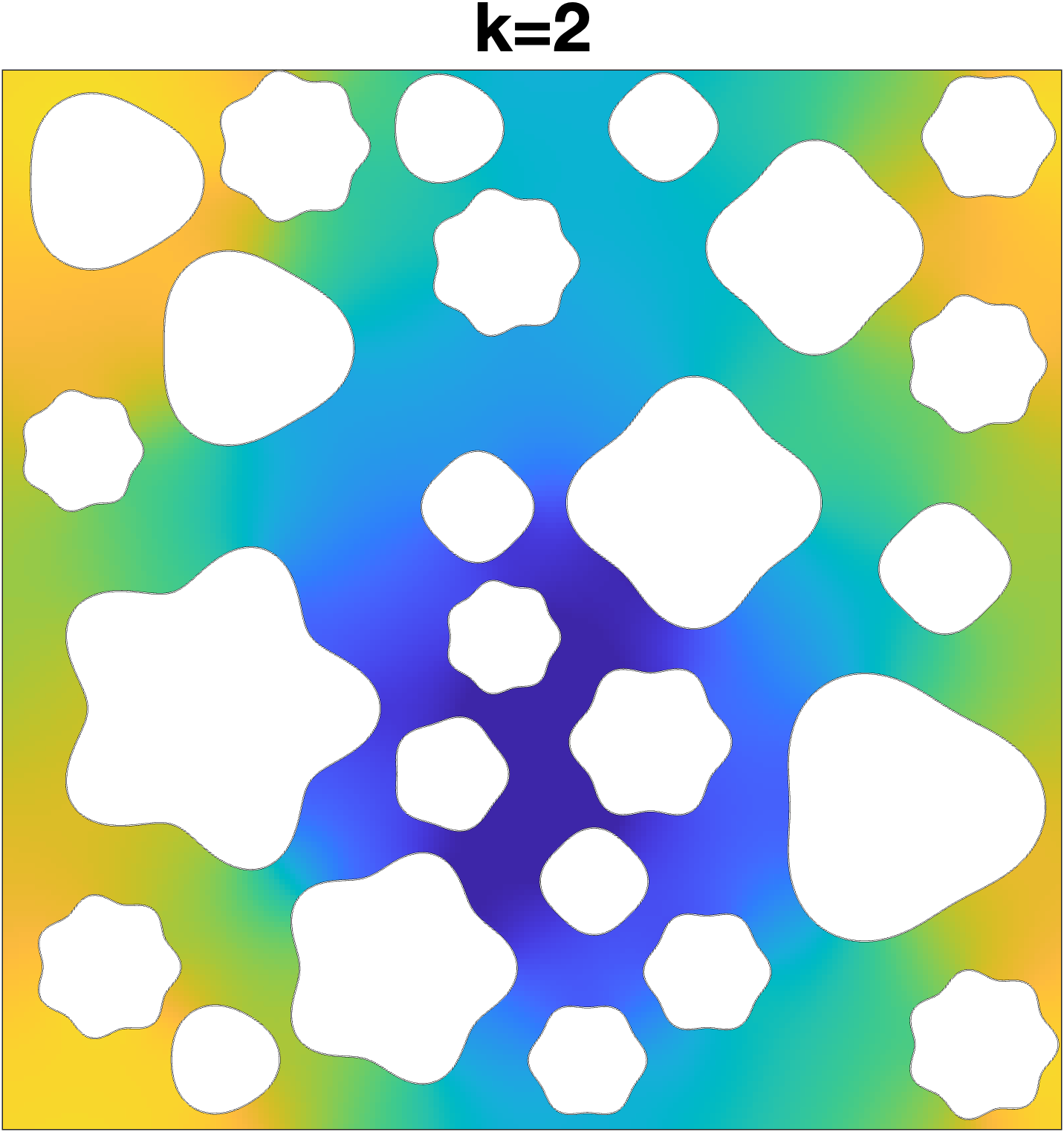}
\includegraphics[width=.38\textwidth]{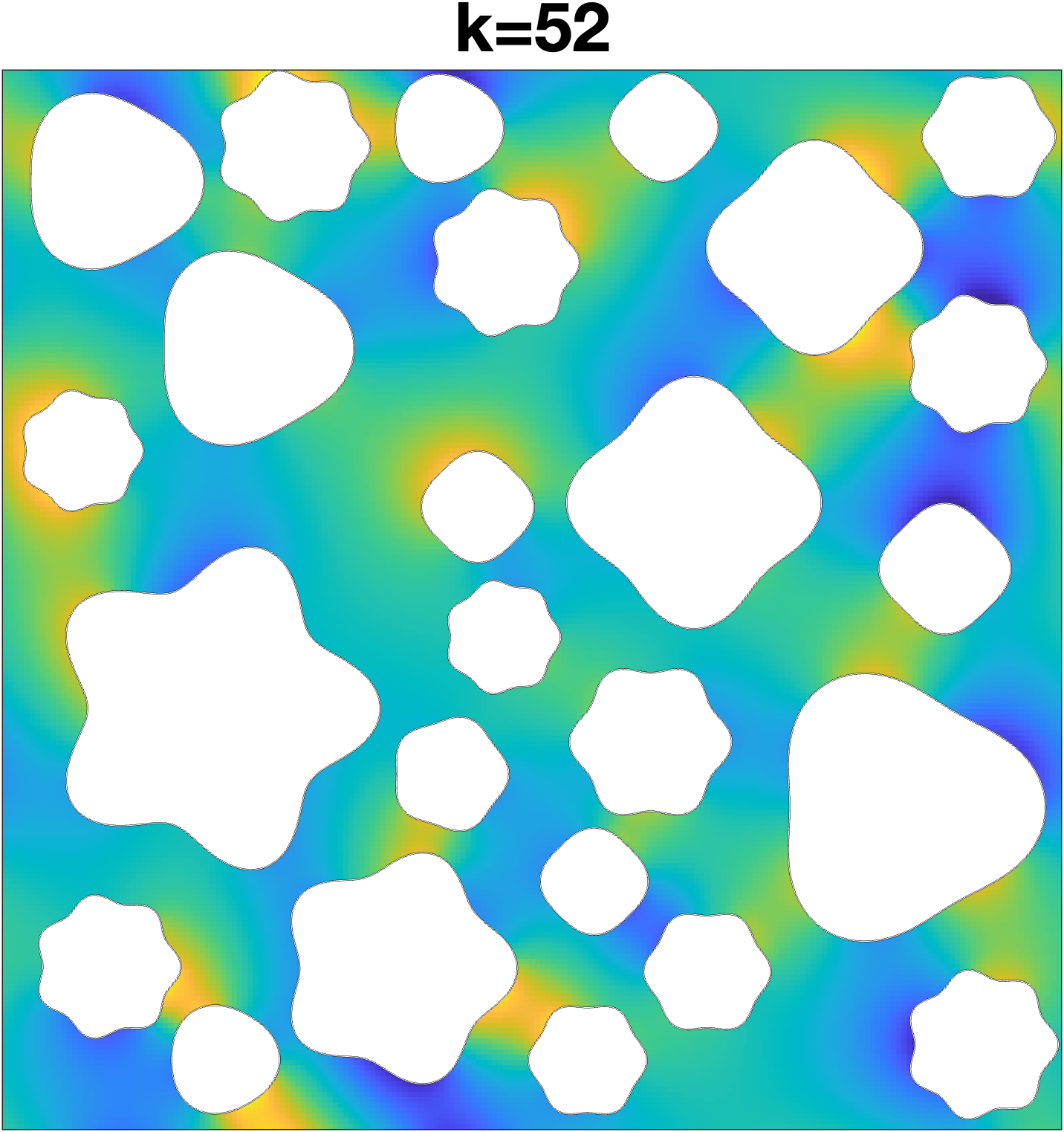}\\
\vspace{3mm}
\includegraphics[width=.38\textwidth]{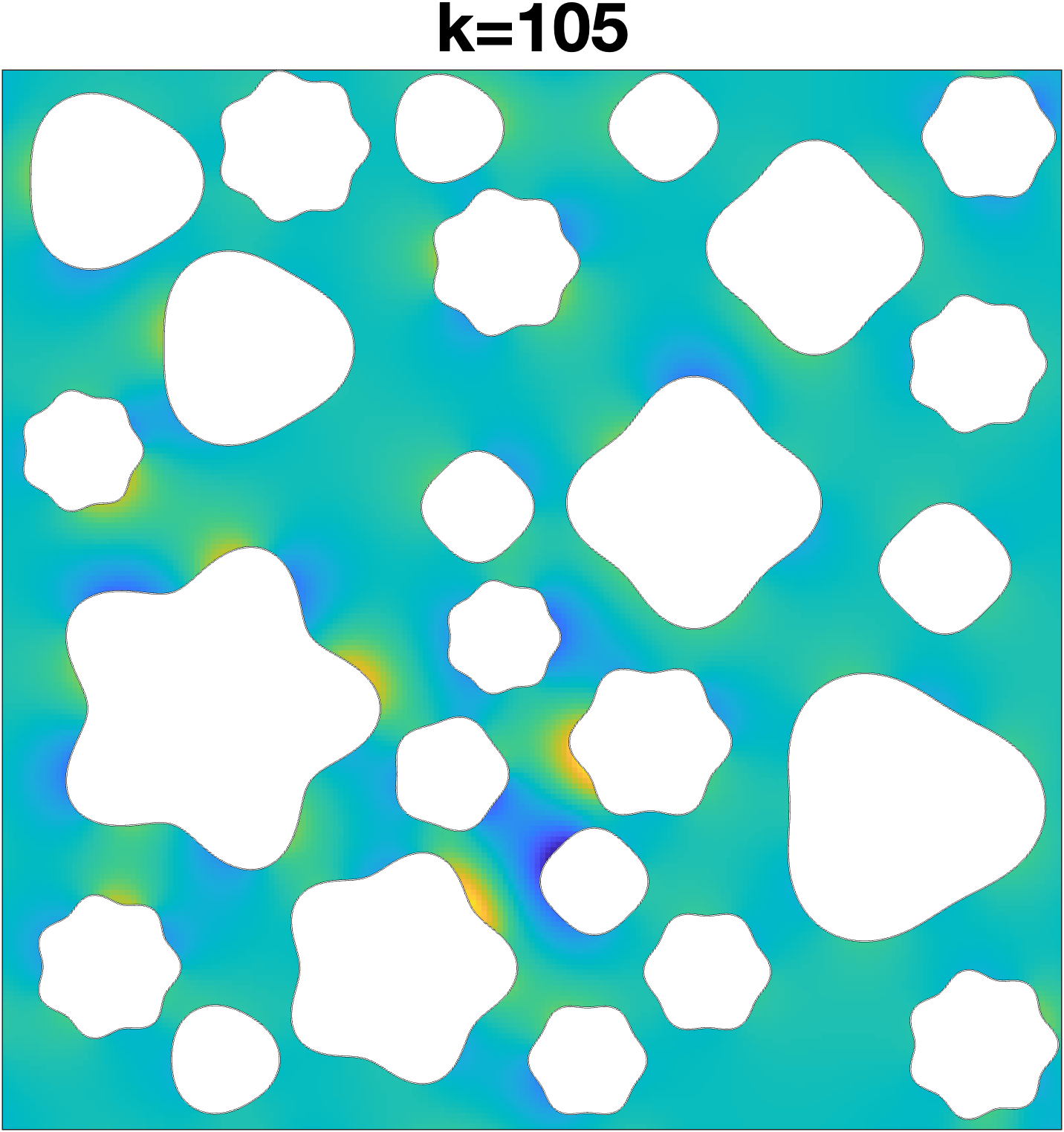}
\includegraphics[width=.38\textwidth]{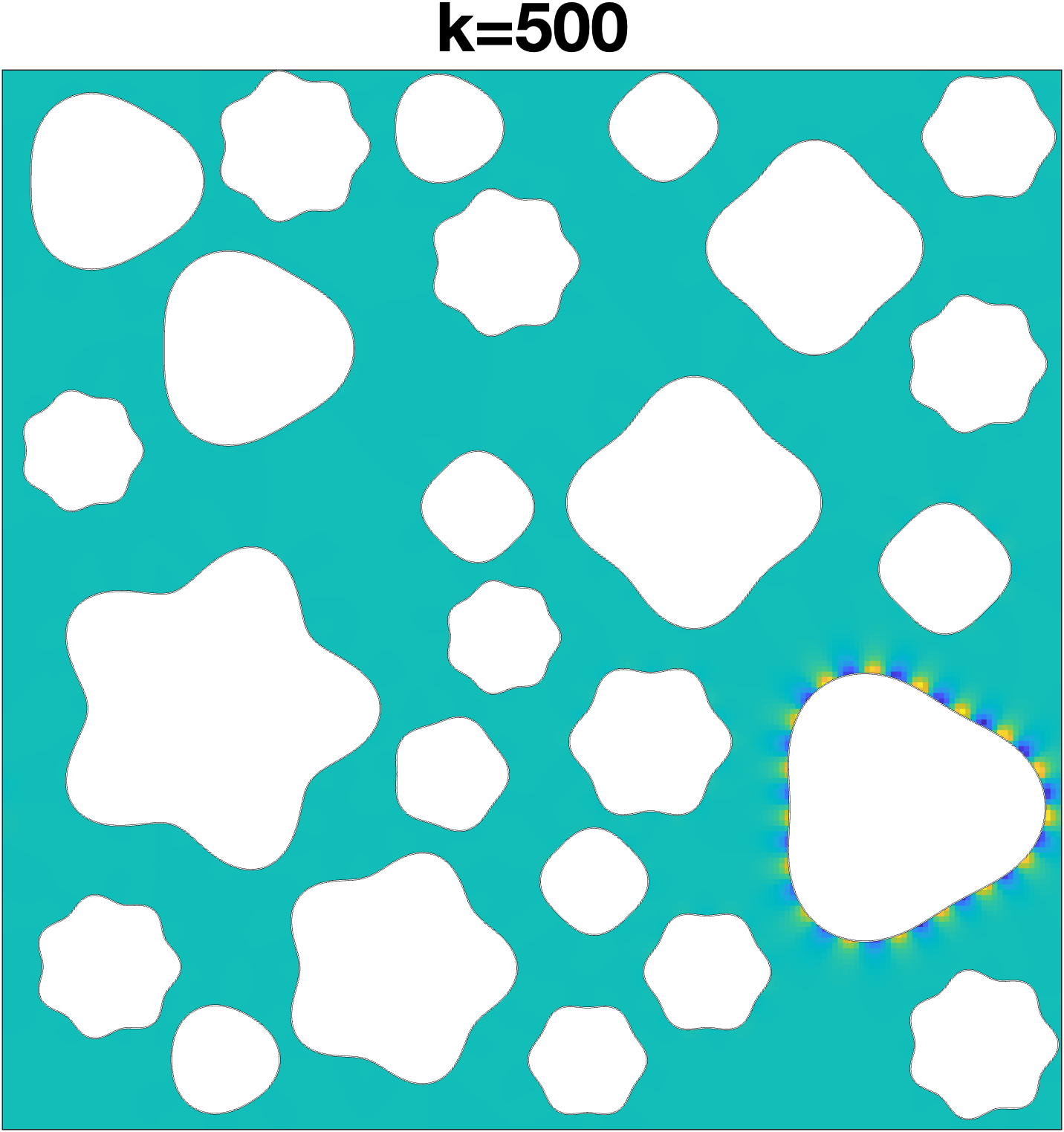}
\end{center}
\caption{Eigenfunctions $u_k$ of the Steklov EVP~\eqref{eqn: steklov} corresponding to $\sigma_k$, for $k=2$, $52$, $105$, and $500$, are computed on a square torus with $M=25$ holes. 
The corresponding eigenvalues are tabulated in Table~\ref{tabl: steklov many err}. 
Each hole is discretized with $N=200$ boundary points, resulting in a total of $5000$ degrees of freedom.
Using an a posteriori estimate, the eigenvalue approximations have a relative error less than $10^{-14}$. 
As anticipated, the eigenfunctions $u_k$ corresponding to large $\sigma_k$ concentrate near $\partial \Omega$ and oscillate with a wavelength $\approx \frac{2\pi}{\sigma_k}$ \cite{girouard2017spectral,polterovich2019nodal}. See Sec.~\ref{sec: steklovbvp numerical}, Example 6.}
\label{fig:many holes}
\end{figure}

\vspace{3mm} 
\noindent {\bf Example 6. Steklov EVP with $M=25$ holes.}
We solve the Steklov EVP on a square torus with $M=25$ holes; the results are shown in Fig.~\ref{fig:many holes}. Following the Neumann BVP setup, for each hole, we randomly choose each center $a_i$ for $i=1,2,\ldots, M$, the oscillation factor $\omega_i$ with $3\leq \omega_i\leq7$, and the maximum radius of the hole \( r_i \). The boundary of each hole \(\partial D_i \) is parametrized by~\eqref{eqn: neumann bdry shape}.

We use $N_i=200$ ($N=50000$ total) grid points per boundary component to discretize the boundary integral equation~\eqref{e: steklov bie}.
The eigenfunctions $u_k$ in~\eqref{e:uRep, stek}, corresponding to $\sigma_k$, for $k=2$, $52$, $105$, and $500$, are plotted in Fig.~\ref{fig:many holes}. 
As anticipated, the eigenfunctions \( u_k \) associated with large values of \( \sigma_k \) become concentrated near \( \partial \Omega \) and exhibit oscillations with a wavelength approximately equal to \( \frac{2\pi}{\sigma_k} \) \cite{girouard2017spectral,polterovich2019nodal}. 
The corresponding eigenvalues are reported in Table~\ref{tabl: steklov many err}. We note that an a posteriori error estimate can be used to bound the relative errors in $\sigma_k$ by $\Vert\partial_\nu u_k - \sigma_k u_k \Vert_{L^2(\partial \Omega)}$; see \cite{BOGOSEL2016265} and \cite[Prop.~4.1]{kaoharmonic2023}. In Table~\ref{tabl: steklov many err}, we report $\Vert \partial_\nu u_k - \sigma_k u_k \Vert_{L^\infty(\partial \Omega)}$, approximated using a finer discretization on the boundary. We estimate that the errors in $\sigma_k$ are on the order of $10^{-14}$. 

\begin{table}[t!]
\centering
\begin{tabular}{|c|c|c|}
\hline
$k$ & $\sigma_k$ & $\Vert \partial_\nu u_k - \sigma_k u_k \Vert_{L^\infty(\partial \Omega)}$ \\
\hline
$2$ & $0.9989357$ & $9.760 \times 10^{-15}$ \\ \hline
$52$ & $19.57834$ & $1.960 \times 10^{-14}$ \\ \hline
$105$ & $35.13811$ & $1.664 \times 10^{-14}$ \\ \hline
$500$ & $138.0663$ & $1.386 \times 10^{-14}$ \\ \hline
\end{tabular}
\vspace{0.5em}
\caption{For the Steklov EVP on a square torus with 25 holes, we tabulate values of  $\sigma_k$ and $\Vert \partial_\nu u_k - \sigma_k u_k \Vert_{L^\infty(\partial \Omega)}$ for different values of $k$. See Sec.~\ref{sec: steklovbvp numerical} Example 6 and Fig.~\ref{fig:many holes}.}
\label{tabl: steklov many err}
\end{table}

\section{Discussion} \label{s:disc}
In this paper, we have developed and analyzed layer potential methods to represent harmonic functions on finitely-connected tori. 
The layer potentials are expressed in terms of a doubly-periodic and non-harmonic  Green's function $G$ in~\eqref{eqn: green}. 
Extending results for Euclidean domains, we establish in Lemmas~\ref{lem: properties single layer} and~\ref{lem: properties double layer} that the single- and double-layer potential operators are compact linear operators and we derive relevant limiting properties at the boundary. 
Here, we also show that when the boundary has multiple connected component, the Fredholm operator of the second kind, $\K^*-\onehalf I $, possesses a non-trivial null space, for which we construct a basis.  In 
Theorems~\ref{thm: dirichlet},~\ref{thm: neumann}, and~\ref{thm: steklov}, we use the layer potentials to represent solutions to the Dirichlet BVP~\eqref{eqn: dirichlet}, Neumann BVP~\eqref{eqn: neumann}, and Steklov EVP~\eqref{eqn: steklov}, respectively. 
Finally, we implement the developed methods and demonstrate their accuracy across several numerical examples; see Sec.~\ref{sec: Numerical Experiments}.

There are several interesting extensions of this work. 
First, our method can be improved by incorporating techniques for close evaluation \cite{barnett2014boundary,Gedney2003}, 
using fast multipole acceleration to solve linear integral equations \cite{GREENGARD1987325},  
and handling less regular boundaries (e.g., boundary with corners) \cite{Goodwill_2023,grisvard2011elliptic,mclean2000strongly}. 
While we solve the Laplace equation on finite-connected tori here, extending these ideas to the Helmholtz and Stokes equations on finitely-connected tori is a natural next step. 
In this paper we focus on solving BVPs and EVPs on genus one surfaces; from this perspective, the work extends research on solving the Laplace-Beltrami equation on a multiply-connected sphere \cite{Kropinski_2013}. It is interesting to extend layer potential methods to higher genus surfaces; recently, progress has been made in this direction for the method of particular solutions  \cite{nahon2024}. Finally, as we focus on two-dimensional tori, it is interesting to consider developing layer potential methods to represent harmonic functions on higher-dimensional tori.

\section*{Acknowledgments}
We thank the anonymous referees for their careful reading of the manuscript and their helpful suggestions.

\bibliographystyle{siamplain}
\bibliography{references}

\clearpage
\appendix 
\section{Proofs of lemmas in Sections~\ref{sec:LayerPotentialProperties}--\ref{sec: BIE}}
\label{sec: proofs of lem thms-properties}

\subsection{Preliminary Estimates} This subsection provides the propositions and lemmas required for the proofs in Secs.~\ref{sec:LayerPotentialProperties}--\ref{sec: BIE}. Specifically, the following proposition summarizes the asymptotic properties of $G(z-\xi)$ and its normal derivatives (e.g., \eqref{eqn: asymptotic G} and~\eqref{eqn: asymptotic Gnu}), which are utilized throughout the proofs of Lemmas~\ref{lem: exact double layer}--\ref{lem: properties single layer}.
\begin{prop}[Asymptotic properties]
\label{prop: asymptotic}
   Let \(\Omega\) satisfy Assumption~\ref{assumption1} and \(z_0, \xi \in \partial D_j\) be sufficiently close. Define $z = z_0 \mp h\nu(z_0)$ with $h = |z-z_0|$ ($-$ for $z\in \Omega$ and $+$ for $z\in D_j$). Then, the following properties hold:  
    \begin{align}
        G(z - \xi) & = -\frac{1}{2\pi} \log|z - \xi| + c + O(|z - \xi|^2), \label{e: green approx}\\
        \partial_{\nu_\xi} G(z - \xi) &= \frac{(z-\xi) 
 \circ \nu(\xi)}{2\pi |z - \xi|^2} +   O(|z-\xi|), \label{eqn: dgdxi close}\\
  \partial_{\nu_z} G(z - \xi) &=-\frac{(z-\xi) 
 \circ \nu(z_0)}{2\pi |z - \xi|^2} 
     +  O(|z-\xi|), 
\label{eqn: dgdz close}\\
\left \vert\partial_{\nu_\xi} G(z - \xi) -\partial_{\nu_\xi} G(z_0 - \xi)\right\vert &\leq \frac{1}{2\pi}\left[ \frac{2h}{h^2 + |z_0-\xi|^2}   + 3L_{1,j}\right]+  O(h), \label{eqn: dlayer inequality}
\end{align}
where $O(|z-\xi|)$ terms are absolutely convergent,  $c = \frac{\log|\vartheta_1'(0)|}{2\pi}$, and $L_{1,j}$ are constants for each $j=1,\ldots,M$.
\end{prop} 

\begin{proof}
To obtain~\eqref{e: green approx}, we use the Weierstrass sigma expansion:
\begin{equation*}
G(z - \xi) = -\frac{1}{2\pi}\log|\sigma(z - \xi)| + c + O(|z - \xi|^2),
\end{equation*}
where $c = \frac{\log|\vartheta_1'(0)|}{2\pi}$.
 This follows from a well-known identity~\cite{chandrasekharan2012elliptic,lin2010elliptic},
\begin{align*}
\sigma(z) = e^{\eta_1 z^2/2} \frac{\vartheta_1(z)}{\vartheta_1'(0)}. 
\end{align*}
Furthermore, using the product definition of $\sigma$, we may approximate 
\begin{align*}
    \log|\sigma(z - \xi)| &= \log|z - \xi| + \log\left|\prod_{\omega \neq 0} \left(1 - \frac{z - \xi}{\omega}\right) e^{(z - \xi)/\omega + (z - \xi)^2/2\omega^2} \right| \\
    &= \log|z - \xi|+ O(|z - \xi|^3), 
\end{align*}
where the infinite product is absolutely and uniformly convergent for \(|z| \leq R\) and has order \(O(|z - \xi|^3)\) (see \cite{daniels1883note}), from which~\eqref{e: green approx} follows.

To obtain~\eqref{eqn: dgdxi close}-\eqref{eqn: dgdz close}, we use the Laurent series of the Weierstrass zeta~\cite{chandrasekharan2012elliptic}: \begin{equation}
\zeta(z) = \frac{1}{z} - \sum_{n=2}^{\infty} c_n z^{2n-1}.
\label{eqn: laurent weierstrass zeta}
\end{equation} 
The series term $O(|z|)$ is absolutely convergent, from which \eqref{eqn: dgdxi close}-\eqref{eqn: dgdz close} follows.

Lastly, to obtain \eqref{eqn: dlayer inequality}, we utilize the $C^2$ regularity of $\partial D_j$ (\cite[Lemma 6.16]{kress1989linear}): there exist positive constants \(L_{1,j}, L_{2,j}\) such that for all \(z_0, \xi \in \partial D_j\),
\begin{align}
    |(z_0 - \xi) \circ \nu(\xi)| &\leq L_{1,j} |z_0 - \xi|^2, \quad
    |\nu(z_0) - \nu(\xi)| \leq L_{2,j} |z_0 - \xi|. \label{lemma: C2 prop}
\end{align}
For \(z_0, \xi \in \partial D_j\) sufficiently close, we choose $h= |z-z_0| <1/(2 L_{1,j})$. Using \eqref{lemma: C2 prop} and~\cite[Lemma 3.20]{folland2020introduction}, we get
\[|z-\xi|^2 \geq \frac{1}{2}\left(|z-z_0|^2+|z_0-\xi|^2\right)\quad \implies \quad\frac{1}{|z-\xi|^2} \leq \frac{2}{|z-z_0|^2 +|z_0-\xi|^2}.
\]
We then compute
{\small
\begin{align*}
 \left \vert\partial_{\nu_\xi} G(z - \xi) -\partial_{\nu_\xi} G(z_0 - \xi)\right\vert  & \leq \frac{1}{2\pi} \left[\left\vert \frac{(z-\xi) 
 \circ \nu(\xi)}{|z - \xi|^2} \right\vert  + \left\vert \frac{(z_0-\xi) 
 \circ \nu(\xi)}{|z_0 - \xi|^2} \right\vert \right]+  O(|z-z_0|)\\
&\leq \frac{1}{2\pi} \left[ \frac{2|z-z_0|+2L_{1,j}|z_0-\xi|^2 }{|z-z_0| + |z_0-\xi|^2}   + L_{1,j}\right]+  O(|z-z_0|)
\\
&\leq \frac{1}{2\pi}\left[ \frac{2|z-z_0|}{|z-z_0|^2 + |z_0-\xi|^2}   + 3L_{1,j}\right]+  O(|z-z_0|).
\end{align*}}
\end{proof}

Next, we introduce a complex double-layer potential and its properties, which are used to prove Lemmas~\ref{lem: properties double layer}(4)--(6) and Lemma~\ref{lem: characteristic}. One can understand this potential as the analogue to the Cauchy integral~\cite{kress1989linear,MIKHLINdirichlet,Nasser_2011}.

\begin{definition}[Complex double-layer potential] 
Given $\phi \in C(\partial \Omega)$, we define the complex double-layer potential for $z\in \mathcal{P}_\tau \setminus \partial\Omega$ by
\begin{equation}
\mathcal{D}_f [\phi] (z) \coloneq \int_{\partial\Omega} \frac{\phi(\xi)}{2\pi i} \left[\zeta(z - \xi) - \eta_1(z - \xi) \right] |d\xi| 
- \frac{i}{b}\int_{\partial\Omega}(z- \xi)\phi(\xi)  \Re[d\xi], 
\label{eqn: complex double layer}
\end{equation}
where $\mathcal{P}_\tau\subset \mathbb{C}$ is the fundamental parallelogram,
\begin{equation}
\label{e:parallelo}
    \mathcal{P}_\tau \coloneq \{x + \tau y \in \mathbb{C}: (x,y)\in [0,1]^2\}. 
\end{equation} 
Here, $\zeta(z)$ is the Weierstrass zeta function with quasi-period $\eta_1 \coloneq 2\zeta(1/2)$. The associated boundary operator $\K_f \colon C^2(\partial\Omega) \to C^2(\partial\Omega)$ is defined for $z_0 \in \partial\Omega$ by 
\begin{equation}
    \K_f[\phi](z_0) \coloneq \mathrm{p.v.}\int_{\partial\Omega} \frac{\phi(\xi)}{2\pi i} \left[\zeta(z_0 - \xi) - \eta_1(z_0 - \xi) \right]|d\xi| 
- \frac{i}{b}\int_{\partial\Omega}(z- \xi)\phi(\xi) \Re[d\xi], 
\end{equation}
where $\mathrm{p.v.}$ denotes the Cauchy principal value.
\end{definition}
We note that the last term of $\mathcal{D}_f[\phi](z)$ in~\eqref{eqn: complex double layer} can be expanded as:$$-\frac{i}{b}\int_{\partial\Omega}(z- \xi)\phi(\xi) \, \Re[d\xi] =  -\frac{iz}{b}\left[\int_{\partial\Omega}\phi(\xi) \, \Re[d\xi]\right] + \frac{i}{b}\left[\int_{\partial\Omega}\xi\phi(\xi) \, \Re[d\xi]\right].$$Since $\phi \in C(\partial\Omega)$, the integrals in the brackets are well-defined constants in $\mathbb{C}$.
\begin{lemma}[Properties of the Complex Double-Layer Potential]
\label{lem: Complex layer} Let \(\Omega\) satisfy Assumption~\ref{assumption1}. 
The complex double-layer potential \eqref{eqn: complex double layer} satisfies the following: 
\begin{enumerate}
\item 
For $\phi \in C(\partial \Omega)$, \( \mathcal{D}_f[\phi] \) is analytic on $\mathcal P_\tau \setminus \partial\Omega$; therefore, both $\Re  \mathcal{D}_f[\phi]$ and $\Im \mathcal{D}_f[\phi]$ are harmonic on $P_\tau \setminus \partial\Omega$. Notably, $\Re \mathcal{D}_f[\phi] = \mathcal{D}[\phi]$ is a function on $\mathbb T_\tau\setminus\partial\Omega$, whereas $\Im \mathcal{D}_f[\phi]$ is not necessarily doubly-periodic. 

\item  For $\phi \in C(\partial \Omega)$, the complex double layer potential satisfies
\[
\mathcal{D}_f[\phi](z_0^{\pm}) =\K_f[\phi](z_0) \mp \onehalf\phi(z_0), \qquad z_0 \in \partial \Omega.
\]
\end{enumerate}
\end{lemma}
Note that Lemma~\ref{lem: Complex layer}(2) follows from the Sokhotski-Plemelj formula in~\cite[Thm. 7.8]{kress1989linear}, since $\zeta(z-\xi)$ admits a Laurent series with leading order $(z-\xi)^{-1}$~\cite{lin2010elliptic}.
\begin{proof}  
\noindent {\bf (1)} 
Following \cite{Nasser_2011}, we express $\mathcal{D}[\phi]$ as the real part of a complex function. Let \(\xi = \xi(t)\) parametrizes \(\partial\Omega\) with the normal vector $\nu(\xi) = -i\frac{\xi'(t)}{|\xi'(t)|}$. Then, we have
\begin{align*}
    \partial_{\nu_\xi} G(z - \xi) &= -(G_{x}(z - \xi) + iG_{y}(z - \xi)) \circ \nu(\xi)\\
    & = \Re\left[\frac{\xi'(t)}{i|\xi'(t)|} \left( \frac{1}{2\pi } (\zeta(z - \xi) - \eta_1(z - \xi)) + \frac{i\Im(z-\xi)}{b} \right)\right]. 
\end{align*}
Integrating against $\phi(\xi)$ along $\partial\Omega$ yields
\begin{align*}
    \mathcal{D}[\phi](z) = \Re\left[\int_{\partial\Omega} \left(\frac{1}{2\pi i}(\zeta(z - \xi) - \eta_1(z - \xi)) +\frac{\Im(z-\xi)}{b}\right)\phi(\xi)\,d\xi\right].
\end{align*}
Note that the first two terms on the right-hand side are analytic, whereas the term $\frac{\Im(z-\xi)}{b}$ is non-analytic. Since the latter is continuous, we apply the Cauchy-Riemann equations to identify its harmonic conjugate:
\begin{equation*}
-\int_{\partial\Omega} \frac{\Re(z-\xi)}{b}\phi(\xi)\,\Re\left[d\xi\right] = -\Im \left[ \frac{i}{b}\int_{\partial\Omega}(z- \xi)\phi(\xi) \, \Re[d\xi] \right].
\end{equation*}
Considering all the terms, we get $\Re \mathcal{D}_f[\phi] (z) = \mathcal{D}[\phi] (z)$ for $\mathcal{D}_f[\phi] (z)$ in \eqref{eqn: complex double layer}.

\noindent {\bf (2)} 
The result follows from the Sokhotski-Plemelj formula \cite[Thm. 7.8]{kress1989linear}, as $\zeta(z-\xi)$ has a Laurent series~\eqref{eqn: laurent weierstrass zeta} with leading order $(z-\xi)^{-1}$. Thus, for $ \mathbb{T}_\tau\ni z\to z_0\in \partial\Omega$,
\begin{align*}
    \mathcal{D}_f[\phi](z)- \K_f[\phi](z_0) = \frac{1}{2\pi i} \int_{\partial\Omega} \phi(\xi)\left[\frac{1}{z-\xi}-\frac{1}{z_0-\xi} \right] \, |d\xi| + O(z-z_0) \to \mp \frac{\phi(z_0)}{2}.  
\end{align*}

\end{proof}

\subsection{Proof of Lemma~\ref{lem: exact double layer}}
\begin{figure}[ht!]
\small{
\begin{center}
\includegraphics[width=0.2\textwidth]{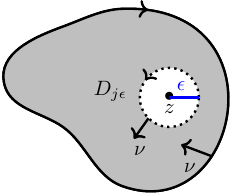}
\hspace{5em}
\includegraphics[width=0.40\textwidth]{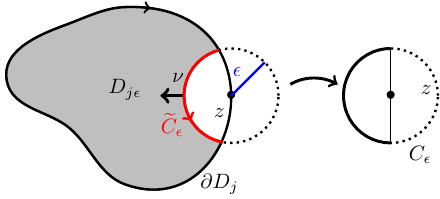}
\end{center}}
\caption{Diagrams for the proof of Lemma~\ref{lem: exact double layer}: {\bf~(left)} \(z \in D_j\); {\bf (right)}  \(z \in \partial D_j\).}
\label{fig: lemma 4 diagrams}
\end{figure}

\begin{proof}  We first prove~\eqref{eqn: gen gauss sublemma}. Without loss of generality, fix $j $. We use Green's formula (cf. \cite[Ch. 15.6]{henrici1986applied} ) and \eqref{eqn: not harmonic G}. \\

\noindent {\bf Case $z \in \mathbb{T}_{\tau} \setminus \overline{D_j}$:} 
\begin{equation*}
    \int_{\partial D_j} \partial_{\nu_\xi}G(z-\xi) \, |d\xi| = -\int_{D_j} \Delta_\xi G(z-\xi) \, d\xi  = -\frac{|D_j|}{b}. 
\end{equation*}\\
\noindent {\bf Case $z \in D_j$:} Let $\varepsilon > 0$ and  $D_{j\varepsilon} = D_j \setminus B_{\varepsilon}(z)$ for an open ball $B_{\varepsilon}(z) \subset D_j$ (see Fig.~\ref{fig: lemma 4 diagrams} (left)). The Green's formula yields: 
\begin{align*}
    -\frac{\vert D_{j\varepsilon} \vert}{b} &= -\int_{D_{j\varepsilon}} \Delta_\xi G(z-\xi) \, d\xi   = \int_{\partial D_{j}} \partial_{\nu_\xi}G(z-\xi) \, |d\xi| + \int_{\partial B_{\varepsilon}(z)} \partial_{\nu_\xi}G(z-\xi) \, |d\xi|,
\end{align*}
and~\eqref{eqn: asymptotic Gnu} implies
\begin{align*}
        \int_{\partial B_{\varepsilon}(z)} \partial_{\nu_\xi}G(z-\xi) \, |d\xi| &= \int_{\partial B_{\varepsilon}(z)} \frac{(z-\xi)\circ (\xi-z)}{2\pi |z-\xi|^3} 
        +O(|z-\xi|) \, |d\xi|= -1 + O(\varepsilon^2).
\end{align*} Letting \(\varepsilon \to 0\) yields the desired result.

\noindent {\bf Case $z \in \partial D_j$:} 
Let $\varepsilon>0$ and define (see Fig.~\ref{fig: lemma 4 diagrams} (right)) 
\begin{equation*}
\begin{split}
    & \widetilde{C}_\varepsilon \coloneq \partial D_{j\varepsilon} \cap \partial B_{\varepsilon}(z) 
    \,\, \textrm{and} \,\,    
    C_\varepsilon \coloneq \{y \in \partial B_{\varepsilon}(z) \colon  \nu(z)\circ ({y-z}) \geq 0\};
\end{split}
\end{equation*}
Then, Green's formula and~\eqref{eqn: asymptotic Gnu}
yields
\begin{align*}
{-}\frac{\vert D_{j\varepsilon}\vert}{b} &={-} \int_{D_{j\varepsilon}}\Delta_\xi G(z-\xi)\,d\xi 
 = \int_{\partial D_{j\varepsilon} \setminus \widetilde{C}_\varepsilon} \partial_{\nu_\xi}G(z-\xi)\,|d\xi| + \int_{\widetilde{C}_\varepsilon} \partial_{\nu_\xi}G(z-\xi)\,|d\xi| \\
&=\int_{\partial D_{j\varepsilon} \setminus \widetilde{C}_\varepsilon} \partial_{\nu_\xi}G(z-\xi)\,|d\xi| -\frac{1}{2} + O(\varepsilon).
\label{eqn: exact double at boundary D_j}
\end{align*} For the last equality, we use  \cite[Prop. 3.19]{folland2020introduction},
\begin{equation*}
    \int_{\widetilde{C}_\varepsilon} |d\xi| = \int_{C_\varepsilon} |d\xi| + O(\varepsilon^2) = \pi \varepsilon + O(\varepsilon^2).
\end{equation*}
Again, letting \(\varepsilon \to 0\) yields the desired result.

Finally, we use~\eqref{eqn: gen gauss sublemma} to get~\eqref{eqn: gen gauss}. We sum over $j$ to obtain \eqref{eqn: gen gauss}. For $z \in \Omega$, the sum is $-\sum |D_j|/b = -|D|/b$. If $z \in D_k$, the $k$-th term contributes an additional $1$, yielding $1 - |D|/b$. If $z \in \partial D_k$, the $k$-th term contributes $1/2$, yielding $1/2 - |D|/b$. 
\end{proof}

\subsection{Proof of Lemma~\ref{lem: properties double layer}}
\begin{proof} We follow the proofs for the analogous statements for double-layer potentials on Euclidean domains, as found in \cite[Ch. 3]{folland2020introduction}, \cite[Ch. 6-7]{kress1989linear}, \cite[Ch. 4]{MIKHLINdirichlet}.\\ 
\noindent{\bf (1)} The operator \(\mathcal{D}\) is linear, and its periodicity follows from the doubly-periodic kernel \(\partial_{\nu_\xi} G(z - \xi)\), since $G$ is itself doubly-periodic~\cite{lin2010elliptic}. By~\eqref{eqn: not harmonic G}, we get
\begin{equation*}
        \Delta \mathcal{D}[\phi](z)  = \int_{\partial \Omega} \phi(\xi) \Delta_z \partial_{\nu_\xi}G(z-\xi)\,|d\xi| = 0, \qquad   z \in \mathbb{T}_\tau\setminus\partial\Omega.
\end{equation*}

\noindent {\bf (2)} We prove the result for \( z_0^{+} \); the case for \( z_0^{-} \) follows similarly. Let $z_0 \in \partial D_j$ and $z\in \Omega$ for some $j= 1,2,\ldots, M$. Applying~\eqref{eqn: gen gauss}, we can write
\begin{equation*}
\mathcal{D}[\phi](z) - \K[\phi](z_0) = \int_{\partial \Omega} [\phi(\xi) - \phi(z_0)] \left[ \partial_{\nu_\xi}G(z-\xi) - \partial_{\nu_\xi}G(z_0-\xi) \right] \, |d\xi| \ - \  \onehalf\phi(z_0).
\end{equation*}
The proof is completed by showing that:
\begin{equation}
    \lim_{z \to z_0} \int_{\partial \Omega} [\phi(\xi) - \phi(z_0)] \left[ \partial_{\nu_\xi}G(z-\xi) - \partial_{\nu_\xi}G(z_0-\xi) \right] \, |d\xi| = 0.
    \label{eqn: limit zero}
\end{equation}
\noindent\textit{Proof of \eqref{eqn: limit zero}:} For $\varepsilon > 0$, decompose the integral in~\eqref{eqn: limit zero} into $(I)+(II)$:
\begin{align*}
    (I)   &=  \int_{\partial\Omega \setminus B_\gamma(z_0)} [\phi(\xi) - \phi(z_0)] \left[ \partial_{\nu_\xi}G(z-\xi) - \partial_{\nu_\xi}G(z_0-\xi) \right] \, |d\xi|,\\
    (II) &= \int_{\partial D_j \cap B_\gamma(z_0)} [\phi(\xi) - \phi(z_0)] \left[ \partial_{\nu_\xi}G(z-\xi) - \partial_{\nu_\xi}G(z_0-\xi) \right] \, |d\xi|, 
\end{align*}
where $B_\gamma(z_0)$ is a ball of \textit{fixed} radius $\gamma>0$.  We show both $|\text{(I)}|$ and $|\text{(II)}|$ can be made less than $ \varepsilon/2$ for $z$ sufficiently close to $z_0.$ 

To bound $|\text{(I)}|$, let  $S = \partial\Omega \setminus B_\gamma(z_0)$. Since $S$ is compact, \(\partial_{\nu_\xi} G(z - \xi)\) is uniformly continuous at $z_0$ for \(\xi \in S\). 
Therefore, choosing $z$ close to $z_0$ will allow 
\begin{align*}
\vert(I)\vert
&\leq \left\Vert \partial_{\nu_\xi}G(z-\xi) - \partial_{\nu_\xi}G(z_0-\xi) \right\Vert_{L^{\infty}(S)} \left\Vert \phi(\xi) - \phi(z_0) \right\Vert_{L^{\infty}(S)} \int_{S} |d\xi| \\
&\leq \frac{\varepsilon}{4\Vert\phi\Vert_{L^{\infty}(\partial \Omega)}|\partial \Omega|}  2\Vert\phi\Vert_{L^{\infty}(\partial \Omega)} |\partial\Omega| \leq \frac{\varepsilon}{2}.
\end{align*} 
To bound $\vert(II)\vert $, let $z = z_0-h\nu(z_0)$, with $h = \vert z- z_0 \vert$, and define $\tilde{S} =\partial D_j \cap B_\gamma(z_0) $. Since $\mathcal{D}[\phi]$ is continuous in $\Omega$, it suffices to consider $z$ along the normal. From~\eqref{eqn: dlayer inequality}, 
\begin{align*}
|\text{(II)}| 
&\leq\Vert\phi(z_0)-\phi(\xi)\Vert_{L^{\infty}(\tilde{S})} \left[ \frac{1}{2\pi}\int_{\tilde{S}} \frac{2h}{h^2 + |z_0-\xi|^2} + 3L_{1,j}\, |d\xi| +  O(h)\right].
\end{align*}
We consider terms inside the bracket. For the term containing the integral, we use polar coordinates $r = |z_0 - \xi|$~\cite[Ch. 0]{folland2020introduction}\ to obtain:
\begin{align*}\frac{1}{2\pi} \int_{\tilde{S}} \frac{2h}{h^2 + r^2} + 3L_{1,j} \, |d\xi| \leq 2 \int_0^\gamma \frac{h}{h^2 + r^2} \, dr + \frac{3L_{1,j}|\partial D_j|}{2\pi} \to \pi + \frac{3L_{1,j}|\partial D_j|}{2\pi}\end{align*}as $h \to 0$. Consequently, as $z \to z_0$, the above expression is bounded by a constant $C$, while the remaining $O(h)$ term can also be made smaller than $C$. Finally, allowing  $\Vert \phi(z_0) - \phi(\xi) \Vert_{L^\infty(\tilde{S})} < \varepsilon/4C$ ensures $|II| \leq \varepsilon/2$.

\noindent {\bf (3)} We show that \(\partial_{\nu_\xi} G(z-\xi)\) can be continuously extended as $z\to\xi \in \partial\Omega$, implying that \(\K\) is compact~\cite[Thm. 2.27]{kress1989linear}. Suppose $z,\xi \in \partial D_j$ are parameterized by \(z(s) = x(s) + iy(s)\), \(\xi = z(t)\) for $s,t \in [0,2\pi)$. Using~\eqref{eqn: asymptotic Gnu} and~\cite[Sec. 2]{kress1991boundary}, we obtain
\begin{align*}
    \lim_{z \to \xi} \partial_{\nu_\xi} G(z - \xi)
    = \lim_{s \to t} -\frac{(z(s)-z(t))\circ iz'(t)}{2\pi |z'(t)| |z(s) - z(t)|^2} 
    = \frac{y'(t) x''(t) -x'(t) y''(t)}{4\pi |z'(t)|^3} = -\frac{\kappa(\xi)}{4\pi}.
\end{align*}
Here, the outer normal vector \(\nu(\xi)\) and the signed curvature \(\kappa(\xi)\)  are given by
\[
    \nu(\xi) = \frac{y'(t) - ix'(t)}{|z'(t)|}  
    \quad \textrm{and} \quad 
    \kappa(\xi) = \frac{x'(t)y''(t) - y'(t)x''(t)}{|z'(t)|^3}.
\]
The compactness of the adjoint $\K^*$ follows as well~\cite[Thm. (0.37)]{folland2020introduction}. We also obtain
\[
    \lim_{z \to \xi} \partial_{\nu_z} G(z - \xi) = \lim_{z \to \xi} \overline{\partial_{\nu_z} G(z - \xi)} = -\frac{\kappa(\xi)}{4\pi}.
\]

\noindent {\bf (4)} We first show that $\dim N(\K-\onehalf I ) \leq M-1$. Suppose $\phi \in N(\K-\onehalf I )$. By Lemma~\ref{lem: properties double layer}(1)--(2), $\mathcal{D}[\phi]$ solves the homogeneous Dirichlet BVP on $\Omega$, which, by the maximum principle, implies $\mathcal{D}[\phi] = 0$ on $\overline{\Omega}$. In addition, by Lemma~\ref{lem: Complex layer}(1), $\mathcal{D}_f[\phi]$ is analytic on $\mathcal P_\tau \setminus \overline{D}$ with $\Re \mathcal{D}_f[\phi] = \mathcal{D}[\phi] = 0$. Thus, the Cauchy--Riemann equations imply $\mathcal{D}_f[\phi] \equiv iC$ on $\mathcal P_\tau \setminus \overline{D}$ for some $C \in \mathbb R$, yielding $\mathcal{D}_f[\phi] (z^{+}_0) \equiv iC$ for all $z_0\in \partial D$. Furthermore, Lemma~\ref{lem: Complex layer}(2) implies for all $z_0  \in \partial D$,
 $$       
 \mathcal{D}_f[\phi] (z^{-}_0) =  \frac{\phi(z_0)}{2} +\K_f[\phi] (z_0) = \frac{\phi(z_0)}{2} + \left[\mathcal{D}_f[\phi] (z^{+}_0) + \frac{\phi(z_0)}{2}\right] =  \phi(z_0) + iC.
$$
Since $\mathcal{D}_f[\phi]$ is analytic on each $D_j$ and its imaginary part is constant,  the real part $\phi(z_0)$ must be constant on each $\partial D_j$. i.e., \(\phi = \sum_{j=1}^M c_j \mathbbm{1}_j\) for some $c_j \in \mathbb R$. 

We show that if \(\phi\equiv c\), then $c=0$, implying that $\dim N(\K-\onehalf I ) \leq M-1$.  Let \(\phi \equiv c \in \mathbb{R}\). Since \(\phi \in N(\K - \onehalf I)\), it follows from Lemma~\ref{lem: exact double layer} that 
\begin{align*}
(\K-\onehalf I )[\phi](z_0)& = -c\sum_{j=1}^M \frac{|D_j|}{b} = 0 \quad \text{for} \quad z_0 \in \partial\Omega. 
\end{align*}
Since $\sum_{j=1}^M \frac{|D_j|}{b}\neq 0$ due to Assumption~\ref{assumption1}, it follows that  $c = 0 \equiv \phi$. 

Next, we show that $\dim N(\K-\onehalf I ) \geq M-1$ by verifying that $\{\phi_{jM}\}_{j=1}^{M-1}$in~\eqref{eqn: characteristic phi} forms a basis for $N(\K-\onehalf I )$. The set $\{\phi_{jM}\}_{j=1}^{M-1}$ is linearly independent. Furthermore, using
$\K [\mathbbm{1}_j]= \frac{\mathbbm{1}_j}{2}- \frac{|D_j|}{b}\mathbb{}$ from~\eqref{eqn: gen gauss sublemma},  we have
\begin{align*}
  (\K-\onehalf I )[\phi_{jM}]= \frac{b}{|D_j|}\K [\mathbbm{1}_j] - \frac{b}{|D_M|}\K [\mathbbm{1}_M] - \onehalf \phi_{jM} =\onehalf \phi_{jM}-\onehalf \phi_{jM} =0.
\end{align*}
Thus,  we conclude $\dim N(\K-\onehalf I )=M-1$ and $N(\K-\onehalf I ) = \textrm{span}\{\phi_{jM}\}_{j=1}^{M-1}$. 

\noindent {\bf (5)}
Since $\K$ is a compact linear operator by Lemma~\ref{lem: properties double layer}(3), the Fredholm alternative and Lemma~\ref{lem: properties double layer}(4) imply $\dim N(\K^*-\onehalf I ) =  M-1.$
Let $\{\psi_{jM}\}_{j=1}^{M-1}$ be a basis of this null space.  We show \eqref{eqn: Spsi span}. Indeed, by Lemma~\ref{lem: exact double layer}, we obtain 
\begin{align*}
    & 0 = \langle(\K^*-\onehalf I )\psi_j, \mathbbm{1}\rangle = \langle \psi_j, (\K-\onehalf I )[\mathbbm{1}]\rangle = \frac{-|D|}{b}\int_{\partial\Omega}\psi_j(\xi)\,|d\xi|,
\end{align*}
which implies $\int_{\partial\Omega}\psi_j(\xi)\,|d\xi| = 0$ for all $j$, given that $-|D|/b\neq0$. Additionally, Lemma~\ref{lem: properties single layer} (2) implies that $\partial_\nu\mathcal{S}[\psi_j](z^-_0)$ vanishes for all $j$ for all $z_0\in\partial\Omega$. Hence, Lemma \ref{lem: properties single layer} (1) implies that $\mathcal{S}[\psi_j](z)$ solves a homogeneous Neumann BVP on $D$; thus, \eqref{eqn: Spsi span A} holds. Equation \eqref{eqn: Spsi span B} then follows from the continuity of $S[\psi_j]$ on $\overline{D}$.

To show linear independence of $\{ \mathbf{s}_1,\mathbf{s}_2, \ldots, \mathbf{s}_{M-1}, \mathbf{1}\}$, suppose $\sum_{j=1}^{M-1}c_j \mathbf{s}_j+ c_M\mathbf{1} = 0$. We show that $c_j = 0$ for all $j =1,2,\ldots,M$. Let
\begin{align*}
     u(z) =\sum_{j=1}^{M-1} c_j\mathcal{S}[\psi_j](z)  \quad \text{for } z\in \mathbb{T}_\tau \,\,\implies \,\, u(z) = -c_M \quad \text{on } \overline{D} \,\,(\text{because of ~\eqref{eqn: Spsi span}}). 
\end{align*}
Note that $\Delta u(z) = 0$ on $\Omega$ and $u(z)$ is continuous throughout $\mathbb{T}_\tau$. Thus, $u(z^{+}_0) = -c_M$ implies 
$u(z) \equiv -c_M$ on $\mathbb{T}_\tau.$ 
Finally, apply Lemma~\ref{lem: properties single layer}(2)  for all $z_0\in \partial\Omega$ to obtain
\begin{equation*}
    0 = \partial_\nu u(z^{+}_0)- \partial_\nu u(z^{-}_0)= \sum_{j=1}^{M-1} c_j\left[\partial_\nu\mathcal{S}[\psi_j](z^+_0)-\partial_\nu\mathcal   {S}[\psi_j](z^-_0)\right]=\sum_{j=1}^{M-1}c_j\psi_j(z_0).
\end{equation*}
However, $\{\psi_j\}_{j=1}^{M-1}$ forms a basis of $N(\K^{*}-\onehalf I )$, implying $c_j =0$ for all $j = 1,2, \ldots, M$.

\noindent {\bf (6)} 
We show that $\dim N(\K^{*} + \onehalf I ) = 0$, which implies  $\dim N(\K + \onehalf I ) = 0$ from the Fredholm alternative. Let $\phi\in N(\K^* +\onehalf I )$.
By Lemma~\ref{lem: exact double layer}, 
\begin{align*}
    0 =\langle(\K^*+\onehalf I )[\phi],\mathbbm{1}\rangle= \langle\phi,(\K+\onehalf I )[\mathbbm{1}]\rangle =  \left(1-\frac{|D|}{b}\right) \int_{\partial\Omega}\phi(\xi)\,|d\xi|.
\end{align*}
Since $|D|<b$ (the area of the holes, $|D| = \sum_{j=1}^{M} |D_j|$, is less than the area of the torus, $|\mathbb{T}_\tau| = b$), we have $\int_{\partial\Omega}\phi(\xi)\,|d\xi| = 0$. By Lemma~\ref{lem: properties single layer}(1)--(2), $\mathcal{S}[\phi]$ solves the homogeneous Neumann BVP on $\Omega$, implying $\mathcal{S}[\phi] \equiv C$ on $\overline{\Omega}$, for some $C\in \mathbb{R}$.  In addition, by continuity (Lemma~\ref{lem: properties single layer}(1)) and the maximum principle, $\mathcal{S}[\phi] \equiv C$ on each  $\overline{D}_j$ for $j = 1,2, \ldots, M$. Finally, the jump relation (Lemma~\ref{lem: properties single layer}(2)) yields
\begin{equation*}
    0 = \partial_\nu \mathcal{S}[\phi](z^+_0)- \partial_\nu \mathcal{S}[\phi](z^-_0) =\phi(z_0), \quad \forall z_0 \in \partial\Omega.
\end{equation*}
Thus, we conclude that $\dim N(\K^* + \onehalf I ) = 0$. 
\end{proof}

\subsection{Proof of Lemma~\ref{lem: properties single layer}}
\begin{proof}
We follow proofs for the analogous statements for single-layer potentials on Euclidean domains, as found in \cite{folland2020introduction,kress1989linear,kress1991boundary,Nasser_2011}. \\ 
\noindent{\bf (1)}  The operator \(\mathcal{S}\) is linear and inherits doubly-periodic property from its kernel \(G(z - \xi)\)~\cite{lin2010elliptic}. For $\phi\in C_0(\partial\Omega)$, a direct calculation yields
\begin{equation*}
        \Delta \mathcal{S}[\phi](z)  = \int_{\partial \Omega} \phi(\xi) \Delta G(z-\xi)\,|d\xi|= \frac{1}{b} \int_{\partial \Omega} \phi(\xi) \,|d\xi| = 0, \quad z\in \mathbb{T}_\tau\setminus\partial\Omega. 
\end{equation*} 

To establish continuity on $\mathbb{T}_\tau$, it is sufficient to show continuity at $\partial \Omega$, following the argument in~\cite[Prop. 3.25]{folland2020introduction}. For $z_0\in \partial\Omega$ and $\varepsilon>0$, fix $\gamma >0$. For $z\in B_\gamma(z_0)$,
\begin{align*}
    |\mathcal{S}[\phi](z)-\S[\phi](z_0)| \leq&  \Vert\phi\Vert_{L^\infty(\partial\Omega)} \int_{\partial\Omega \setminus B_{\gamma}(z_0)} |G(z_0-\xi)-G(z-\xi)|\, |d\xi|\\
&+\Vert\phi\Vert_{L^\infty(\partial\Omega)}\int_{\partial\Omega\cap B_{\gamma}(z_0)}|G(z-\xi) -  G(z_0-\xi)|\, |d\xi|.
\end{align*}
The first integral can be bounded by $\varepsilon/(2\Vert\phi\Vert_{L^\infty(\partial\Omega)})$ for $z$ close to $z_0$ due to the uniform continuity of $G$ away from the singularity. For the second term, letting $S=\partial\Omega\cap B_{\gamma}(z_0)$ and using~\eqref{e: green approx}, we simplify
\begin{align*}
    \int_{S}|G(z-\xi) -  G(z_0-\xi)|\, |d\xi|& \leq \int_{S} \frac{\log|z-\xi|^{-1} + \log|z_0-\xi|^{-1}}{2\pi}\, |d\xi|+O(|z-z_0|^2).
\end{align*}
As $z \to z_0$, $O(|z-z_0|^2)$ can be made less than $\varepsilon/(4\Vert\phi\Vert_{L^\infty(\partial\Omega)})$. Evaluated in polar coordinates, the log-integrals are $O(-\gamma \log \gamma)$, which can be made less than $\varepsilon/(4\|\phi\|_{L^\infty(\partial\Omega)})$ given that $\gamma$ is fixed small enough. Combined, these estimates ensure $|\mathcal{S}[\phi](z) - \S[\phi](z_0)| \leq \varepsilon$, establishing continuity at the boundary.

\noindent {\bf (2)}  We prove the result for \( z_0^{+} \), following ~\cite[Thm. 3.28]{folland2020introduction}; the case for \( z_0^{-} \) follows similarly. Let $z_0 \in \partial D_j$ and fix \( \gamma > 0 \) small so that for $z\in B_\gamma(z_0)\cap \Omega$, we may express \( z = z_0 - h\nu(z_0) \). Note that \( \K^{*}[\phi](z) \) is continuous on \( \partial \Omega \) due to Lemma~\ref{lem: properties double layer}(3), so it suffices to show that
\begin{align*}
&\lim_{h \to 0} \int_{\partial \Omega} \phi(\xi) \left[ J(z, \xi) - J(z_0, \xi) \right]\, |d\xi| = 0,\,\, \text{for}\,\,  J(z, \xi) := \partial_{\nu_z} G(z - \xi) + \partial_{\nu_\xi} G(z - \xi). 
\end{align*}
Here,  we evaluate the above integral at $z_0$ in the sense of $\mathrm{p.v.}$ 
Applying Lemma~\ref{lem: properties double layer}(2) leads to
\begin{align*}
        \partial_\nu \mathcal{S}[\phi](z_0^+) &= \lim_{h \to 0} \int_{\partial \Omega} \phi(\xi) \left[ J(z_0, \xi) - \partial_{\nu_\xi} G(z - \xi) \right] \, |d\xi| = \onehalf\phi(z_0) + \K^{*}[\phi](z_0).
\end{align*}

\noindent\textit{Proof of the claim:} 
Let $\varepsilon > 0$. We split the integral into $S = \partial\Omega \setminus B_\gamma(z_0)$ and $\tilde{S} = \partial D_j \cap B_\gamma(z_0)$ and show that $|(\text{I})|, |\text{(II)}\leq \varepsilon/2$:
\begin{align*}
   \text{(I)}+ \text{(II)} = \int_{S} \phi(\xi) \left[ J(z, \xi) - J(z_0, \xi) \right]  \, |d\xi| +  \int_{\tilde{S}} \phi(\xi) \left[ J(z, \xi) - J(z_0, \xi) \right] \, |d\xi|. 
\end{align*}
To bound $|(I)|$, since $J$ is uniformly continuous away from the singularity, choosing $z$ sufficiently close to $z_0$ ensures \(\Vert J(z, \xi) - J(z_0, \xi) \Vert_{L^\infty(S)} \leq \varepsilon/(2\Vert \phi \Vert_{L^{\infty(\partial\Omega)}}|\partial\Omega|)\). 

Next, to bound \( |\text{(II)}|\), we use approximations \eqref{eqn: dgdxi close} and \eqref{eqn: dgdz close}  to obtain
\[J(z, \xi)= - \frac{(z - \xi) \circ [\nu(z_0) - \nu(\xi)]}{2\pi |z - \xi|^2} + O(|z - \xi|).\]
The same bound is derived for $z_0$ by replacing $z$. Therefore, it follows that
\begin{align*}
      &\left\vert \int_{\tilde{S}} \phi(\xi) \left[ J(z, \xi) - J(z_0, \xi) \right] \, |d\xi| \right\vert \\
      &\leq  \frac{\Vert\phi \Vert_{L^{\infty}(\partial \Omega)}}{2\pi} \int_{\tilde{S}} \left[ \left\vert\frac{(z - \xi) \circ [\nu(z_0) - \nu(\xi)]}{|z - \xi|^2}\right\vert + \left\vert \frac{(z_0 - \xi) \circ [\nu(z_0) - \nu(\xi)]}{|z_0 - \xi|^2}\right\vert \, |d\xi|\right]\\
      &+ O(|z-z_0|),
\end{align*} where $O(|z-z_0|)$ can be made less than $\varepsilon/4$ by letting $z$ sufficiently close to $z_0$.  Since $z$ is normal through $z_0$, $|z-\xi|\geq \frac{1}{C} |z_0-\xi|$ for some constant.
By \eqref{lemma: C2 prop}, we have
\begin{align*}
  \left\vert \frac{(z-\xi) 
 \circ [\nu(z_0)-\nu(\xi)]}{|z - \xi|^2} \right\vert &\leq   \frac{ 
 |\nu(z_0)-\nu(\xi)|}{|z - \xi|} \leq CL_{2,j}\\
  \left\vert \frac{(z_0-\xi) 
 \circ [\nu(z_0)-\nu(\xi)]}{|z_0 - \xi|^2} \right\vert &\leq  \frac{|\nu(z_0)-\nu(\xi)|}{|z_0-\xi|}\leq L_{2,j},
\end{align*}
Hence, choosing a sufficiently small $\gamma$ such that $|\tilde{S}| \leq 2\pi\varepsilon/(4\Vert \phi \Vert_{L^{\infty}(\partial \Omega)}(1+C)L_{2,j})$ allows $\vert\text{(II)}\vert\leq\varepsilon/2$

\noindent {\bf (3)}  We show that \( \S \)  is an integral operator of order zero, implying its compactness \cite[Eq. (3.10) and Props. 3.10]{folland2020introduction}. It is evident that \(G(z - \xi)\) remains continuous for \(z \neq \xi\). While $G(z - \xi)$ is continuous for $z \neq \xi$, for $z$ near $\xi \in \partial\Omega$, \eqref{e: green approx} implies:
\begin{equation*}
G(z-\xi) = -\frac{1}{2\pi} \log|z - \xi| + c + O(|z - \xi|^2) = A(z,\xi)\log|z-\xi| + B(z,\xi),
\end{equation*}
where \(A(z, \xi) = -\frac{1}{2\pi}\) and \(B(z, \xi)\) is a bounded remainder containing the constant and higher-order terms.
\end{proof}

\subsection{Proof of Lemma~\ref{lem: characteristic}}
\begin{proof} 
\noindent{\bf (1)} We see that $\X$ is linear and compact, as it has finite rank. By definition, $\X$ is a constant function on each component $\partial D_j$: 
\begin{equation}
    \X[\phi](z_0) =\begin{cases} \int_{\partial D_j}\phi(\xi)\,|d\xi| & \text{if } z_0 \in \partial D_j \text{ for } j  =1,2,\ldots,M-1,\\
    0 & \text{if } z_0 \in \partial D_M.\end{cases}
\end{equation}

\noindent {\bf (2)} The compactness and linearity of $\K + \X - \onehalf I $ follows from Lemma~\ref{lem: characteristic}(1) and Lemma~\ref{lem: properties double layer}(3). Therefore, the Fredholm alternative implies that injectivity ensures the surjectivity. To establish injectivity, let $\phi\in N(\K + \X - \onehalf I )$; we show that $\phi\equiv 0$. Substituting $\phi$ to $\K + \X - \onehalf I$, we obtain
\begin{equation}
\label{e:dirich home M-1}
   \left(\K - \onehalf\right)[\phi] = \sum_{j=1}^{M-1} c_j\mathbbm{1}_j, \quad\text{ with } c_j = 
        - \int_{\partial D_j} \phi(\xi) \, |d\xi|.
\end{equation}
We make the following claim:

\medskip 
\noindent \textit{Claim}: $c_j = 0,\, \forall j=1,2,\ldots, M-1$. 
\medskip

Postponing the proof of the claim, this implies $\phi\in N(\K - I/2)$, allowing us to express $\phi = \sum_{j=1}^{M-1} d_j \phi_{jM}$ using the basis in Lemma~\ref{lem: properties double layer}(4), \eqref{eqn: characteristic phi}. In particular, since $c_j = 0$ for all $j= 1,2,\ldots,M-1$ and~\eqref{e:dirich home M-1}, we obtain
\[
0 = c_j = -\int_{\partial D_j}\phi(z_0) |dz_0|= -d_j\frac{b|\partial D_j|}{ |D_j|} \implies d_j = 0.
\]
Thus, $\phi \equiv  0$, confirming injectivity. 

\medskip
\noindent \textit{Proof of the Claim}: By Lemma~\ref{lem: Complex layer}(1), $\mathcal{D}_f[\phi]$ is analytic on $\mathcal{P} \setminus \partial\Omega$, allowing us to apply the Cauchy-Riemann equations: 
\begin{equation}
-\partial_\nu \Im \mathcal{D}_f[\phi](z^+_0) = \partial_\sigma \Re \mathcal{D}_f[\phi](z^+_0)= \partial_\sigma\mathcal{D}[\phi](z^+_0)= 0, \quad \forall z_0 \in \partial \Omega.
\label{eqn: zero partial nu dirch}
\end{equation}
The tangential derivative $\partial_\sigma \mathcal{D}[\phi](z^+_0)$ vanishes since $\mathcal{D}[\phi](z^+_0)$ is piecewise constant on $\partial D_j$ according to~\eqref{e:dirich home M-1}. In addition, we know that $\mathcal{D}_f[\phi](z)$ is analytic, so $\Im \mathcal{D}_f[\phi](z)$ is harmonic. The Green's formula over the domain $\mathcal{P} \setminus D$ yields 
\begin{align}
\int_{\mathcal{P} \setminus D} |\nabla \Im \mathcal{D}_f[\phi](z)|^2 \, dA = \int_{\partial \mathcal{P}} \Im \mathcal{D}_f[\phi](z) \partial_\nu \Im \mathcal{D}_f[\phi](z) \, |dz|,
\label{eqn: constant Im f}
\end{align}
where the boundary integrals vanish due to \eqref{eqn: zero partial nu dirch}.

\begin{figure}[t!]
\begin{center}
\includegraphics[width=0.4\textwidth]{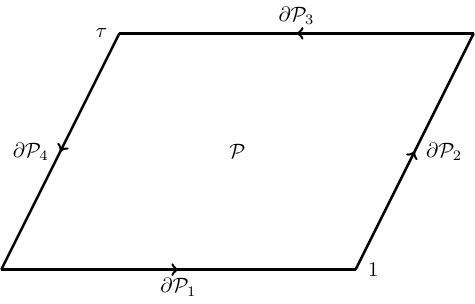}
\end{center}
\caption{An illustration of notation introduced in the proof of Lemma~\ref{lem: characteristic}(2) for the parallelogram $\mathcal{P}$ defined in~\eqref{e:parallelo}.}
\label{fig:parallelo}
\end{figure}

We show the right-hand side of~\eqref{eqn: constant Im f} vanishes by splitting and calculating the integral into four parts
$\partial P = \bigcup_{i=1}^4 \partial P_i$ (see Fig. \ref{fig:parallelo}). Although $\Im \mathcal{D}_f[\phi](z)$ is not periodic, its normal derivative $\partial_\nu \Im \mathcal{D}_f[\phi](z)$ is doubly periodic due to $\partial_\nu \Im \mathcal{D}_f[\phi](z) = -\partial_\sigma \mathcal{D}[\phi](z)$. Therefore, we simplify the right hand side of~\eqref{eqn: constant Im f} as
\begin{align*}
    \int_{\partial P_1\cup\partial P_3 } 
    &= -\int_{0}^1\partial_\sigma \Re \mathcal{D}_f[\phi](t)\Im \mathcal{D}_f[\phi](t)  \,dt +\int_{0}^1\partial_\sigma \Re \mathcal{D}_f[\phi](t+\tau)\Im \mathcal{D}_f[\phi](t+\tau) \,dt\\
    &=\int_{0}^1\partial_\sigma \Re \mathcal{D}_f[\phi](t)\left[\Im \mathcal{D}_f[\phi](t+\tau)-\Im \mathcal{D}_f[\phi](t)\right]  \,dt =C_2\int_{0}^1\partial_\sigma \Re \mathcal{D}_f[\phi](t)  \,dt\\
    &= 0,
\end{align*}
where $C_2$ is a constant and $\Re \mathcal{D}_f[\phi](1) =\Re \mathcal{D}_f[\phi](0)$. A similar calculation shows that $\int_{\partial P_2\cup\partial P_4} = 0$. Therefore, we conclude that the right-hand side of~\eqref{eqn: constant Im f} is $0$ so that $\Im \mathcal{D}_f[\phi]$ is constant on $\mathcal{P}\setminus D$. 

Finally, since $\mathcal{D}_f[\phi](z)$ is analytic on $\mathcal{P}\setminus\partial\Omega$, the Cauchy-Riemann equations imply that $\Re \mathcal{D}_f[\phi](z) =\mathcal{D}[\phi](z)$ is constant on $\Omega$. By~\eqref{e:dirich home M-1}, $\mathcal{D}[\phi](z^+_0) =\K[\phi](z_0) - \frac{\phi(z_0)}{2} =  c_M = 0$ on $\partial D_M$, which implies that $c_1 =\cdots= c_{M-1} =0$, establishing the claim. 
\end{proof}

\end{document}